\documentclass{amsart}
\usepackage{amssymb}
\usepackage{amscd}
\usepackage{amsmath}
\usepackage{enumerate}
\usepackage[dvips]{graphicx}
\newtheorem{theo}{Theorem}
\newtheorem{lema}[theo]{Lemma}
\newtheorem{cor}[theo]{Corollary}
\newtheorem{prop}[theo]{Proposition}

\newtheorem{definition}[theo]{Definition}

\newtheorem{remark}[theo]{Remark}

\newtheorem{example}[theo]{Example}
\newtheorem{problem}[theo]{Problem}

\newcommand{\CC}{{\mathbb{C}}}
\newcommand{\NN}{{\mathbb{N}}}
\newcommand{\PP}{{\mathbb{P}}}
\newcommand{\RR}{{\mathbb{R}}}
\newcommand{\SSS}{{\mathbb{S}}}
\newcommand{\ZZ}{{\mathbb{Z}}}
\newcommand{\Zb}{{\mathbf{Z}}}
\newcommand{\Zc}{{\mathbf{Z'}}}

\newcommand{\calO}{{\mathcal{O}}}

\newcommand{\calT}{{\mathcal{T}}}

\newcommand{\calV}{{\mathcal{V}}}
\newcommand{\calW}{{\mathcal{W}}}
\newcommand{\calX}{{\mathcal{X}}}
\newcommand{\calY}{{\mathcal{Y}}}
\newcommand{\calZ}{{\mathcal{Z}}}
\newcommand{\comp}{{\circ}}

\newcommand{\mm}{{\mathbf{m}}}

\newcommand{\mult}{{\mathrm{mult}}}

\newcommand{\Id}{{\mathrm{Id}}}

\begin{document}
\title[Topological equisingularity]{Topological equisingularity of function germs with $1$-dimensional critical set}
\author{Javier Fern\'andez de Bobadilla}
\address{Departamento de Matem\'aticas Fundametales. Facultad de Ciencias. U.N.E.D. c/ Senda del Rey 9.
28040 Madrid. Spain.}
\email{javier@mat.uned.es}
\thanks{Supported by Ramon y Cajal contract. Supported by the Spanish project MTM2004-08080-C02-01.}
\date{}
\begin{abstract}
We focus on topological equisingularity of families of holomorphic function germs with 1-dimensional critical set. We introduce the notion of equisingularity at the critical set and prove that any family which is equisingular at the critical set is topologically equisingular. We show that if a family of germs with 1-dimensional critical set has constant generic Le numbers then it is equisingular at the critical set, and hence topologically equisingular (answering a question of D. Massey). We use this to modify the definition of singularity stem present in the literature, introducing and characterising topological stems (being this concept closely related with Arnold's series of singularities). We provide another sufficient condition for topological equisingularity for families whose reduced critical set is deformed flatly. Finally we study how the critical set can be deformed in a topologically equisingular family and provide examples of topologically equisingular families whose critical set is a non-flat deformation with singular special fibre and smooth generic fibre.
\end{abstract}


\maketitle

\section{Introduction}

In this paper we answer several topological equisingularity questions for holomorphic function germs with $1$-dimensional critical set.

The topological equisingularity notions that we will handle are the following.
Denote by $B_\epsilon$ and $\SSS_\epsilon$ the ball and the sphere of radius $\epsilon$ centered at the origin of $\CC^n$.
Let $f:(\CC^n,O)\to\CC$ be a holomorphic germ. The {\em embedded link} of $f$ is the topological type of the pair
$(\SSS_\epsilon,f^{-1}(0)\cap\SSS_\epsilon)$ for $\epsilon$ sufficiently small; the {\em abstract link} is the topological type
of $f^{-1}(0)\cap\SSS_\epsilon$. Two germs $f$ and $g$ are topologically $R$-equivalent if there is a germ of self-homeomorphism 
$\varphi$ of $(\CC^n,O)$ such that $g=f\comp\varphi$; if we only have that $\varphi(f^{-1}(0))=g^{-1}(0)$ then we say that 
$f$ and $g$ have the same topological type. By the conical structure of singularities having the same embedded link implies 
having the same topological type. A family $f_t:(\CC^n,O)\to\CC$ of holomorphic germs depending
continuously on a parameter $t$ varying in a manifold $T$ is {\em topologically $R$-equisingular} if there is
a family of self-homeomorphisms $\varphi_t$ depending continuously in $t$ and parametrised over $T$ such that 
$f_{t}\comp\varphi_{t}=f_{t_0}$ for a fixed $t_0\in T$ and any $t\in T$.

Topological equisingularity has proved to be a very subtle subject, with several long standing open questions. For instance 
Zariski's Multiplicity Question asks whether two germs with the same topological type must have the same multiplicity 
(Question~A,~\cite{Zar}). 

Suppose that $f_t$ has an isolated singularity at the origin and that the Milnor number of $f_t$ is independent on $t$.
Assume $n\neq 3$. 
L\^e and Ramanujam~\cite{LR} proved that the diffeomorphism type of the Milnor fibration of $f$ and of the embedded link are
independent on $t$. Later King~\cite{Ki2} and Timourian~\cite{Ti} proved that for any $t\in T$ there is a neighborhood $U$ of $t$ 
in $T$ such that the restriction of $f_t$ over $U$ is topologically $R$-equisingular. As the Milnor number is a topological
invariant, if $n\neq 3$, it is an invariant
characterising topological equisingularity for germs with isolated singularities. Answering whether families with constant Milnor
number are topologically equisingular for $n=3$ is a major open problem.

Topological equisingularity for germs having not necessarily isolated singularities turns out to be much more
difficult, and the theory is less developed than in the case of germs with isolated singularities.
As an illustration of this, and also to motivate the results of this paper, let us mention some of the latest developments.
In~\cite{Ma1},~\cite{Ma2} Massey introduced the L\^e numbers,
a set of polar invariants attached to a germ and to a coordinate system. Suppose that the family $f_t$ depends holomorphically on a 
complex parameter $t$ and that the critical set of $f_t$ has dimension $s$ (with $n\geq s+4$) for any $t\in T$,
Massey~\cite{Ma2} proved that if the L\^e numbers of $f_t$ with respect to a 
coordinate system (satisfying a certain genericity property) are constant then the diffeomorphism type of the Milnor fibration 
of $f_t$ is independent on $t$ (it was proved by the author in~\cite{Bo2} 
that the homotopy type of the abstract links remains constant if
$n\geq 4$, without the condition that $n\geq s+4$). 
It remained open whether the embedded topological type of $f_t$ is independent on $t$. This was 
answered in the negative in~\cite{Bo2}, where counterexamples 
with critical set of dimension $3$ were provided (a counterexample of Zariski's
Question~B, of~\cite{Zar} was constructed as an application). Probably, the main reason of the failure of the L\^e numbers in 
controlling the embedded topological type is the fact that they control primarily the Milnor fibration, but, as it is shown 
in~\cite{Bo2} in the case of non-isolated singularities the relation of the Milnor fibration with the embedded topological type 
is, in general, weak (in the case of isolated singularities the Milnor fibration determines the embedded topological type).
On the other hand, in~\cite{BG}, it was shown that the generic L\^e numbers are not topological invariants (even in the case of critical set of dimension $1$). A different approach for topological equisingularity of families of functions satisfying certain conditions on
their Newton polytope can be found in~\cite{DG}.

In this paper we adopt a new viewpoint towards topological equisingularity. Although a general 
formulation of it is possible (and the reader is encouraged to figure it out) we introduce it for the case of functions 
with $1$-dimensional critical set (which are the object of interest along this paper). Functions with $1$-dimensional critical set occupy a distinguished position among those with non-isolated singularities, since they are at the limit of series of singularities; we will be more explicit about this below.  
The critical set of such a function $f$ (the support of the 
sheaf of vanishing cycles) is a curve which is stratified by the Milnor number of the singularity of the generic transverse
hyperplane sections. Roughly speaking, a family $f_t$ is {\em equisingular at the critical set}
if there is a family of self-homemorphisms
of $(\CC^n,O)$ topologically trivialising the family of critical sets and respecting the stratification by transverse Milnor 
number {\em outside the origin} (see~Definition~\ref{equiatsing}). 
We may view equisingularity at the critical set as an intermediate notion between 
topological equisingularity itself and the constancy of any candidate set of invariants intended to imply topological equisingularity.
The strategy is to prove that equisingularity at the critical set implies (any notion of) topological equisingularity and to 
find numerical invariants implying equisingularity at the critical set.

Now we summarise the main results of this paper. In the families that one finds in practice, the dependence on the parameter is 
always at least differentiable. We will work in this paper with families $f_t$ of holomorphic germs with $1$-dimensional critical 
set depending smoothly on a
parameter $t$ varying on a manifold $T$. We will gain, in case that the family is topologically equisingular, some additional 
regularity on the trivialising family of homeomorphisms. Often in this paper the space of paremeters $T$ is assumed to be 
(diffeomorphic to) a cube $[0,1]^n$. This is due to the fact that we prove topological equisingularity theorems which are {\em global
over the parameter space}, and in order to achieve it is necessary that the parameter space has trivial topology. Anyway, if the
objetive was to obtain equisingularity results {\em local over the parameter space}, the assumption that it is a cube 
does not decrease generality.\\

\noindent\textbf{Theorem A}. (See Theorems~\ref{indeplinkmilnorfib}~and~\ref{porfin} for more details.)
{\em If $f_t$ is equisingular at the critical set then
\begin{enumerate}[(i)]
\item The diffeomorphism type of the Milnor fibration of $f_t$ is independent on $t$.
\item The topological type of the embedded link and the topological right-equivalence type of $f_t$ are independent on $t$.
\item If $T$ is a cube then $f_t$ is topologically $R$-equisingular over $T$. The family of homeomorphisms $\varphi_t$ trivialising
$f_t$ satisfies that the restrictions of $\varphi_t$ to the complement of the critical set, and to 
the critical set minus the origin
is a diffeomorphism, and the dependence on $t$ of $\varphi_t$ restricted to these strata is differentiable. 
\end{enumerate}
$\quad$}\\

\noindent\textbf{Theorem B}. (See Theorem~\ref{leconsteqsing}.)
{\em If $f_t$ depends holomorphically on a parameter and any of the following conditions hold:
\begin{enumerate}
\item there is a coordinate system $\mathbf{Z}$ such that the L\^e numbers at the origin of $f_t$ 
with respect to $\mathbf{Z}$ are defined and independent of $t$, 
\item the generic L\^e numbers at the origin of $f_t$ are independent of $t$,
\end{enumerate}
then $f_t$ is equisingular at the critical set.}\\

We immediately deduce\\

\noindent\textbf{Corollary C}. (See Theorem~\ref{wanted}.)
{\em Suppose that $f_t$ depends holomorphically on a parameter, if one 
of the following conditions hold 
\begin{enumerate}
\item there is a coordinate system $\mathbf{Z}$ such that the L\^e numbers at the origin of $f_t$ 
with respect to $\mathbf{Z}$ are defined and independent of $t$,
\item the generic L\^e numbers at the origin of $f_t$ are independent of $t$, 
\end{enumerate}
then the embedded topological type of $f_t$ at the origin and the the diffeomorphism type of the 
Milnor fibration is independent of $t\in T$. Moreover, if $T$ is (diffeomorphic to) a cube, 
the restriction of the family over $T$ is topologically $R$-equisingular}.\\

This recovers Massey's Theorem~\cite{Ma2} for $1$-dimensional critical set and answers positively Massey's question about the 
embedded topological type, in contrast with what happens for higher dimensional critical set~\cite{Bo2}.

Isolated hypersurface singularities are related to non-isolated ones through the idea of ``series 
of singularities''. In his classification of low codimensional $\mu$-classes of isolated singularites
\cite{Ar}, \cite{AGV}, Arnol'd observed that the list of singularities splits into series. Although
it was impossible to give a precise definition of what a series of singularities is, it became clear that
series are associated with singularities of infinite codimension (non-isolated singularities). Inspired 
by Arnol'd remark C.T.C Wall formulated a conjecture 
(see the Conjecture at the introduction of~\cite{Wa}), or rather, a guiding principle for 
classification of singularities. Later D. Mond~\cite{Mo} introduced the idea of singularity stem, as a 
first step to understand the notion of series of singularities in the context of mappings from $\CC^2$ to
$\CC^3$ (the idea is that a singularity stem is a
non-isolated singularity which lies in the limit of a series of singularities). 
Following a suggestion of Montaldi, Pellikaan~\cite{Pe} defined inductively the notion of stem of degree 
$d$, characterised stems of degree $1$ as functions with irreducible $1$-dimensional critical set and 
transversal type $A_1$, and proved that any stem of finite degree is a function with $1$-dimensional 
critical set. He was also able to give bounds on the degree of the stem depending on the number of 
irreducible components of the critical set of the function and the transversal Milnor number. 
Here we modify the concept of stem and define topological stem (see Definition~\ref{topstem}). 
With this notion we prove\\

\noindent\textbf{Theorem D}. (See Theorem~\ref{charstem}.)
{\em A holomorphic function germ $f:(\CC^n,O)\to\CC$ is a topological stem of positive finite degree if and only if 
its critical set is $1$-dimesional at the origin. Moreover the degree of the stem is bounded above by 
the generic first L\^e number at the origin of $f$.}\\

Our modification of the definition of stem consists essentially in replacing differentiable 
$R$-equivalence by topological $R$-equivalence. 
This is natural since series in Arnol'd classification of singularities are 
in fact topological series (series of $\mu$-classes), and thus it is therefore reasonable that the object
to find at the limit of the series admits a topological definition as well. It is interesting to notice 
that Theorem~D is an application of our sufficient conditions for topological $R$-equisingularity, which
is not at reach using any of the previously known sufficient conditions for equisingularity.
 
A consequence of our theory and of the results of~\cite{BuG} is the following theorem, giving an 
alternative sufficient condition of topological equisingularity, for the case
that the reduced critical set is deformed flatly:\\

\noindent\textbf{Theorem E}. (See Theorem~\ref{gbb}.)
{\em Let $f:\CC\times(\CC^n,O)\to (\CC,0)$ be a family of holomorphic germs at the origin, holomorphically depending on a parameter $t$, 
and having $1$-dimensional critical set at the origin. Let $\Sigma:=V(\partial f/\partial z_1,...,\partial f/\partial z_n)$.
Suppose that all the irreducible components of $\Sigma$ at $(0,O)$ are of dimension $2$ and that the restriction 
\[\pi:\Sigma_{red}\to \CC\]
is a flat morphism with {\em reduced} fibres (where $\Sigma_{red}$ is $\Sigma$ with reduced structure and $\pi$ is the restriction 
of the projection of $\CC\times\CC^n$ to the first factor) such that the Milnor number at the origin $\mu((\Sigma_{red})_t,O)$
(in the sense of~\cite{BuG}) is independent of $t$. If, in a neighborhood of $(0,O)$ in $\CC\times\CC^n$,
the generic transversal Milnor number of $f_t$ at any point $(t,x)$ of
$\Sigma_{red}\setminus\CC\times\{O\}$ only depends on the irreducible component of $\Sigma_{red}$ to which $(t,x)$ belongs 
then the family $f$ is topologically $R$-equisingular.}\\

Finally we study how the reduced critical set may vary in a family which is topologically $R$-equisingular. 
It seems that there are
no further restrictions appart from the obvious topological ones. 
We observe a great amount of flexibility in the critical set of a topologically equisingular family.
We exhibit a family $f_t$ (see Example~\ref{example2}) which is topologically
$R$-equisingular and the reduced critical set is singular for $t=0$ and smooth for $t\neq 0$ 
(thus the deformation of reduced critical sets is non flat). To the author's knowledge it is
the first time that this phenomenon is observed. 
Example~\ref{example2} is a family for which the sufficient conditions for topological equisingularity
given in Theorem~E fail: we really need the finer Theorem A to establish the equisingularity. In
Problem~\ref{problem} we propose to construct more families which are equisingular at the critical set
starting out of deformations of parametrisations (this is the way how Example~\ref{example2} was constructed). The motivation for it 
is based on the following observation: the change at multiplicity of the critical set of $f_t$ 
implies, using L\^e-Iomdine formulas (see~\cite{Ma3}), that adding a high power of a linear function to $f_t$ will not yield a 
$\mu$-contant family. This suggests that, in the realm of function germs with $1$-dimensional critical set,
 we have more space to find a possible counterexample of Zariski's multiplicity
question than in the realm of isolated singularities.

For the proof of Theorem A we need a improved version of the Theorems of King and Timourian as follows:\\

\noindent\textbf{Preliminary Theorem}. {\em Let $f_t$ be a family of germs with isolated singularities at the origin, with constant
Milnor number, and depending smoothly on a parameter $t$ which varies on a cube $T$. Then the family is 
topologically $R$-equisingular over $T$, and the restriction of the trivialising family of homeomorphisms $\varphi_t$,
to $\CC^n\setminus\{O\}$ is a family of diffeomorphisms depending smoothly on the parameter $t$.}\\

The methods appearing in~\cite{Ki1},~\cite{Ki2}~and~\cite{Ti} are certainly local in the base. 
We need to introduce a new device, which we call a {\em cut}, which allow us to prove the result above,
which is global in the base, and also
enables to prove the differentiability outside the origin. Moreover the method with which we prove the Preliminary Theorem
provides tools that generalise well to prove Theorem~A.

A natural question is to determine the relation between equisingularity at the critical set and topological equisingularity.
We conjecture that they coincide:\\

\noindent\textbf{Conjecture A}. {\em A family $f_t$ of function germs with $1$-dimensional critical set depending smoothly on $t$
is equisingular at the critical set if and only if $f_t$ is topologically equisingular to $f_{t'}$, for any parameters $t,t'$}.\\

It is easy to reduce the above conjecture to the following one:\\

\noindent\textbf{Conjecture B}. {\em Let $f$ be a function germ with smooth $1$-dimensional critical set such that the
transversal Milnor number of $f$ at any point of its critical set is independent on the point. Then any 
function germ $g$ topologically equisingular to $f$ has smooth critical set and constant transversal Milnor number}.\\

Actually I expect the following stronger statement to be true:\\

\noindent\textbf{Conjecture C}. {\em Let $g$ be a function germ with irreducible $1$-dimensional critical set. Let $\mu$ denote the 
transversal Milnor number of $f$ at a generic point of its critical set. If the reduced Euler characteristic
of the Milnor fibre of $f$ is equal to $(-1)^n\mu$, then $f$ is a function germ with smooth $1$-dimensional critical set such that the
transversal Milnor number of $f$ at any point of its critical set is equal to $\mu$}.\\

In fact, by a recent result of L\^e-Massey~\cite{LM} and Tibar~\cite{Tib}, in conjectures $B$ and $C$ it suffices to prove that 
the critical set of $g$ is smooth. Conjecture $C$ is a numerical characterisation of function germs with $1$-dimensional critical set 
which are topologically a product.

The material is distributed as follows.

In Section~2 we introduce cuts and prove the Preliminary Theorem. In Section~3 we define a particular kind of neighborhoods of the
origin in $\CC^n$ which are essential in the proof of Theorem A. In Section~4 we introduce equisingularity at the critical set.
In Section~5 we prove parts $(i)$ and $(ii)$ of Theorem~A. In Section~6 we prove part $(iii)$ of 
Theorem~A. In Section~7 we prove Theorem B. In Section~8 we define topological stems and prove Theorem~D.
Finally, in Section~9 we prove Theorem~E, study how the reduce critical set deforms in a topologically equisingular 
family, and construct Example~\ref{example2}.

I would like to thank M. Pe Pereira and T. Gaffney for their careful reading and useful corrections and suggestions. This paper is dedicated to Mar\'ia Jos\'e Hern\'andez Navarro and to 
Juan Jos\'e Olaz\'abal Valverde.

\section{The case of isolated singularities revisited}
\label{globalisationLR}

\subsection{Cobordisms}
A subspace $X\subset \RR^N$ is {\em smoothly stratified} if it admits an stratification whose strata are smooth submanifolds.
A stratified subspace $Y\subset X$ is a subspace which is a union of strata.
Given any subset $X$ of $\RR^N$ we denote by $\dot{X}$ the subset of interior points, and by $\partial X$
its frontier $X\setminus\dot{X}$. Observe that if $X$ is a stratified space, then both $\dot{X}$ and $\partial{X}$ are forcedly
stratified subspaces. A {\em stratified cobordism between $X_0$ and $X_1$} is, 
by definition, a triple $(X,X_0,X_1)$ given by a stratified space $X$ and 
two stratified subspaces $X_0$, $X_1$ of $\partial X$, satisfying the following condition: for any stratum $A\subset X$ the 
''cobordism closure'' 
\[\hat{A}:=A\cup(\partial A\cap (X_1\cup X_2))\] has a structure of manifold with boundary (being the boundary $\partial A\cap (X_1\cup X_2)$).
A stratified cobordism is {\em homotopically (homologically) trivial} if the 
inclusions of $X_0$ and $X_1$ in $X$ are homotopy (homology) equivalences. Given a stratified space $Y$, the product 
$Y\times [0,1]$ inherits a natural stratification which makes the triple $(Y\times [0,1],Y\times\{0\},Y\times\{1\})$ a stratified
cobordism. A stratified cobordisms $(Y,Y_0,Y_1)$ is {\em trivial} if there is a homeomorphism
\[\varphi:(X,X_0,X_1)\to (Y\times [0,1],Y\times\{0\},Y\times\{1\})\]
such that for any stratum $A\subset X$ there is a stratum $B\subset Y$ such that 
\[\varphi|_{\hat{A}}:\hat{A}\to B\times [0,1]\]
is a diffeomorphims between manifolds with boundary. Two special cases of stratified cobordisms will be of special interest.
The first is the case in which $X$ is a compact manifold with two boundary components $X_1$ and $X_2$. Then our definitions coincide
with the usual ones. We will say that in this case we have a (non-stratified) {\em cobordism} between $X_1$ and $X_2$. 
The second is when $X$ is a manifold with corners and the boundary $\partial X$ admits a decomposition
$\partial X=X_0\cup X_1\cup Y$, where $X_0$, $X_1$ and $Y$ are manifolds with boundary and we have
$\partial X_0=X_0\cap Y$, $\partial X_1=X_1\cap Y$. Then $\partial Y=\partial X_0\cup \partial X_1$ and $(Y,Y\cap X_0,Y\cap X_1)$ 
is a cobordism. We say that $(X,X_0,X_1)$ is a {\em cobordism with boundary} between $X_0$ {\em and} $X_1$, and that 
$(Y,Y\cap X_0,Y\cap X_1)$ is the {\em boundary cobordism}.\\

\noindent\textbf{Setting and notation}.
Thorought the rest of this section we denote by 
\[\pi:E\to U\] 
be a smooth complex vector bundle of rank $n$ over a manifold $U$.
Let 
\[f:E\to\CC\] 
be a smooth germ at $\{O\}\times U$ whose restrictions to the fibres of $\pi$ are holomorphic. Assume that 
the Milnor number $\mu(f|_{E_p})$ is independent of $p\in U$. Denote by 
\[\phi:E\to U\times\CC\]
the mapping $\phi:=(\pi,f)$. Fix a hermitian metric for the vector bundle $E$ and let $\rho$ be the associated (fibrewise)
''distance to the origin''
function. Consider $K\subset U$ and let $\theta:K\to\RR_+$ be any positive continuous function. We define the 
$\theta$-neighborhood of $K$ as 
\[B(K,\theta):=\{x\in E:\pi(x)\in K\quad\text{and}\quad\rho(x)\leq\theta(\pi(x))\}.\]
Given subsets $A\subset U$ and $B\subset E$, we will denote the intersection $B\cap\pi^{-1}(A)$ by $B_A$ or by $B|_A$.

\subsection{Cuts}
\label{inicial}

\begin{definition}
\label{cuts}
A {\em cut for} $f$ {\em over a submanifold} (possibly non-closed and with corners) $V\subset U$ 
with amplitude $\delta$ is a closed smooth hypersurface $H$ of 
$\phi^{-1}(V\times D_\delta)$ with the following properties:
\begin{enumerate}
\item For any $(t,s)\in V\times D_\delta$ the level set $\phi^{-1}(t,s)$ is smooth outside the zero section of $\pi$ 
and transverse to $H$ at any of their intersection points.
\item There is a unique connected component of $\phi^{-1}(V\times D_\delta)\setminus H$ containing the zero-section of $E$. 
Its closure $Y_{int}(H,V,\delta)$ is 
a smooth manifold with corners, called the {\em interior component}, 
and the restriction 
\[\pi:Y_{int}(H,V,\delta)\to V\] 
is a smooth locally trivial fibration with contractible fibres.
\item For any $t\in V$ and any radius $\epsilon_t$ for the Milnor ball of $f_t$ the space 
\[\phi^{-1}(t,0)\cap Y_{int}(H,V,\delta)\setminus\dot{B}((O,t),\epsilon_t)\] 
is a smooth trivial cobordism. 
\end{enumerate}
\end{definition}

Usually $V$ will be either an open subset of $U$ or the closure of an open subset.
If $H$ is a cut over $V$ of amplitude $\delta$ then for any $V'\subset V$ submanifold and $\delta'<\delta$, the hypersurface 
$H\cap\phi^{-1}(V'\times D_{\delta'})$ is a cut over $V'$ with amplitude $\delta'$, we will abuse notation and denote 
by $Y_{int}(H,V',\delta')$ the interior component of the new cut. 
We need the following relation in the set of cuts: let $H$ and $H'$ be cuts over open subsets $V$ and $V'$ with amplitudes 
$\delta$ and $\delta'$ respectively. We say that $H\preceq H'$ if 
\[Y_{int}(H,V\cap V',\min\{\delta,\delta'\})\subset Y_{int}(H',V\cap V',\min\{\delta,\delta'\}).\]
The strict version of the relation is: we say that $H\prec H'$ if
\[Y_{int}(H,V\cap V',\min\{\delta,\delta'\})\subset \dot{Y}_{int}(H',V\cap V',\min\{\delta,\delta'\}).\]
Observe that $H_V\preceq H_W$ if $V\cap W=\emptyset$. Beware that the relation is not an ordering, and that given
two cuts related in 
both directions, they are not necessarily equal.

The following Lemma shows that cuts appear naturally
\begin{lema}
\label{milnorcut}
Suppose $n\neq 3$.
For a fixed $t\in U$ we consider a pair $(\epsilon,\delta)$ of radii for the Milnor fibration of $f_t$, There exists an open 
neighborhood $V$ of $t\in U$ such that $\phi^{-1}(V\times D_\delta)\cap \partial B(V,\epsilon)$ is a cut over $V$ of amplitude 
$\delta$. 
\end{lema}
\begin{proof}
Conditions (1) and (2) are easy. Condition (3) is contained in the proof of L\^e-Ramanujam Theorem~\cite{LR}.
\end{proof}

\begin{lema}
\label{pointtolocal}
Suppose $n\neq 3$. If we have an open subset $V\subset U$, a positive $\delta$ and a closed smooth hypersurface $H_V$ of $\varphi^{-1}(V\times D_\delta)$, and $H_V$ satisfies the
conditions of a cut at a point $t\in V$, then there is an open neighborhood $V'$ of $p$ in $V$ such that
$H_V\cap\varphi^{-1}(V'\times D_\delta)$ is a cut over $V'$ of amplitude $\delta$.
\end{lema}
\begin{proof}
The proof is an application of $h$-cobordism Theorem as in the proof of L\^e-Ramanujam Theorem~\cite{LR}.
\end{proof}

\subsection{Extension of cuts}
\label{segunda}

\begin{definition}
\label{extension}
Let $V_1\subset V_2$ be submanifolds of $U$ and $H_1$, $H_2$ cuts of the same amplitude $\delta$ 
over $V_1$ and $V_2$ respectively. We say that $H_2$ is an extension of $H_1$, if we have $H_2\cap \phi^{-1}(V_1\cap D_\delta)=H_1$.
\end{definition}

We consider the following situation:
define $I:=[0,1]$. Denote by $C$ the $d$-dimensional cube $I^d$, with its natural structure of manifold
with corners, and
$U:=(-\epsilon,1+\epsilon)^d$ for a certain $\epsilon>0$.

All the cuts below are chosen with the same amplitude $\delta$.
Consider two cuts $H_+$ and $H_-$ defined over $U$, with $H_-\prec H_+$. Let $A$ be a contractible union of faces of $C$. We consider 
a cut $H_0$ over a neighborhood $V$ of $A$, which satisfies $H_-\prec H_0\prec H_+$. We consider also a number of pairs 
\[(B_1,\mu_1),...(B_{k_1},\mu_{k_1})\quad\quad (C_1,\nu_1),...(C_{k_2},\nu_{k_2}),\]
where the $B_i$'s and the $C_i$'s are faces of $C$, each $\mu_i$ denotes a cut over a {\em contractible} 
neighborhood $B'_i$ of $B_i$ in $U$ 
and each $\nu_i$ denotes a cut over a {\em contractible}
neighborhood $C'_i$ of $C_i$ in $U$. Suppose that the following relations hold: 
$H_-\prec\mu_i$ for any $i$, $H_-\prec\nu_i$ for any $i$, $\mu_i\prec\mu_j$ if $i<j$,
$\mu_i\prec H_0$ for any $i$, $\mu_i\prec\nu_j$ for any $i,j$, $\mu_i\prec H_+$ for any $i$, $H_0\prec\nu_i$ for any $i$,
$\nu_i\prec \nu_j$ if $i>j$, and $\nu_i\prec H_+$ for any $i$.

We will need the following extension Lemma:

\begin{lema}
\label{extcut}
In the setting above, after possibly shrinking $U$ and $V$ to smaller neighborhoods of the corresponding sets, 
and possibly decreasing $\delta$, there exists a cut 
$H'_0$ over $U$ with amplitude $\delta$ 
which extends $H_0$ and satisfies $H_-\prec H_0$, $\mu_i\prec H_0$ for any $i$, $H_0\prec \nu_i$ for any $i$, and $H_0\prec H_+$.
\end{lema}
\begin{proof}
Define 
\[Y:=\overline{Y_{int}(H_+,U,\delta)\setminus Y_{int}(H_-,U,\delta)}.\] 
By Ehresmann Fibration Theorem and the contractibility 
of $U\times D_\delta$ the mapping 
\[\phi:Y\to U\times D_\delta\]
is a trivial fibration. From the defining properties of a cut and an easy manipulation with cobordisms we deduce that each 
fibre of $\phi$ is a trivial cobordism. 

We will use the following fact: 

\textsc{Fact 1}: let $A$ be a contractible closed union of faces of $C$ and $V$ be
a neighborhood of $A$ in $U$. There exist 
a compact neighborhood $D_1$ of $A$ in $V$, whose boundary $\partial D_1$ is a smooth hypersurface in 
$U$, an open neighborhood $D_2$ of $C$ in $U$ and a smooth mapping 
\begin{equation}
\label{auxiliar}
\xi:D_2\setminus \dot{D}_1\to\partial D_1
\end{equation}
which is a trivial fibration over its image, which is contractible,
with fibre diffeomorphic to the half-open interval $[0,1)$ (the closed extreme 
of the fibres corresponding with the intersection of the fibre with $\partial D_2$).
Let 
\begin{equation}
\label{yacasi}
\Psi:D_2\setminus\dot{D}_1\to [0,1)\times\xi(D_2\setminus\dot{D}_1)
\end{equation}
be a diffeomorphism giving a trivialisation.

This fact follows from easy geometrical considerations. As an example we do in detail 
the case in which $A=[0,1]^{d_1}\times \{0\}^{d_2}$, with 
$d_1+d_2=d$. Let $\sigma$ denote the usual distance function in $\RR^d$. We may consider 
\[D_1:=\{p\in\RR^d:\sigma(p,A)\leq\epsilon_1\},\]
\[D_2:=\{p\in\RR^d:\sigma(p,C)<\epsilon_2\},\]
for $0<<\epsilon_2<<\epsilon_1<<1$, and define $\xi(x_1,...,x_d)$ to be the only intersection point
of $\partial D_1$ with the segment joining $(x_1,...,x_d)$ with $(x_1,...,x_{d_1},0,...,0)$. We encourage the reader to draw a 
picture of this construction.

We shall proceed in two steps.

\textbf{Step 1}: we shall inductively reduce to the case in which $k_1+k_2=0$. Suppose that $k_1+k_2>0$. Assume that $k_2>0$.

Choose any point $t\in C_1$. From the defining properties of a cut, an easy manipulation with cobordisms, Ehresmann Fibration 
Theorem and the fact that $C'_1$ is contractible, it follows that the restrictions
\[\phi:\overline{Y_{int}(H_+,C'_1,\delta)\setminus Y_{int}(\nu_1,C'_1,\delta)}\to C'_1\times D_\delta,\]
\[\phi:\overline{Y_{int}(\nu_1,C'_1,\delta)\setminus Y_{int}(H_-,C'_1,\delta)}\to C'_1\times D_\delta,\]
are smooth trivial fibrations with fibre a trivial smooth cobordism. Hence, it is easy to construct a 
smooth vector field $\calY$ in $Y|_{C'_1\times D_\delta}$  which is tangent to the fibres of $\phi$ such that 
its associated flow defines a diffeomorphism $\psi:(H_-)|_{C'_1\times D_\delta}\times [0,1]\to Y_{C'_1}$ satisfying the 
equality
\begin{equation}
\label{referir}
\overline{Y_{int}(\nu_1,C'_1,\delta)\setminus Y_{int}(H_-,C'_1,\delta)}=\psi((H_-)|_{C'_1\times D_\delta}\times [0,1/2]).
\end{equation}

We are going to extend $\calY$ to a vector field $\calX$ defined in $Y$.
As $C_1$ is a face, we may choose $D_1$ to be the neighborhood of $C_1$ in $C'_1$, and $D_2$ the neighborhood of $C$ in $U$ 
so that the mapping~(\ref{auxiliar}) exists as predicted by Fact~1.
Using that $\phi:Y\to U\times D_\delta$ is a trivial fibration, it is easy to prove that there exists a smooth diffeomorphism
\[\Phi:Y_{D_2\setminus\dot{D}_1}\to [0,1)\times Y_{\xi(D_2\setminus\dot{D}_1)}\]
such that, if $q_i$ denotes the projection of $[0,1)\times Y_{\xi(D_2\setminus\dot{D}_1)}$ to the $i$-th factor ($i=1,2$), then 
\begin{enumerate}[(i)]
\item the mappings
\[(\Psi,Id_{D_\delta})\comp\phi:Y_{D_2\setminus\dot{D}_1}\to [0,1)\times\xi(D_2\setminus\dot{D}_1)\times D_\delta\]
and 
\[ (q_1,\phi|_{Y_{\xi(D_2\setminus\dot{D}_1)}})\comp\Phi:Y_{D_2\setminus\dot{D}_1}\to [0,1)\times\xi(D_2\setminus\dot{D}_1)\times D_\delta\]
coincide,  
\item the restriction $q_2\comp\Phi|_{Y_{\xi(D_2)}}$ is the identity; moreover, the mappings 
\[(0,Id_{Y_{D_1}}):Y_{D_1}\to [0,1)\times Y\]
and $\Phi$ glue to a smooth mapping.
\end{enumerate}

We define a smooth vector field $\calX$ in $Y_{D_2}$ piecewise as follows: the restriction of $\calX$ to 
$Y_{D_1}$ is equal to the restriction of $\calY|_{Y_{D_1}}$; the restriction of $\calX$ to $Y_{(D_2\setminus\dot{D}_1)}$
is the pullback by $\Psi$ of the vector field in $[0,1)\times Y_{\xi(D_2\setminus\dot{D}_1)}$ whose $[0,1)$-component is zero and 
whose $Y_{\xi(D_2\setminus\dot{D}_1)}$-component is equal to the restriction $\calY|_{Y_{\xi(D_2\setminus\dot{D}_1)}}$.
Observe that $\calX$ is
smooth, tangent to the fibres of $\phi$, extends $\calY$, and its flow trivialises the cobordism given by each fibre of $\phi$.

We redefine $U:=D_2$. Then the flow of $\calX$ induces a diffeomorphism 
\[\varphi:H_-\times [0,1] \to Y.\]

Let $\{\theta_1,\theta_2\}$ be a partition of unity subordinated to $\{C'_1\times D_\delta,(U\setminus C_1)\times D_\delta\}$. 
Define $\beta:=1/2\theta_1+\theta_2$.
Let $\sigma:Y\to [0,1]$ be the smooth function $\sigma:=p_2\comp \varphi^{-1}$, with $p_2$ the projection of $H_-\times [0,1]$ to 
the second factor.
We claim that if $C'_1$ is a sufficiently small neighborhood of $C_1$ in $U$, then:
\begin{enumerate}[(a)]
\item for any $1\leq i\leq k_1$ and any $x\in \mu_i|_{B_i}$ we have $\sigma(x)<\beta(\pi(x))$,
\item for any $2\leq j\leq k_2$ and any $x\in \nu_j|_{C_i}$ we have $\sigma(x)<\beta(\pi(x))$, and 
\item for any $x\in H_0|_{A}$ we have $\sigma(x)<\beta(\pi(x))$.
\end{enumerate}

Consider $x\in\mu_i|_{B_i\cap D_1}$.
Observe that we have Equality~(\ref{referir}) together with the fact that $\mu_i\prec\nu_1$ implies that $\sigma(x)<1/2$. 
By continuity this inequality occurs in a neighborhood of $x$ in $\mu_i$.
On the other hand $\sigma(x)<1$ for any $x\in\mu_i$. Hence, using that $B_i\times D_\delta$ is compact, it is easy to see that ,
by shrinking $C'_1$, we have 
$\sigma(x)<\beta(\phi(x))$ for any $x\in \mu_i$. This shows property (a). The proofs of properties (b) and (c)
are analogous.

Now, by continuity, shrinking $B'_i$, $C'_j$ and $V$, we have 
\begin{enumerate}[(a)]
\item for any $1\leq i\leq k_1$ and any $x\in \mu_i|_{B'_i}$ we have $\sigma(x)<\beta(\pi(x))$,
\item for any $2\leq j\leq k_2$ and any $x\in \nu_j|_{C'_i}$ we have $\sigma(x)<\beta(\pi(x))$, and 
\item for any $x\in H_0|_{V}$ we have $\sigma(x)<\beta(\pi(x))$.
\end{enumerate}

By construction the subspace  $\nu'_1:=\{\varphi(x,\beta(\pi(x))),x\in H_-\}$ is a closed smooth hypersurface of 
$\pi^{-1}(U)\cap f^{-1}(D_\delta)$ which constitutes a cut over $U$ and satisfies
$H_-\prec\nu'_1$, $\mu_i\prec\nu'_1$ for any $i$, $H_0\prec\nu'_1$ and $\nu_i\prec\nu'_1$ for any $i\neq 1$. 

Then we can redefine $H_+:=\nu'_1$ and we have  
$k_1+k_2$ strictly smaller. We have worked out the case in which $k_2>0$. The case in which $k_1>0$ is 
identical. Inductively we are reduced to the case $k_1+k_2=0$.

\textbf{Step 2}: We assume $k_1=k_2=0$.

As $A$ is a contractible closed union of faces of $C$, a simpler version of the argument of Step 1 allows to construct, 
after eventually 
shrinking $U$ and $V$, a cut $H'_0$ over $U$ which extends $H_0$ and satisfies
\[H_-\prec H'_0\prec H_+.\]
\end{proof}

\subsection{Existence of cuts}
\label{extcutii}

In this section we let $L$ denote either a point or the circle $\SSS^1$. Again $C$ denotes the cube $I^d\subset\RR^d$. Let 
$U:=L\times (-\epsilon,1+\epsilon)^d$ for a certain $\epsilon>0$. We shall make the following assumption:\\

\noindent\textbf{Assumption A}:
there exists $\zeta>0$ such that $\rho^{-1}(\epsilon)$ meets $\phi^{-1}(t,0)$ {\em transversely} for any 
$t\in L\times\{(0,...,0)\}$ and any positive $\epsilon\leq\zeta$.\\

\begin{remark}
\label{unpunto}
Observe that, as the restriction of $\rho$ to each fibre of $\pi$ is real analytic, the assumption is satisfied automatically if 
$L$ is a point (see~\cite{Lo}), pages 21-25.
\end{remark}

\begin{lema}
\label{cutcircle}
Let $\theta:U\to\RR_+$ be a positive continuous function. There exists a neighborhood $V$ of $L\times\{(0,...,0)\}$ in
$U$, and positive $\epsilon$ and $\delta$ such that 
$H_V:=\rho^{-1}(\epsilon)\cap\phi^{-1}(V\times D_\delta)$ is a cut over $V$ of amplitude $\delta$ such that 
$Y_{int}(H_V,V,\delta)$ is contained in $B(U,\theta)$.
\end{lema}
\begin{proof}
The proof is easy after Lemmas~\ref{milnorcut} and \ref{pointtolocal}.
\end{proof}

\begin{prop}
\label{usefulconstr}
Let $\theta:U\to \RR_+$ be any continuous function. 
After possibly shrinking $U$ to a smaller neighborhood of $L\times C$ in $L\times\RR^d$, 
there exits positive $\epsilon$, $\delta$, a neighborhood $V$ of $L\times\{(0,...,0)\}$ in $U$ 
and a cut $H$ over $U$ of amplitude $\delta$ such that 
\begin{enumerate}
\item the space $Y_{int}(H,U,\delta)$ is contained in $B(U,\theta)$,
\item the cut $H$ extends the cut $H_V:=\rho^{-1}(\epsilon)\cap\phi^{-1}(V\times D_\delta)$ obtained in Lemma~\ref{cutcircle}.
\end{enumerate}
\end{prop}

\begin{proof}
We will assume $L=\SSS^1$. The proof for the case in which $L$ is a point is a simplification of the 
case presented here.

For any $t\in L\times C$ 
we choose a positive $\epsilon_t<\theta(t)$ such that $\rho|_{\phi^{-1}(t,0)\cap\rho^{-1}((0,\epsilon_t])}$ is a
submersion. By Lemma~\ref{milnorcut} there exists a neighborhood $U_t$ of 
$t$ in $U$ and a positive $\delta_t$ such that $H_+(U_t):=\phi^{-1}(U_t\times D_{\delta_t})\cap\rho^{-1}(\epsilon_t)$ is a cut
over $U_t$ with amplitude $\delta_t$ and $Y_{int}(H_+(U_t),U_t,\delta_t)$ is contained in $B(U,\theta)$.

By compactness of $L$ and $C$ there exists a finite cover $\{U_{t}\}_{t\in A}$ of $L\times C$ which 
contains a collection $\{U_{t}\}_{t\in A_0}$, with $A_0\subset A\subset U$, which covers $L\times\{(0,...,0)\}$ 
and satisfies that $A_0\subset L\times\{(0,...,0)\}$.

Let $\gamma:I\to\SSS^1$ be the parametrisation defined by $\gamma(t)=e^{2\pi it}$. Given a partition
\[0=\alpha_0<\alpha_1<...<\alpha_r=1\]
we define $I_j:=[\alpha_{j-1},\alpha_j]$ and also, for any multi-index $J=(j_0,...,j_{d})$, the cube 
\[C_J:=\gamma(I_{j_0})\times\prod_{r=1}^d I_{j_r},\]
which is contained in $L\times C$. Choose the partition so fine that any cube $C_J$ 
is contained at least in an open set of the cover $\{U_{t}\}_{t\in A}$, and moreover, 
if $C_J$ meets $L\times\{(0,...,0)\}$, then it is contained in an open set of the cover 
$\{U_{t}\}_{t\in A_0}$. We assign to each $C_J$ a fixed set $U_{t(J)}$ of the 
cover containing it, taking care that if $C_{J}$ meets $L\times\{(0,...,0)\}$ then $t(J)$ belongs to $A_0$.
We denote $U_{t(J)}$ by $U_{J}$. 
Define $\epsilon_{J}:=\epsilon_{t(J)}$. Define 
$\delta_{J}:=\delta_{t(J)}$. Choose $\delta_{min}$ and $\epsilon_{min}$ strictly smaller than any 
$\delta_J$ and $\epsilon_J$ respectively.

By Lemma~\ref{cutcircle} there exists a neighborhood $V$ of $L\times\{(0,...,0)\}$ in
$U$, positive $\epsilon<\epsilon_{min}$, $\delta<\delta_{min}$ such that $H_V:=\rho^{-1}(\epsilon)\cap\phi^{-1}(V\times D_\delta)$ 
is a cut over $V$ of amplitude $\delta$.

Choose a positive $\eta$ strictly smaller than $\epsilon$.
For any $t\in L\times C$ we choose a positive $\xi_t<\eta$ such that 
$\rho|_{\phi^{-1}(t,0)\cap\rho^{-1}((0,\xi_t])}$ is a submersion. 
By Lemma~\ref{milnorcut}, for any $t\in L\times C$, there exists a neighborhood $V_{t}$ of 
$t$ in $U$ and a positive $\delta'_t<\delta_{min}$ 
such that $H_-(V_t):=\phi^{-1}(U_t\times D_{\delta'_t})\cap\rho^{-1}(\xi_t)$ is a cut over 
$V_t$ with amplitude $\delta'_t$. 
We pick $V_t$ so small that it is contained in any $U_{J}$ which contains $t$.

By compactness we can choose a finite cover $\{V_{t}\}_{t\in B}$ of $L\times C$.

For each $j$, given any subdivision 
\[\alpha_j=\alpha_{j,0}<\alpha_{j,1}<...\alpha_{j,s}=\alpha_{i+1},\]
we define the intervals $I_{j,k}:=[\alpha_{j,k-1},\alpha_{j,k}]\subset I$ 
and, for multi-indexes $J=(j_0,...,j_d)$ and $K=(k_0,...,k_d)$ we define the cubes  
\[C_{J}^{K}:=\gamma(I_{j_0,k_0})\times \prod_{r=1}^d I_{i_r,j_r}.\] 
Each cube $C_{J}$ splits in an 
union of the cubes $C_{J}^{K}$, where $K$ vary.

We choose the subdivisions so fine that each cube 
$C_{J}^{K}$ is contained in an open 
subset of the cover $\{V_{t}\}_{t\in B}$. We assign to each cube $C_J^K$ 
a fixed open subset $V_{t_J^K}$ containing it. We denote 
$V_{t_J^K}$ by $V_J^K$, and define 
$\xi_J^K:=\xi_{t_J^K}$,
$\epsilon_J^K:=\epsilon_J$ and 
$\delta_J^K:=\delta'_{t_J^K}$.

Given $\zeta>0$ we define the intervals $\tilde{I}_{i,j}:=(\alpha_{i,j-1}-\zeta,\alpha_{i,j}+\zeta)$ and 
the open cubes
\[\tilde{C}_{J}^{K}:=\gamma(\tilde{I}_{j_0,k_0})\times\prod_{r=1}^d\tilde{I}_{j_r,k_r}.\]
Choose $\zeta$ small enough so that
\begin{enumerate}[(i)]
\item the closure of $\tilde{C}_{J}^{K}$ is contained in 
$V_{J}^{K}$, and hence in $U_{J}$,
\item for any choice of indexes, the sets $\tilde{C}_{J}^{K}$ and 
$\tilde{C}_{J'}^{K'}$ meet if and only if 
$C_{J}^{K}$ and $C_{J'}^{K'}$ meet.
\end{enumerate}

Condition $(i)$ and the compactness of the closure of $\tilde{C}_{I}^{J}$ imply that we can shrink
$\delta_{J}^{K}$ such that
\[H_+(\tilde{C}_{J}^{K}):=\phi^{-1}(\tilde{C}_{J}^{K}\times D_{\delta_{J}^{K}})\cap\rho^{-1}(\epsilon_{J}^{K})\]
and
\[H_-(\tilde{C}_{J}^{K}):=\phi^{-1}(\tilde{C}_{J}^{K}\times D_{\delta_{J}^{K}})\cap\rho^{-1}(\xi_{J}^{K})\]
are cuts over $\tilde{C}_{J}^{K}$ with amplitude $\delta_{J}^{K}$.

We choose $\delta>0$ smaller than $\delta_{J}^{K}$ for any choice of multi-indexes $(J,K)$.

Given a fixed choice of multi-indexes $(J,K)$ we consider the following cuts:
\begin{itemize}
\item for any indexes $(J',K')$ different from the fixed ones, and such that
$\epsilon_{J'}^{K'}<\epsilon_{J}^{K}$ we consider
the cut 
\[\nu_{J'}^{K'}(\tilde{C}_{J}^{K}):=\phi^{-1}((\tilde{C}_{J}^{K}\cap\tilde{C}_{J'}^{K'}) \times D_{\delta})\cap\rho^{-1}(\epsilon_{J'}^{K'})\]
over $\tilde{C}_{J}^{K}\cap\tilde{C}_{J'}^{K'}$ with amplitude $\delta$,
\item for any indexes $(J',K')$ different from the fixed ones, and such that
$\xi_{J'}^{K'}>\xi_{J}^{K}$ we consider
the cut 
\[\mu_{J'}^{K'}(\tilde{C}_{J}^{K}):=\phi^{-1}((\tilde{C}_{J}^{K}\cap\tilde{C}_{J'}^{K'}) \times D_{\delta})\cap\rho^{-1}(\xi_{J'}^{K'})\]
over $\tilde{C}_{J}^{K}\cap\tilde{C}_{J'}^{K'}$ with amplitude $\delta$.
\end{itemize}

Because of Condition (i) and the compactness of the closure of $\tilde{C}_J^K$ it is clear that a small modification of
$\epsilon_J^K$ and $\xi_J^K$ keeps all the properties obtained above. An adequate
perturbation making the $\epsilon^K_J$'s and the $\xi^K_J$'s pairwise different gives rise to the following properties:

The cuts $\nu_{J'}^{K'}(\tilde{C}_{J}^{K})$ obtained by letting $(J',K')$ vary are linearly related by 
$\preceq$ and, moreover, if two multi-indexes $(J',K')$ and $(J'',K'')$ are such that 
$\tilde{C}_{J}^{K}$ and $\tilde{C}_{J''}^{K''}$ meet then, if the cuts $\nu_{J'}^{K'}(\tilde{C}_{J}^{K})$
and $\nu_{J''}^{K''}(\tilde{C}_{J}^{K})$ are defined, they are necessarily related by $\prec$ in one of 
the directions. The same happens for the cuts 
$\mu_{J'}^{K'}(\tilde{C}_{J}^{K})$ obtained by letting
 $(J',K')$ vary. We have moreover  
\[\mu_{J'}^{K'}(\tilde{C}_{J}^{K})\prec \nu_{J''}^{K''}(\tilde{C}_{J}^{K}),\]
\[H_-(\tilde{C}_{J}^{K})\prec \mu_{J'}^{K'}(\tilde{C}_{J}^{K}),\]
\[\mu_{J'}^{K'}(\tilde{C}_{J}^{K})\prec H_+(\tilde{C}_{J}^{K}),\]
\[H_-(\tilde{C}_{J}^{K})\prec \nu_{J'}^{K'}(\tilde{C}_{J}^{K}),\]
\[\nu_{J'}^{K'}(\tilde{C}_{J}^{K})\prec H_+(\tilde{C}_{J}^{K})\]
for any choice of multi-indexes for which the expression makes sense.

Shrink $V$ to a neighborhood of $L\times\{(0,...,0)\}$ in $U$ which is small enough that it only meets
$\tilde{C}_{J}^{K}$ if $C_{J}^{K}$ meets 
$L\times\{(0,...,0)\}$. If $C_{J}^{K}$ meets $L\times\{(0,...,0)\}$
we consider the cut 
\[H_{\tilde{C}_{J}^{K}\cap V}:=\phi^{-1}((\tilde{C}_{J}^{K}\cap V)\times D_{\delta})\cap H_V\]
obtained by restriction over $\tilde{C}_{J}^{K}\cap V$ of the previously constructed cut $H_V$.

By our constructions (choosing the perturbation of the $\epsilon^K_J$'s and the $\xi^K_J$'s small enough) we have:
\begin{equation}
\label{rel1}
H_-(\tilde{C}_{J}^{K})\prec H_{\tilde{C}_{J}^{K}\cap V},
\end{equation}
\begin{equation}
\label{rel2}
H_{\tilde{C}_{J}^{K}\cap V}\prec H_+(\tilde{C}_{J}^{K}),
\end{equation}
\begin{equation}
\label{rel3}
\mu_{J'}^{K'}(\tilde{C}_{J}^{K})\prec H_{\tilde{C}_{J}^{K}\cap V},
\end{equation}
\begin{equation}
\label{rel4}
H_{\tilde{C}_{J}^{K}\cap V}\prec \nu_{J'}^{K'}(\tilde{C}_{J}^{K}).
\end{equation}

We shall construct the cut $H$ extending $H_V$ step by step over the cubes $C_{J}^{K}$
using Lemma~\ref{extcut}. 
Order the cubes $C_{J}^{K}$ lexicographically in $(j_0,...,j_d,k_0,...,k_d)$.

We start with the first cube, that is $(J,K)=(0,...,0)$. 
Observe that 
\[A:=(\SSS^1\times\{(0,...,0)\})\cap U_{(0,...,0)}^{(0,...,0)}=\gamma(I_0)\times \{0\}^{d}\] 
is a contractible union of 
faces of $C_{0,...,0}^{0,...,0}$. We apply Lemma~\ref{extcut} and obtain a cut  
$H(\tilde{C}_{(0,...,0)}^{(0,...,0)})$ which extends 
$H_{\tilde{C}_{(0,...,0)}^{(0,...,0)}\cap V}$ and satisfies the relations analogous to
(\ref{rel1})-(\ref{rel4}) with respect to the cuts 
$H_-(\tilde{C}_{(0,...,0)}^{(0,...,0)})$, $H_+(\tilde{C}_{(0,...,0)}^{(0,...,0)})$,
$\mu_{J'}^{K'}(\tilde{C}_{(0,...,0)}^{(0,...,0)})$ and 
$\nu_{J'}^{K'}(\tilde{C}_{(0,...,0)}^{(0,...,0)})$, for 
$(J',K')$ varying. We redefine  $V:=V\cup\tilde{C}_{0,...,0}^{0,...,0}$ and the 
cut $H_V$ as the union of the previous $H_V$ and $H(\tilde{C}_{(0,...,0)}^{(0,...,0)})$. 
By construction, the new cut $H_V$ satisfies the following property: 
for any $(J,K)$ the cut $H_{\tilde{C}_{J}^{K}\cap V}$ given by the restriction of $H_V$ over 
$\tilde{C}_{J}^{K}\cap V$ satisfies relations (\ref{rel1})-(\ref{rel4})

The inductive step runs similarly.
\end{proof}

\begin{remark}
\label{basepoint}
Choose any point $p\in C$.
The above proposition is valid (with the same proof, up to some notational complication) 
if we let $L\times\{p\}$ play the role of $L\times\{(0,...,0)\}$ both in Assumption A and in the statement of the Proposition. 
\end{remark}

\subsection{Topological equisingularity of $\mu$-constant families} 

Let $\theta_1:U\to \RR_+$ be any continuous function. We choose a positive $\delta_1$ and a cut $H_1$ as predicted by  Proposition~\ref{usefulconstr}. In order to simplify the notation we denote $Y_{int}(H_1,U,\delta)$ simply by $Y_1$.
We view $U$ embedded in $Y_1$ as the zero section of the vector bundle $\pi:E\to U$. Define $Y_1^*:=Y_1\setminus f^{-1}(0)$.

\begin{lema}
\label{vectorfield}
There exists a smooth 
vector field $\calX$ in $Y|_{L\times C}\setminus (L\times C)$ with the following properties:
\begin{enumerate}
\item it is tangent to the fibres of $\pi$, 
\item There exists a vector field $\calW$ in $D_\delta^*$, which is radial, pointing to the origin and of modulus 
$|\calW(z)|\leq |z|^2$, such that $df(\calX)(x)=\calW(f(x))$ for any $x\in Y_1^*$,
\item the vector field $\calX$ is tangent to $f^{-1}(0)$,
\item any integral curve converges to the
origin of the fibre of $\pi$ in which it lies in infinite time.
\end{enumerate}
\end{lema}
\begin{proof}
We will construct $\calX$ as the amalgamation of two vector fields $\calY$ and $\calZ$. We define each of them separately.
Denote by $\calV$ the vector field in $D_\delta^*$ which is radial, pointing at the origin and of modulus $||\calV(z)||=||z||^{2}$.
Consider in $U\times D_\delta^*$ the vector field $\calV'$ characterised by being tangent to the fibres of the projection of
$U\times D_\delta^*$ to the first factor, and a lift 
by the projection of $U\times D_\delta^*$ to the second factor of the vector field $\calV$.
As $\phi:Y_1^*\to U\times D_\delta^*$ is submersive we can define the vector field $\calY$ in 
$Y_1^*$ to be a lifting of $\calV'$ by the mapping $\phi$.

For the definition of $\calZ$ some auxiliary constructions are needed. We will shrink $U$ when necessary without explicitly
mentioning it.
There is a continuous function $\theta_2:U\to\RR$ such that the inclusion
$B(U,\theta_2)\subset Y_1$ is satisfied. Applying Proposition~\ref{usefulconstr} to the function $1/2\theta_2$ 
we obtain new positive $\delta_2$ (which we choose to be smaller than $\delta_1/2$)
and a cut $H_2$ of amplitude $\delta_2$ such that $Y_2:=Y_{int}(H_2,U,\delta_2)$
is contained in $B(U,1/2\theta_2)$. 

We iterate this procedure to obtain an infinite sequence of cuts $H_i$ over $L\times C$ 
with amplitude $\delta_i$, which gives rise to an infinite, nested, sequence
\begin{equation}
\label{nestedset}
Y_1\supset Y_2\supset ...\supset Y_i\supset ...
\end{equation}
of closed neighborhoods of $L\times C$, satisfying that $\cap_{i=1}^\infty Y_i$ is equal to the zero section.

Define $Z_i:=Y_i\cap f^{-1}(D_{\delta_{i+1}})\setminus \dot{Y}_{i+1}$.
By Property (1) of Definition~\ref{cuts} the restriction
\[\phi_i:=\phi|_{Z_i}:Z_i\to L\times C\times D_{\delta_{i+1}}\]
is a locally trivial fibration. By Property (3) of Definition~\ref{cuts} and an easy manipulation with cobordisms, 
its generic fibre is a trivial cobordism.

By the second property that Proposition~\ref{usefulconstr} predicts for the cut $H_i$ there is a neighborhood 
$V_i$ of $L\times\{(0,...,0)\}$ in $U$ and a positive $\epsilon_i$ such that 
\[H_i\cap\pi^{-1}(V_i)=\phi^{-1}(V_i\times D_{\delta_i})\cap \rho^{-1}(\epsilon_i).\]
We have that $\{\epsilon_i\}_{i\in\NN}$ is a decreasing sequence converging to $0$. 

Observe that we have 
\[Z'_i:=Z_i\cap\pi^{-1}(V_{i+1})=\phi^{-1}(V_{i+1}\times D_{\delta_{i+1}})\cap\rho^{-1}([\epsilon_{i+1},\epsilon_i]).\]
Using  Assumption A, and possibly shrinking 
$V_{i+1}$ and $\delta_{i+1}$ we obtain that the restriction of 
$\rho$ to each fibre of $\phi_i|_{Z'_i}$ is a submersion. Therefore there exists a vector field 
$\calZ'_i$ on $Z'_i$ which is tangent to the fibres of 
$\phi_i|_{Z'_i}$, non-zero at any point, and such that $d\rho(\calZ'_i)=-1$. Consequently,
 its flow takes $(Y_{i}\cap H_{i}\cap f^{-1}(D_{\delta_{i+1}}))_{V_i}$ to $(Y_{i+1}\cap H_{i+1})_{V_i}$ in finite time 
$T_i=\epsilon_i-\epsilon_{i+1}$.

We may suppose that $V_i$ is of the form $L\times B_i$, with $B_i$ a small ball centered around 
$(0,...,0)$ in $\RR^d$. Using that $C$ is a cube and the fact that $\phi_i$ is a locally trivial 
fibration, we may (in the same way that we extended $\calY$ to $\calX$ in Step 1 of the proof of 
Lemma~\ref{extcut}) extend $\calZ'_i$ to a 
vector field $\calZ_i$ on $Z_i$ 
which is tangent to the fibres of 
$\phi_i$, non-zero at any point, and whose flow takes $Y_{i}\cap H_{i}\cap f^{-1}(D_{\delta_{i+1}})$ to
$Y_{i+1}\cap H_{i+1}$ in time $T_i$. We rescale the vector field in order to ensure that 
$T_i:=1$. 

As, for any $i\in\NN$, both $\calZ_i$ and $\calZ_{i+1}$ are transverse to $H_{i+1}$, and both of them 
point into $Y_{i+1}$. Using an adequate partition of unity, we may glue the vector fields 
$\{\calZ_i\}_{i\in\NN}$ to a vector field $\calZ$ defined in
$\cup_{i\in\NN} Z_i$, nowhere vanishing and tangent to the fibres of $\phi$. Up to a rescaling 
of $\calZ$ we may assume that its flow takes $Y_{i}\cap H_{i}\cap f^{-1}(\dot{D}_{\delta_{i+1}})$ to
$Y_{i+1}\cap H_{i+1}$ in time $T_i=1$.

Let $\rho_1:D_\delta\to [0,1)$ be a smooth function vanishing at $0$ and positive in $D_\delta^*$.
Let $\rho_2:U\to\RR$ be an smooth function with support contained in $\cup_{i\in\NN} Z_i$
and which is identically $1$ in a neighborhood of 
$f^{-1}(0)\cap (\cup_{i\in\NN} Z_i)$ in $\cup_{i\in\NN} Z_i$. Define the vector field 
$\calX:=(\rho_1\comp f)\calY+\rho_2\calZ$ on $Y$. 
The first three properties that $\calZ$ must satisfy are clear by construction (for the second one we take $\calW=\rho_1\calV'$).

For the fourth observe that as $\sum_{i=1}^{\infty}T_i=\infty$ 
the convergence to the origin in infinite time is clear for any curve contained
in $f^{-1}(0)$. Let $\gamma(t)$ be an integral curve lying off $f^{-1}(0)$. By property (2) and the fact that $\calX$ is tangent to
the fibres of $\pi$ (which are compact) it must accumulate in infinite time to a
point in $f^{-1}(0)$, but since $\calX$ is defined and non-zero at $f^{-1}(0)\setminus E_0$, the only possibility is that $\gamma$
satisfies property (4).
\end{proof}

\noindent\textbf{Notation}. Let $H_1$ $\theta_1$, $Y_1$, $\delta_1$ as in the beginning of the section. Denote the cut 
$H_1$ by $H$, the space $Y_1$ by $Y$, and $\delta_1$ by $\delta$ for simplicity of notation. 
Consider the natural projections $\sigma_1:L\times C\to L$, $\sigma_2:L\times C\to C$.
We view $L$ as a subset of $L\times C$ via its natural identification with $L\times\{(0,...,0)\}$.
View $L\times C$ as a subset of $Y$ identifying it with the zero section of $E|_{L\times C}$.
Given a subset $A\subset E$, we define $\partial_\pi A:=\cup_{t\in U}\partial A_t$, that is, 
the union of the frontiers of the fibres over $U$.

The globalisation of Timourian's~\cite{Ti} and King's~\cite{Ki2} result that we will need is the 
following:

\begin{theo}
\label{globalisation} 
There exist a homeomorphism
\begin{equation}
\label{homeoisolated}
\Psi:Y\to Y_L\times C
\end{equation}
such that
\begin{enumerate}
\item We have the equality $\pi=(\sigma_1\comp\pi|_{Y_L}\comp p_1,p_2)\comp\Psi$, where $p_i$ is the projection of $Y_L\times C$ 
to the $i$-th factor ($i=1,2$).
\item We have the equality $f|_{Y_L}\comp p_1\comp\Psi=f$.
\item The restriction of $\Psi$ to $Y\setminus (L\times C)$ is smooth.
\end{enumerate}
\end{theo}
\begin{proof}
By Ehresmann fibration Theorem the mappings 
\begin{equation}
\label{lafibracion}
\phi|_{Y\cap f^{-1}(\partial D_\delta)}:Y\cap f^{-1}(\partial D_\delta)\to U\times \partial D_{\delta},
\end{equation}
\begin{equation}
\label{lafrontera}
\phi|_{H\cap Y}:H\cap Y\to U\times D_{\delta}
\end{equation}
are locally trivial fibrations, (whose fibre are diffeomorphic respectively to the Milnor fibre of $f_t$,
and the abstract link of $f_t$ for any $t\in U$).

Observe that we have  
\[\partial_\pi Y=(Y\cap f^{-1}(\partial D_\delta))\cup (H\cap Y),\]
\[\partial_\pi (Y\cap f^{-1}(\partial D_\delta))=\partial_\pi (H\cap Y)=(Y\cap f^{-1}(\partial D_\delta))\cap (H\cap Y)=Y\cap f^{-1}(\partial D_\delta)\cap H.\]
Using the fibrations above and the fact that $C$ is a cube (and hence contractible) 
we will construct a diffeomorphism
\begin{equation}
\psi:\partial_\pi Y \to C\times \partial_{\pi}Y_L
\end{equation}
such that, if let $\psi$ play the role of $\Psi$, properties (1) and (2) are true 
at the points where $\psi$ is defined. We define first
the restriction of $\psi$ to $Y\cap f^{-1}(\partial D_\delta)$, and after 
extend it to $H\cap Y$.

Lift 
a radial vector field $\calZ_1$ steming from $(0,...,0)$ in $C$, first to a vector field $\calZ_2$ in 
$L\times C$ whose projection to $L$ vanishes, 
and then to a vector field $\calZ_3$ in $Y\cap f^{-1}(\partial D_\delta)$ which is tangent to the 
fibres of $f$. The existence of the lifting
$\calZ_2$ is obvious, and the existence of the lifting $\calZ_3$ follows because the 
mapping~(\ref{lafibracion}) is locally trivial. The inverse of the restriction of $\psi$ to 
$Y_1\cap f^{-1}(\partial D_\delta)$ can be obtained easily from the flow of $\calZ_3$.

Now we extend $\psi$ to $H\cap Y$. We consider a lift $\calZ_4$ of $\calZ_2$ in $H\cap Y$ which is tangent to the 
fibres of $f$, and which coincides with $\calZ_3$ at their common domain $Y\cap f^{-1}(\partial D_\delta)\cap H$
(this lift exists by the
local triviality of the mapping~(\ref{lafrontera})). The inverse of the desired
extension is constructed easily from the flow of $\calZ_4$.

In order to obtain $\Psi$ we have to extend $\psi$ to the interior of $Y$.

The set $Z:=\partial_\pi Y\times [0,\infty)$
is a manifold with corners 
The properties of the vector field $\calX$ constructed in Lemma~\ref{vectorfield} easily imply that,
if $\varphi$ denotes its flow, then the restriction
\[\varphi|_Z:Z\to Y\setminus L\times C\]
is a diffeomorphism.

Define the restriction $\Psi|_{Y\setminus L\times C}$ as the composition of the following sequence of mappings
\[Y\setminus L\times C\stackrel{\varphi|_Z^{-1}}{\longrightarrow}\partial_\pi Y\times [0,\infty)\stackrel{(\psi ,id_{[0,\infty))}}{\longrightarrow} C\times\partial_\pi Y_L\times [0,\infty)\stackrel{(id_C,\varphi|_Z)}{\longrightarrow} C\times Y_L,\]
and the restriction $\Psi|_{L\times C}$ as 
\[\Psi|_{L\times C}:=Id_{L\times C}.\]

By construction each of the restrictions $\Psi|_{Y\setminus L\times C}$ and $\Psi|_{L\times C}:=Id_{L\times C}$ is smooth. 
The mapping $\Psi$ is continuous at $L\times C$ because $\psi$ satisfies
Property (1) of the statement.
Property (1) of the statement is also 
satisfied by $\Psi$ because $\psi$ satisfies it and $\calX$ satisfies Property~(1) of
Lemma~\ref{vectorfield}. Property (2) of the statement is satisfied by $\Psi$ because $\psi$ satisfies it and $\calX$ 
satisfies Property~(2) of Lemma~\ref{vectorfield}.
\end{proof}

\begin{cor}
Let $f:C\times (\CC^n,O)\to\CC$ be a smooth family of holomorphic germs defining isolated singularities at the origin with constant 
Milnor number. Let $p$ be any base point in $C$. Define the space $Y$ as in Theorem~\ref{globalisation}. There exists a 
homeomorphism
\[\Psi:Y\to C\times Y|_p\]
satisfying properties $(1-3)$ of Theorem~\ref{globalisation}.
\end{cor}
\begin{proof}
We only have to notice that Assumption A is satisfied (Remark~\ref{unpunto}).
\end{proof}

\begin{remark}
The above result improves the results of King and Timourian (\cite{Ki2},\cite{Ti}) in the sense that it is global in the base 
(we do not need to pick a small neighborhood of $p\in C$ in order to find the trivialisation), and also because the trivialisation 
$\Psi$ is smooth at all the strata where it can possibly be (smoothness at the origin can not be expected
since isolated singularities have non-trivial moduli).
\end{remark}

\subsection{Families over Riemann surfaces}

Now we prove Assumption A in an important situation in which $L=\SSS^1$.

Let $S$ be a compact 
surface diffeomorphic to $\SSS^1\times I$ embedded as a smooth submanifold in 
a Riemann surface $\tilde{S}$. As we are only interested in a neighborhood of $S$ we can assume that
$\tilde{S}$ is diffeomorphic to $\SSS^1\times (-\epsilon,1+\epsilon)$ for a certain positive $\epsilon$; we 
fix a product decomposition and denote by
\[\alpha:\tilde{C}\to(-\epsilon,1+\epsilon)\]
the projection to
the second factor. We suppose that $U=\tilde{C}$ and that the vector bundle $\pi:E\to U$ and the function $f$ are 
{\em holomorphic}. Recall that we have defined $\phi:=(\pi,f)$.

\begin{prop}
\label{preassumptionA}
There is a finite subset $J\subset I$ such that,
for any $p\in I\setminus J$, there exists $\zeta>0$ satisfying that
$\rho^{-1}(\epsilon)$ meets $\varphi^{-1}(t,0)$ transversely for any positive $\epsilon\leq\zeta$
and any $t\in\alpha^{-1}(p)$, In other words, up to a finite number of exceptional $p\in I$, 
assumption A is satisfied over $\alpha^{-1}(p)=\SSS^1\times\{p\}$. 
\end{prop}
\begin{proof}
Recall that the critical set of the function $f$ coincides with the zero section 
of $\pi$ and the restrictions of $f$ to the fibres of the bundle $\pi$ have 
an isolated singularity at the origin with Milnor number not depending on the particular fibre. 

By the holomorphicity of $\pi$ and $f$, the locus $Z$ in which the stratification 
\begin{equation}
\label{estrati}
\{f^{-1}(0)\setminus \tilde{C},\tilde{C}\}
\end{equation}
fails to satisfy Whitney conditions is the complex analytic subset where the $\mu^*$ sequence jumps.
The set $Z$ is of dimension $0$, and hence discrete. Thus $Z\cap S$ is finite.
Define $J:=\alpha(Z)\cap I$.

Suppose $p\in I\setminus J$. If the statement of the proposition fails there exists a sequence 
$\{x_n\}_{n\in\NN}$ in $E_{\alpha^{-1}(p)}$, converging to a 
point $x$ in the zero section over $\alpha^{-1}(p)$ such that $x_n$ is a critical point for the 
restriction of $\rho$ to $f^{-1}(0)\cap E_{\pi(x_n)}$. Then it is easy to check that the tangent 
space $T_{x_n}f^{-1}(0)$ of $f^{-1}(0)$ at $x_n$ contains the orthogonal complement 
(by the hermitian inner product) of the line $\overline{\pi(x_n),x_n}$ in $E_{\pi(x_n)}$. 
Choosing a subsequence we may assume that $T_{x_n}f^{-1}(0)$ converges to a complex hyperplane $H$ of 
$T_x E$ and that $\overline{\pi(x_n),x_n}$ converges to a line $l$ contained in $\pi^{-1}(\pi(x))$, 
which satisfies that $H$ contains the hermitian orthogonal complement of $l$ in $\pi^{-1}(\pi(x))$. 
As $p$ does not belong to $J$, the stratification~(\ref{estrati}) satisfies Whitney conditions at $x$, 
and therefore we have that $l$ is also included in $H$. Then $H$ is in fact equal (by dimensional reasons) to $\pi^{-1}(\pi(x))$.

Now we choose a sequence $\{y_n\}_{n\in\NN}$, with $f(y_n)\neq 0$ and $y_n$ so close to 
$x_n$ that the sequence of hyperplanes $T_{y_n}f^{-1}(f(y_n))$ has limit $H$. Such sequence 
shows that Thom's $A_{f}$-condition is not satisfied at $x$ by the open stratum of the 
stratification 
\begin{equation}
\label{estrati2}
\{E\setminus f^{-1}(0),f^{-1}(0)\setminus\tilde{C},\tilde{C}\},
\end{equation}
which is in contradiction with the fact that the transversal Milnor number of $f$ is constant along 
$\tilde{C}$.
\end{proof}

We recall our original setting: consider a smooth manifold $U=\SSS^1\times (-\epsilon,\epsilon)^d$, 
a smooth complex vector bundle $\pi:E\to U$ and a complex smooth function
$f:E\to\CC$ whose restriction to the fibres of $\pi$ is holomorphic. Suppose that $f$ has an isolated singularity at the origin of
each fibre of $\pi$, with Milnor number independent of the fibre.

The next corollary allows to apply the topological equisingularity results obtained in this section in many situations.
\begin{cor}
\label{assumptionA}
Suppose $L=\SSS^1$.
Assume that there is a smoothly embedded path $\gamma:(-\zeta,\zeta)\hookrightarrow [0,1]^n$ (with image
$B:=\gamma(-\zeta,\zeta)$)
such that $\tilde{S}:=\SSS^1\times B$ has a structure of Riemann surface compatible with the smooth structure, and that 
$E_{\tilde{S}}$ has a structure of holomorphic vector bundle which makes the restriction of $f|_{E_{\tilde{S}}}$ holomorphic.
Then, for a generic point $p\in B$, Assumption A is satisfied if we let $\SSS^1\times\{p\}$ play the role of $\SSS^1\times\{(0,...,0)\}$.
Therefore the statement of Proposition~\ref{usefulconstr} holds letting $\SSS^1\times\{p\}$ play the role of
$\SSS^1\times\{(0,...,0)\}$ 
(Remark~\ref{basepoint}). Hence also all the topological equisingularity results stated above hold for the family $f$. 
\end{cor}

\section{Adapted neighborhoods}

Let $f:(\CC^n,O)\to\CC$ be any holomorphic germ. For the homotopical study of the 
Milnor fibration of $f$, different classes
of systems of neighborhoods of the origin are possible (see the definition of {\em system of Milnor neighborhoods} 
in~\cite{Ma3}~page~105). 

For the topological study of the Milnor fibration and of the embedded link of a function one has to be
more careful.
In~\cite{Lo}~pages~21-25 it is proved that if $\rho:\CC^n\to\RR$ is non-negative continuous, 
real analytic except at the origin, and such that $\rho^{-1}(0)=\{O\}$ (we will refer to it as an 
analytic distance function), then, for $\epsilon$ small enough, the system of $\rho$-balls 
$\{\rho^{-1}([0,\epsilon)\}$ is a system of Milnor neighborhoods for $f$. Moreover it is shown that 
the homeomorphism type of the Milnor fibration and the embedded link of $f$ is independent of the 
chosen analytic distance function. We take as definition up to homeomorphism of 
{\em Milnor fibration and embedded link} of $f$ the ones obtained with an analytic distance function. 
This choice coincides with Milnor's original definition~(see \cite{Mi}).    

We will need a kind of systems of neighborhoods that can not be defined with analytic distance 
functions. We will have to check that the {\em Milnor fibration} and the {\em embedded link} defined with this kind
of neighborhoods are, up to homeomorphisms, the same than the ones defined above.

\begin{definition}
\label{normbund}
Let $M$ be a manifold and $\Sigma\subset M\times\CC^n$ be a smooth submanifold such that the restriction
of the projection $\pi:M\times\CC^n\to M$ to $\Sigma$ is a submersion and for any $m\in M$ the fibre 
$\Sigma_m\subset \{m\}\times\CC^n$ is a complex submanifold of dimension $1$ (a Riemann surface).
Denote by $\calT(M\times\CC^n,M)\to M\times\CC^n$ the relative tangent bundle of $\pi$.
Denote by $\calT(\Sigma,M)\to\Sigma$ the relative tangent bundle of the submersion $\pi|_\Sigma$.
Let $S\subset\Sigma$ be a submanifold. A {\em normal bundle} to 
$\Sigma$ over $S$ is a complex vector subbundle $p:N\to S$ of rank $n-1$ 
of $\calT(M\times\CC^n,M)|_S$ such that $N+\calT(\Sigma,M)|_S=\calT(M\times\CC^n,M)|_S$ , 
(that is, its fibre $N_p$ is transverse to $\calT(\Sigma,M)_p$ for any $p\in S$).
\end{definition}
We endow $N$ with
the hermitian metric coming from a hermitian product in $M\times\CC^n$. If $S$ is compact there exists a positive $\eta$ such 
that the tubular neighborhood
\[T_N(\eta):=B(S,\eta)\subset N\]
of the zero section of $N$ is naturally embedded in $M\times \CC^n$.

\noindent\textbf{Assumption}. Assume from this point that the critical set $\Sigma_f$ of $f$ has dimension at most $1$. 
Denote by $\Sigma_f^{sm}$ the smooth part of $\Sigma_f$.

Observe that $f^{-1}(0)$ is a stratified space with the stratification 
$(f^{-1}(0)\setminus\Sigma_f,\Sigma_f^{sm},\{O\})$.

\begin{definition}
\label{goodneigh}
A {\em system of neighborhoods of the origin adapted to $f$} is a system $\{U_\alpha\}_{\alpha\in A}$ of neighborhoods of the
origin such that for any $\alpha\in A$ there exists $\delta_\alpha>0$ 
(which is called a {\em Milnor disk radius for} $U_\alpha$), satisfying   
\begin{enumerate}
\item the neighborhood $U_\alpha$ is a compact manifold with boundary $\partial U_\alpha$.
\item There is a tubular neighborhood $W$ of $\Sigma_f\cap\partial U_\alpha$ in $\partial U_\alpha$ such that $W$ coincides
with $T_{N}(\eta)$ for a certain normal bundle to $\Sigma_f^{sm}$ called 
\[p:=N\to \Sigma_f\cap\partial U_\alpha,\]
and a certain $\eta>0$.
\item The function $f$ has no critical points in $U_\alpha\cap f^{-1}(D^*_{\delta_\alpha})$, and 
for any $s\in D_{\delta_\alpha}$ the boundary $\partial U_\alpha$ meets $f^{-1}(s)$ transversely 
(in the stratified sense for $s=0$).
\item For any $\beta$ such that $U_\beta\subset \dot{U}_\alpha$, and any
$\delta<\min\{\delta_\alpha,\delta_\beta\}$ there is a homeomorphism
\begin{equation}
\kappa:(U_\alpha\setminus\dot{U}_\beta)\cap f^{-1}(D_\delta)\to\partial U_\beta\cap f^{-1}(D_\delta)\times [0,1]
\end{equation}
whose restrictions to $((U_\alpha\setminus\dot{U}_\beta)\cap f^{-1}(D_\delta))\setminus \Sigma_f)$ and to
 $(U_\alpha\setminus\dot{U}_\beta)\cap\Sigma_f$ are smooth
and which satisfies $f=f\comp p_1\comp\kappa$, 
where $p_1$ is the projection of $\partial U_\beta\cap f^{-1}(D_\delta)\times [0,1]$ to the first factor.
\item For $\delta>0$ small enough the space $U_\alpha\cap f^{-1}(D_\delta)$ is contractible.
\end{enumerate}
\end{definition}

\begin{remark}
The fourth condition of the last definition implies that the space
\begin{equation}
f^{-1}(s)\cap U_{\alpha}\setminus\dot{U}_\beta
\end{equation}
is a trivial (stratified when $s=0$) cobordism for any $s\in D_\delta$.
\end{remark}

Let $\{U_\alpha\}_{\alpha\in A}$ be a system of neighborhoods of the origin adapted to $f$.
By Ehresmann Fibration Theorem,
for any smooth neighborhood $U_\alpha$ of the origin adapted to $f$, the mapping 
\begin{equation}
\label{fibracmilnor}
f:U_\alpha\cap f^{-1}(D^*_{\delta_\alpha})\to D^*_{\delta_\alpha} 
\end{equation}
is locally trivial. Its restriction over $\partial D_{\delta_\alpha}$ is called the 
{\em Milnor Fibration} of $f$ for $U_\alpha$. 

We define the {\em embedded link} of $f$ in $U_\alpha$ as the pair
\begin{equation}
\label{emblink}
\big{(}\partial (U_\alpha\cap f^{-1}(D_\alpha)),f^{-1}(0)\cap\partial (U_\alpha\cap f^{-1}(D_\alpha))\big{)}.
\end{equation}
The {\em abstract link} is, by definition, the space $f^{-1}(0)\cap\partial (U_\alpha\cap f^{-1}(D_\alpha))$.

Property (4) of Definition~\ref{goodneigh} implies that the topological type of the embedded link, 
and the diffeomorphism type of the Milnor fibration are independent of $\alpha\in A$.
Hence the the {\em Milnor fibration and the embedded link associated to a system of neighborhoods adapted to} $f$ 
are well defined objects. 

\begin{lema}
\label{homotopylink}
The homotopy type of the abstract link and of the Milnor fibre is independent on the system of neighborhoods adapted to $f$ 
used to define it. Any system of neighborhoods adapted to $f$ is a system of Milnor neighborhoods for $f$ 
(\cite{Ma3}, Definition/Proposition A.1).  
\end{lema}
\begin{proof}
Suppose that we have two systems of neighborhoods $\{U_\alpha\}_{\alpha\in A}$ and $\{V_\beta\}_{\beta\in B}$ adapted to $f$.
Consider sequences $\{\alpha_n\}_{n\in\NN}$, $\{\beta_m\}_{m\in\NN}$ such that $U_{\alpha_n}\subset\dot{V}_{\beta_n}$ and 
$V_{\beta_{n+1}}\subset\dot{U}_{\alpha_{n}}$. By definition of adapted system of neighborhoods, for any $n$ there exists $\delta$ 
such that the cobordisms $f^{-1}(s)\cap (U_\alpha\setminus\dot{U}_{\alpha_{n+1}})$ and $f^{-1}(s)\cap (V_\beta\setminus\dot{U}_{\beta_{n+1}})$ are trivial for any $s\in D_\delta$. It is easy to deduce that the cobordisms 
\[f^{-1}(s)\cap (V_{\beta_n}\setminus\dot{U}_{\alpha_n})\]
are homotopically trivial for any sufficiently small $s$. This implies the lemma.
\end{proof}

We will prove:
\begin{itemize}
\item the Milnor fibration of $f$ for $\{U_\alpha\}_{\alpha\in A}$ is $C^\infty$-equivalent with the classical Milnor fibration of
 $f$ defined by Milnor.
\item the topological type of the embedded link associated to any neighborhood adapted to $f$
coincides with the classical embedded link defined by Milnor.
\end{itemize}
In the remaining part of the section we exhibit a particular system of neighborhoods of the origin 
adapted to $f$ for which these assertions are true. After this, it is enough to show that the Milnor
fibration and the embedded link are independent of the neighborhood of the origin 
adapted to $f$. This will be achieved in the next sections.

We will use the following Lemmas several times in the sequel.

\begin{lema}
\label{gentranstype}
Let $f:(\CC^n,x)\to\CC$ be a holomorphic germ with smooth $1$-dimensional critical locus $\Sigma_f$. Suppose that the generic 
transversal Milnor number of $f$ at any point $y\in\Sigma_f$ is equal to a constant $\mu$. If $H$ is a hyperplane transverse 
to $\Sigma_f$ through $x$ then $\mu(f|_H,x)=\mu$.
\end{lema}
\begin{proof}
The L\^e numbers of $f$ with respect to a generic coordinate system $\calZ_0$ are $\lambda^1_{f,\calZ_0}(y)=\mu$ and 
$\lambda^0_{f,\calZ_0}(y)=0$ for any 
$y\in\Sigma_f$. As $H$ is transverse to $\Sigma_f$ there is another prepolar coordinate system, $\calZ_1:=\{z_1,...,z_n\}$, 
with $V(z_1)=H$. The first L\^e number $\lambda^1_{f,\calZ_1}(x)$ is equal to the Milnor number $\mu'$ of $f|_{z_1^{-1}(z_1(y))}$ at 
$y$, for $y\in\Sigma_f$ generic. We claim that $\mu'=\mu$. 

As $\mu$ is the generic transversal Milnor number we have $\mu'\geq\mu$.
By definition of the first L\^e number we have the equality 
$\lambda^1_{f,\calZ_1}(x)=\lambda^1_{f,\calZ_1}(y)$ for any $y\in\Sigma_f$. Hence 
as the alternating sum of the L\^e numbers equals the reduced Euler characteristic of the Milnor fibre, which is equal to $\mu$,
if $\mu'>\mu$ then $\lambda^0_{f,\calZ_1}(y)\neq 0$. 
On the other hand, the locus where $\lambda^0_{f,\calZ_1}(y)$
is different from zero is the intersection locus of the relative polar curve $\Gamma^1_{f,\calZ_1}$ with $\Sigma_f$. As, by definition
of $\Gamma^1_{f,\calZ_1}$, it can not contain $\Sigma_f$ we have that $\lambda^1_{f,\calZ_1}(y)=0$ at most points. Hence the claim
is true.

As $\mu$ is the generic transversal Milnor number we have $\mu(f|_H,x)\geq\mu$. If the inequality is strict then the relative polar
curve $\Gamma^1_{f,\calZ_1}$ meets $\Sigma_f$ at $x$, and hence $\lambda^0_{f,\calZ_1}(x)>0$. Then, as $\lambda^1_{f,\calZ_1}(x)=\mu$,
the alternating sum of the L\^e numbers cannot be equal to the reduced Euler characteristic $\mu$. We get a contradiction which forces
the equality $\mu(f|_H,x)=\mu$.
\end{proof}

\begin{remark}
Here we sketch an altenate proof of the previous Lemma suggested by T. Gaffney: let $H_t$ be a family of
generic hyperplanes, so $f|_{H_t}$ is a $\mu$-constant family. Then Thom's $A_f$ condition holds, and
therefore no hyperplane transverse to the critical line can be a limit of tangent hyperplanes. A result 
of L\^e-Henry states that if $H$ is not a limit of tangent hyperplanes then $f|_H$ has the generic 
transversal Milnor number.
\end{remark} 

\begin{lema}
\label{normalbundle}
Let $M$ be a manifold and $\Sigma\subset M\times\CC^n$ be a smooth submanifold such that $\pi|_\Sigma$ is a submersion and,
for any $m\in M$, the fibre 
$\Sigma_m\subset \{m\}\times\CC^n$ is a complex submanifold of dimension $1$. Consider a compact submanifold 
$S\subset\Sigma$ and a normal bundle $p:N\to S$ to $\Sigma$ over $S$. There exists a normal bundle $p':N'\to\Sigma$ to 
$\Sigma$ whose restriction to $S$ coincides with $p:N\to S$.
\end{lema}
\begin{proof}
Let $\sigma:\PP\to M\times\CC^n$ be the projectivised dual of $\calT(M\times\CC^n,M)\to M\times\CC^n$. 
An element of $\PP$ is a pair $(H,x)$ where $x\in M\times\CC^n$ and $H$ is a hyperplane of $\calT(M\times\CC^n,M)_x$. Define 
$Z\subset\PP|_\Sigma$ as 
\[Z:=\{(H,x)\in \PP|_\Sigma:H\supset\calT(\Sigma,M)_x\}.\]

As $\calT(\Sigma,M)_x$ is a line in $\calT(M\times\CC^n,M)_x$ the fibre $Z_x$ is a hyperplane of the projective space $\PP_x$.
Hence the restriction 
\begin{equation}
\label{aqui}
\sigma_{\PP|_\Sigma\setminus Z}:\PP|_\Sigma\setminus Z\to \Sigma
\end{equation}
is an affine bundle. Any affine bundle admits a differentiable global section, and hence can be given a structure of vector bundle. 

Observe that giving a normal bundle to $\Sigma$ over a submanifold $S$ of $\Sigma$ is the same that giving a smooth section
of the vector bundle~(\ref{aqui}) over $S$. Using a partition of unity it is possible to extend the section smoothly to a global one.
\end{proof}

\subsection{Construction of adapted neighbourhoods}.
Consider any hermitian product in $\CC^n$ and let $\rho$ be the distance function to the origin induced 
by it. Denote by $B_\epsilon$ the $\rho$-ball centered at the origin of radius $\epsilon$. 
Choose $\epsilon_0$ sufficiently small so that 
\begin{enumerate}
\item the stratification 
\[\{B_{\epsilon_0}\setminus f^{-1}(0),B_{\epsilon_0}\cap f^{-1}(0)\setminus\Sigma_f,B_{\epsilon_0}\cap\Sigma_f\setminus\{O\},\{O\}\}\]
satisfies Whitney conditions,
\item the sphere $\SSS_\epsilon$ meets $f^{-1}(0)$ transversely (in a stratified sense) for any 
$\epsilon<\epsilon_0$,
\end{enumerate}
It is well known that the system of neighborhoods $\{B_{\epsilon}\}_{0<\epsilon<\epsilon_0}$ satisfies
all the properties of a system of neighborhoods of the origin adapted to $f$ except the second one.
We are going to modify the neighborhoods so that this property is satisfied as well.

Fix $0<\epsilon<\epsilon_0$. The intersection $\Sigma_f\cap\SSS_\epsilon$ is a disjoin union of circles.
We will modify $\SSS_\epsilon$ locally near each circle. We assume for simplicity of notation that there is 
only one circle $S$. For $\xi>0$ small enough 
$A:=\Sigma_f\cap B_{\epsilon+\xi}\setminus\dot{B}_{\epsilon-\xi}$ is diffeomorphic to 
$S\times [-\xi,\xi]$ by a diffeomorphism whose second component in the product $S\times [-\xi,\xi]$ is the function $\rho-\epsilon$. We identify $A$ and $S\times [-\xi,\xi]$ via the diffeomorphism.

I claim that there is a normal bundle $p:N\to A$ to $A$ such that 
\begin{enumerate}[(i)]
\item for any $s\in S\times [-\xi/3,\xi/3]$ the fibre $N_s$ is contained in the tangent space of the level set of $\rho$ at
 $s$,
\item it is holomorphic over $S\times [2\xi/3,\xi]$.
\end{enumerate}
There is a normal bundle $p_1:N_1\to S\times [-\xi/3,\xi/3]$ to $T$ satisfying the first property:
as $T$ meets the level sets of $\rho$ transversely, if $\xi$ is small enough,
for any $s\in S\times [-\xi/3,\xi/3]$ the hermitian orthogonal complement to the vector $s\in\CC^n$ is transverse
to $T_sA$ (we are using the natural trivialisation of the tangent bundle $\calT\CC^n$) and contained in $T\SSS_{\rho(s)}$.
The fibres of $N_1$ are, by definition, those hermitian orthogonal complements.

Choose a coordinate system $\{z_1,...,z_n\}$ of $\CC^n$ such that, for any $x$ in a punctured neighborhood of the origin, 
the hyperplane $V(z_1-z_1(x))$ is transverse to $T_x \Sigma_f$. Then the subbundle 
\begin{equation}
\label{normhol}
p:N_{z_1}\to\Sigma_f\setminus\{O\}
\end{equation}
whose fibre over $x$ is $V(z_1-z_1(x))$ is a normal bundle to $\Sigma_f$ in the same punctured neighborhood.
The bundle $p_2:N_2\to S\times [2\xi/3,\xi]$ is defined as the restriction of the bundle~(\ref{normhol}) over $S\times [\xi/2,\xi]$.

Any normal bundle $p:N\to A$ to $A$ expending $N_1$ and $N_2$ (which exists by Lemma~\ref{normalbundle}) satisfies requirements
$(i)$ and $(ii)$. The claim is proved.

Pick $\eta>0$ so small that $T_N(\eta)$ is embedded in $\CC^n$ and $T_{N_S}(\eta)$ is 
contained in $B_{\epsilon+\xi}$. The restriction $f|_N$ is a smooth function, whose restriction to $N_{S\times [2\xi/3,\xi]}$
and to the fibres of $p$ is holomorphic. As a consequence of condition
(1) above, the generic transversal Milnor number of $f$ at any point $t\in A$ is independent on the point 
$t$. Hence, by Lemma~\ref{gentranstype}, the restriction of $f$ to the fibres of $p$ has an isolated 
singularity at the origin with Milnor number independent of the fibre. As we are in the 
situation of Corollary~\ref{assumptionA}, we can use 
Theorem~\ref{globalisation} in order to obtain (after possibly shrinking $\eta$) a continuous flow
\begin{equation}
\label{homeoisol}
\Phi:T_{N_S}(\eta)\times [-\xi,\xi]\to\CC^n
\end{equation}
whose restrictions to the zero section $S\times [-\xi,\xi]$ and its complement are smooth, such that
\begin{enumerate}[(a)]
\item for any $s\in S$ we have 
$\rho\comp\Phi(s,t)=\epsilon+t$,
\item the flow lines 
are contained in the fibres of $f$ and are compatible with the bundle projection $p$ (the projection by 
$p$ of a flow line is a flow line); the flow takes $S\times [-\xi,\xi]$ into $A$.
\end{enumerate}

Conditions (i) and (a) above, and the compatibility of the flow lines with
$p$, imply that, if $\eta$ is chosen small enough, there exists a positive small enough 
$\xi'$ such that any flow line of the restriction
\begin{equation}
\label{homeoisol'}
\Phi':=\Phi|_{T_{N_S}(\eta)\times [-\xi',\xi']}:T_{N_S}(\eta)\times [-\xi',\xi']\to\CC^n
\end{equation}
is transverse to the level spheres of $\rho$ and points outwards. Hence for any $x\in T_{N_S}(\eta)$ there exists at most one value 
$\alpha(x)\in (-\xi',\xi')$ such that $\Phi'(x,\alpha(x))$ belongs to $\SSS_\epsilon$. 
Since $S$ is contained in $\SSS_\epsilon$, we have that
$\alpha(x)$ is zero for any $x\in S$. We have defined a function
\[\alpha:T_{N_S}(\eta)\to\RR\]
which is continuous, constant and equal to $0$ in $S$, and smooth outside $S$. 

The mapping $\Psi:T_{N_S}(\eta)\to\CC^n$ defined by the formula $\Psi(x)=\Phi'(x,\alpha(x))$ is a homeomorphism 
of $T_{N_S}(\eta)$ onto its image
$W$, a neighborhood of $S$ in $\SSS_\epsilon$. Consider a smooth function $\mu:\SSS_\epsilon\to [0,1]$, with support $K$ contained
in $W$, and identically $1$ in a neighborhood of $S$. We deform $\SSS_\epsilon$ as follows: we define the homotopy
\[\varphi:\SSS_\epsilon\times [0,1]\to\CC^n\]
piecewise by $\varphi(x,\lambda):=\Phi'(x,-\lambda\mu(x)\alpha(\Psi^{-1}(x)))$ if $x\in W$ and $\varphi(x,\lambda):=x$ if $x$ is not 
contained in $K$. For fixed $\lambda$ the mapping $\varphi(\centerdot,\lambda)$ is a homemorphism onto its image and coincides 
with $\Id_{\SSS_\epsilon}$ if $\lambda=0$. We denote by $U_{\epsilon,\lambda}$ the region of $\CC^n$ bounded by
$\varphi(\SSS_{\epsilon},\lambda)$.

For any $\lambda\in \{0,1\}$ the neighborhood $U_{\epsilon,\lambda}$ satisfies Properties (1) and (3) of Definition~\ref{goodneigh},
and Property (2) is satisfied for $\lambda=1$.

Using Property (b) of the flow~(\ref{homeoisol}) and the fact that Property (4) of Definition~\ref{goodneigh} is satisfied by the system
$\{B_\epsilon\}_{0<\epsilon<\epsilon_0}$ it is easy to show that
for any $\epsilon'<\epsilon$ there exists $\delta>0$ such that there is a homeomorphism
\begin{equation}
\label{(22)}
\kappa:(U_{\epsilon,1}\setminus\dot{B}_{\epsilon'})\cap f^{-1}(D_\delta)\to(\SSS_{\epsilon'}\cap f^{-1}(D_\delta))\times [0,1]
\end{equation}
whose restriction to $((U_{\epsilon,1}\setminus\dot{B}_{\epsilon'})\cap f^{-1}(D_\delta))\setminus \Sigma_f)$ and to
$(U_{\epsilon,1}\setminus\dot{B}_{\epsilon'})\cap\Sigma_f$ is smooth and satisfies $f=f\comp p_1\comp\kappa$, where $p_1$ is the 
projection of $(\SSS_{\epsilon'}\cap f^{-1}(D_\delta))\times [0,1]$ to the first factor. In particular
\begin{itemize}
\item for any $\epsilon'<\epsilon$, there exists $\delta>0$ 
such that for any $s\in D_\delta$ the space 
\begin{equation}
\label{auxtrivcob}
f^{-1}(s)\cap U_{\epsilon,1}\setminus\dot{B}_{\epsilon'}
\end{equation}
is a trivial (stratified when $s=0$) cobordism,
\item we have a homeomorphism
\begin{equation}
\label{otrohomeo}
(\partial (B_{\epsilon'}\cap f^{-1}(D_\delta)),f^{-1}(0)\cap \partial (B_{\epsilon'}\cap f^{-1}(D_\delta)))\cong\partial (U_{\epsilon,1}\cap f^{-1}(D_\delta)),f^{-1}(0)\cap \partial (U_{\epsilon,1}\cap f^{-1}(D_\delta))).
\end{equation}
\end{itemize}

Property (5) of Definition~\ref{goodneigh} for $U_{\epsilon,\lambda}$ is a consequence of property (b) of the flow~\ref{homeoisol}
and the fact that for any $\epsilon'<\epsilon$ the ball $B_{\epsilon'}$ satisfies it.

We construct a system of neighborhoods as follows. Give a decreasing sequence of positive numbers $\{\epsilon_n\}_{n\in\NN}$ 
converging to the $0$. For any $n$ we produce a deformation $U_{\epsilon_n,\lambda}$ of $B_{\epsilon_n}$ as above. The choices can 
be made so that $U_{\epsilon_n,\lambda}$ is contained in $\dot{B}_{\epsilon_{n-1}}$ for any $n\in\NN$. 

\begin{prop}
\label{ejemplogoodneigh}
The family $\{U_{\epsilon_n,1}\}_{n\in\NN}$ is a system of neighborhoods adapted to $f$. Moreover the topological type of the 
embedded link associated to it, and the diffeomorphism type of the Milnor fibration associated 
with it coincides with the classical ones.
\end{prop}
\begin{proof}
Property (4) of Definition~\ref{goodneigh} (the only one that remains to be checked) follows easily using that it is 
satisfied for the system $\{B_{\epsilon_n}\}$ and that we have the homeomorphism~(\ref{(22)}).

The triviality of cobordism~(\ref{auxtrivcob}) implies also the diffeomorphisms of Milnor fibrations. The fact that the embedded
link associated to $\{U_{\epsilon_n,1}\}_{n\in\NN}$ is homeomorphic to the pair 
\[(\partial B_{\epsilon_n}\cap f^{-1}(D_\delta), f^{-1}(0)\cap (\partial B_{\epsilon_n}\cap f^{-1}(D_\delta)),\]
which, by a Theorem of~\cite{Mi} is homeomorphic to $(\SSS_\epsilon,f^{-1}(0)\cap\SSS_\epsilon)$ (the classical embedded link),
is essentially due to homeomorphism~(\ref{otrohomeo}).
\end{proof}

\section{Equisingularity at the critical set}
\label{eacs}

\noindent\textbf{Setting}. Let $T$ be a connected smooth manifold.
A smooth family of holomorphic germs $f_t:(\CC^n,O)\to\CC$ depending on a parameter $t\in T$, is given by a smooth function
$f:U\to \CC$ defined at an open neighborhood of $\{0\}\times T$ in $T\times\CC^n$
satisfying that the restriction $f_t:=F|_{\CC^n\times\{t\}}$ is holomorphic for any $t\in T$. We view $T\times\CC^n$ as a trivial
vector bundle over $T$ and consider the natural projection $\tau:T\times\CC^n\to T$. 
Let $\Sigma_{t}$ denote the critical set of $f_t$. Define $\Sigma:=\cup_{t\in T}\Sigma_t$.
Consider a hermitian product in $T\times\CC^n$
and denote by $\rho$ the distance to the origin in each fibre.

Now we define equisingularity at the critical set. Roughly speaking, the family  $f_t$ is equisingular at the critical set if
the family of critical sets $\Sigma_t$ is topologically equisingular by a family of homeomorphisms preserving the transversal
Milnor number (outside the origin). This is formalised as follows.

Given any $A\subset T$, $B\subset T\times\CC^n$ and a positive function 
$\theta:A\to\RR_+$ we define $B_A=B|_A$ and $B(A,\theta)$ as in the previous section.
Choose $t\in T$, there exists a positive $\epsilon_t$ such that $\partial B((t,O),\epsilon)$ meets $\Sigma_t$ transversely for any 
$\epsilon\leq\epsilon_t$. Then $\Sigma_t\cap B((t,O),\epsilon_t)$ is homeomorphic to the cone over 
$\Sigma_t\cap\partial B((t,O),\epsilon_t)$,
and there is an irreducible component of the germ $(\Sigma_t,(t,O))$ for each connected component of $\Sigma_t\cap\partial B((t,O),\epsilon_t)$.
Let $\Sigma_t=\cup_{i=1}^{r_i}\Sigma_{t,i}$ be the decomposition of $(\Sigma_t,(t,O))$ in irreducible components. The number $\epsilon_t$ 
can be chosen small enough that, for any $i\leq r$, the generic transversal Milnor number at any point of 
$\Sigma_{t,i}\cap B((t,O),\epsilon_t)\setminus\{(t,O)\}$ is equal to a constant $\mu_{t,i}$.

We will identify $T$ with $T\times\{O\}$.

\begin{definition}
\label{equiatsing}
Let $f_t$ be a smooth family of holomorphic germs. Suppose that $\Sigma_{t}$ is $1$-dimensional at the origin for any $t$.
We say that $f_t$ is equisingular at the critical set if the following conditions are satisfied:
\begin{enumerate}
\item The space $\Sigma\setminus T$ is smooth of real dimension $2+\mathrm{dim}(T)$ at any of its points, and such that 
the restriction $\tau|_{\Sigma\setminus T}$ is a submersion.
\item The previous property implies that for any $t\in T$, there exists a neighborhood $W_t$ of $t$ in $T$ so small that 
$\Sigma_{t'}$ meets $\partial B((t,O),\epsilon_t)$ transversely for any $t'\in W_t$. For any $t'\in W_t$ the space 
\[\Sigma_{t'}\cap B((t,O),\epsilon_t)\]
is homeomorphic to $\Sigma_{t}\cap B((t,O),\epsilon_t)$, (notice that if the homemorphism exists it can be chosen sending 
$(t,O)$ to $(t,O)$).
\item The previous property implies that for any $t,t'\in W_t$, there is a bijective correspondence (induced by inclusion) 
between the connected components of 
\[\Sigma_t\cap B((t,O),\epsilon_t)\setminus\{(t,O)\},\]
\[\Sigma_{W_t}\cap B(W_t,\epsilon_t)\setminus (\{O\}\times W_t)\] 
and
\[\Sigma_{t'}\cap B((t,O),\epsilon_t)\setminus\{(t,O)\}.\]
For any $t'\in W_t$ let $\Sigma_{t',i}$ be the connected component of 
$\Sigma_{t'}\cap B((t,O),\epsilon_t)\setminus\{(t,O)\}$ corresponding to $\Sigma_{t,i}$. 
The generic transversal Milnor number of $f_{t'}$ at any point of $\Sigma_{t',i}\setminus\{(t',O)\}$ is $\mu_{t,i}$.
\end{enumerate}
\end{definition}

\begin{remark}
\label{equiatsing'}
Condition $(2)$ of the previous Definition can be phrased informally saying that the underlying deformation of critical sets 
is topologically equisingular. It is important to notice that it implies that for any $t'$ close to $t$, the function $f_{t'}$ has 
no isolated singularities in $B((t,O),\epsilon_t)$. As $\Sigma_{t'}$ is of real dimension $2$ the condition can 
be rephrased in any of the following formulations:
the cobordism 
\[\Sigma_{t'}\cap B((t,O),\epsilon_t)\setminus\dot{B}((t,O),\epsilon_{t'})\]
\begin{itemize}
\item is trivial in the smooth category,
\item is homologically trivial.
\end{itemize}
\end{remark}

Let $V$ be any neighborhood of the origin adapted to $f_t$. 
A consequence is that $\partial V$ meets $\Sigma_t$ transversely 
for any $\alpha$. Choose $W_t$ so small that
$\Sigma_{t'}$ meets $\partial V$ transversely for any  $t'\in W_t$. Let 
$V'$ be any neighborhood of the origin adapted to $f_{t'}$ such that $V'\subset\dot{V}$.

\begin{remark}
By Property (4) of Definition~\ref{goodneigh}, Property (2) of Definition~\ref{equiatsing}, and an easy 
cobordism manipulation, the cobordism 
\begin{equation}
\label{coborcomp}
\Sigma_{t'}\cap V\setminus\dot{V}'
\end{equation} 
is trivial.
\end{remark}

Here it is a conormality Lemma (in the sense of~\cite{Ma3}, Definition A.7)

\begin{lema}
\label{conormality}
Let $f_t$ be a smooth family of holomorphic germs such that $\Sigma_{t}$ is $1$-dimensional at 
the origin for any $t$. Suppose that $f_t$ is equisingular at the critical set. Given $t\in T$ and
$V$, a smooth neighborhood of the origin adapted to $f_t$, there exists a 
neighborhood $W_t$ of $t$ in $T$, and a positive $\delta$ such that 
\begin{enumerate}
\item there is a normal bundle $\overline{p}:\overline{N}\to \Sigma\cap (W_t\times \partial V)$ 
to $\Sigma$ such that for a certain
positive $\eta$ the space $T_{\overline{N}}(\eta)$ embeds as a tubular neighborhood of $\Sigma\cap (W_t\times \partial V)$
in $W_t\times \partial V$,
\item the fibre $f_{t'}^{-1}(s)$ meets transversely (in the stratified sense when $s=0$) the boundary
$\partial V$ for any $t'\in W_t$ and $s\in D_\delta^*$.
\end{enumerate}
\end{lema}
\begin{proof}
There is a normal bundle to $\Sigma_t$ called $p:N\to\Sigma_{t}\cap\partial V$, and a positive $\eta$ such that
$T_N(\eta)$ is naturally embedded in $\partial V$. There is a neighborhood $W_t$ of $t\in T$ 
such that $\Sigma\cap(W_t\times \partial V)$ is included in $W_t\times T_N(\eta)$, and the restriction 
\[\tau|_{\Sigma\cap(W_t\times \partial V)}:\Sigma\cap(W_t\times \partial V)\to W_t\]
is submersive (hence a locally trivial fibration).

For any $t'\in W_t$ the fibre $\Sigma_{t'}\cap\partial V$ is a disjoint union of circles,
contained in $\{t'\}\times T_N(\eta)$. In each connected component of $\{t'\}\times T_N(\eta)$ we find exactly one circle. 
Define
\[\Xi:W_t\times\CC^n\to\{t\}\times\CC^n\]
by $\Xi(t',x):=(t,x)$. Fix any $t'\in W_t$. The inclussion 
\[\iota:T_N(\eta)\hookrightarrow N\]
induces a fibration 
\[p\comp\iota\comp\Xi|_{\{t'\}\times T_N(\eta)}:\{t'\}\times T_N(\eta)\to\Sigma_t\cap\partial V.\] 
If $W_t$ is small enough and $t'\in W_t$, each circle of $\Sigma_{t'}\cap\partial V$
meets transversely the fibres of $p\comp\iota\comp\Xi|_{\{t'\}\times T_N(\eta)}$. 
Observe that the fibre of $p\comp\iota\comp\Xi|_{\{t'\}\times T_N(\eta)}$ at any point $(t,y)\in\Sigma_{t}\cap\partial V$
is a $n-1$-dimensional complex analytic submanifold of $\{t'\}\times\CC^n$.
Given $(t',x)\in\cap\Sigma_{t'}\cap\partial V$ we denote by $H_{(t',x)}$ the 
complex tangent space at $(t',x)$ of the unique fibre of $p\comp\iota\comp\Xi|_{\{t'\}\times T_N(\eta)}$ to which $(t',x)$ belongs.
Define 
\[\overline{p}:\overline{N}\to \Sigma\cap (W_t\times \partial V)\] 
as the bundle whose fibre over
any point $(t',x)\in \Sigma\cap (W_t\times \partial V)$ is $H_{(t',x)}$.
Clearly $\overline{p}$ is a normal bundle extending $p$, and, if $\eta$ is small enough then 
$T_{\overline{N}}(\eta)$ is naturally embedded in $W_t\times \partial V$ and satisfies Condition (1).

The space $T_{\overline{N}}(\eta)$ has exactly one connected component $C_i$ for each connected 
component of $\Sigma_{t}\cap\partial V$. Let $\mu_i$ be the transversal Milnor number 
corresponding to the $i$-th connected component of $\Sigma_{t}\cap\partial V$.
As $\overline{p}$ is a normal bundle, by Property (3) of 
Definition~\ref{equiatsing} and Lemma~\ref{gentranstype}, we deduce the following: for any  
connected component $C_i$, the restriction of $f$ to any fibre of $\overline{p}_{C_i}$
is an analytic function which has an isolated singularity at the origin with Milnor number $\mu_i$.
By the constancy of the Milnor number we get, perhaps shrinking $W_t$ and diminishing $\eta$,
that the restriction of $f_{t'}$ the the fibres of $C_i$ have no critical points outside the origin.
This implies the existence of a positive $\delta$ such that $f_{t'}^{-1}(s)$ is transverse (in the stratified sense if $s=0$) to 
$\partial V$ in $T_{\overline{N}}(\eta)$, for any $t'\in W_t$ and $s\in D_\delta$.

If $W_t$ and $\delta$ are sufficiently small and $t'\in W_t$ and $s\in D_\delta$, 
 the transversality of $f_{t'}^{-1}(s)$ to
$\partial V$ outside $T_{\overline{N}}(\eta)$ follows because $f_t^{-1}(0)$ is transverse to
$\partial V$ outside $T_{\overline{N}}(\eta)$ and the transversality is an open condition.
\end{proof}

\section{The embedded topological type and the Milnor fibration}

\noindent\textbf{Assumption}. Thorought the rest of the paper we assume $n\geq 5$.

Along this section we let 
$f_t:(\CC^n,O)\to\CC$ be a smooth family of holomorphic germs depending on a parameter $t$ varying in a smooth manifold $T$.
Let $f:T\times\CC^n\to\CC$ be the smooth function defining the family. We assume that the critical set of $f_t$ at the origin is 
of dimension
$1$ for any $t$, and that the family is {\em topologically equisingular at the critical set}. Denote by 
\[\psi:T\times\CC^n\to T\times\CC\]
the mapping $(\tau,f)$.

\noindent\textbf{Notation}. Given $A\subset T\times\CC^n\times$ and $B\subset T$ we define $A_B=A\cap\tau^{-1}(B)$.

\subsection{Cuts}
\label{1inicial}

\begin{definition}
\label{1cuts}
A {\em cut for} $f$ with {\em amplitude} $\delta$ {\em over a submanifold} $V\subset T$ is a closed smooth hypersurface $H$ of
$\psi^{-1}(V\times D_\delta)$ with the following properties:
\begin{enumerate}
\item There is a tubular neighborhood of $\Sigma\cap H$ in $H$ which coincides
with $T_{N}(\eta)$ for a certain normal bundle to $\Sigma$ called $\pi:N\to\Sigma\cap H$, and a certain $\eta>0$.
\item For any $(t,s)\in V\times D_\delta$ the hypersurface $H$ meets the (open) set of smooth points of $\phi^{-1}(t,s)$ transversely.
\item There is a unique connected component $X_{int}(H,V,\delta)$ (the {\em interior component}) 
of $\psi^{-1}(V\times D_\delta)\setminus H$, which contains the zero section. Moreover the restriction 
\[\tau|_{X_{int}(H,V,\delta)}:X_{int}(H,V,\delta)\to V\] 
is a smooth locally trivial fibration, with fibre a contractible compact manifold with corners.
\item For any $t\in V$ and any neighborhood $W$ of the origin adapted to $f_t$ and 
contained in $\dot{X}_{int}(Y,V,H)_t$, for any Milnor disk radius
$\delta'$ for $(f_t,W)$ satisfying $\delta'<\delta$, and any $s\in D_{\delta'}^*$ we have that 
$f_t^{-1}(s)\cap (X_{int}(Y,V,H)_t\setminus\dot{W})$ is a homotopically trivial cobordism.
\item For any $t\in V$ and any neighborhood $W$ of the origin adapted to $f_t$ and 
contained in $\dot{X}_{int}(Y,V,H)_t$ we have that 
$\Sigma_t\cap (X_{int}(Y,V,H)_t\setminus\dot{W})$ is a trivial cobordism.
\item The previous property implies that, for any $t\in V$, there is a bijective correspondence (induced by inclusion) 
between the connected components of $\Sigma_t\cap X_{int}(Y,V,H)_t\setminus \{(t,O)\}$ and 
$\Sigma_{V}\cap X_{int}(Y,V,H)\setminus V$.
For any connected component $\Gamma$ of $\Sigma_{V}\cap X_{int}(Y,V,H)\setminus V $ and any $(t,p)\in \Gamma$, the 
transversal Milnor number of $f_t$ at $p$ is independent on the point $(t,p)$.
\end{enumerate}
\end{definition}

\begin{remark}
By the proof of Lemma~\ref{homotopylink} it is sufficient to check Property (4) for a single neighborhood adapted to $f$.
\end{remark}

\begin{remark}
\label{gnc}
The definition of system of neighborhoods adapted to $f_t$ implies that for any neighborhood
$W$ of the origin adapted to $f_t$, there exists $\delta$ such that 
$\psi^{-1}(\{t\}\times D_\delta)\cap\partial W$ is a cut over $\{t\}$ of amplitude $\delta$.
\end{remark}

The relations $\preceq$ and $\prec$ are defined as in Section~\ref{inicial}.

\subsection{Construction}
\label{construccion}
Let $C=[0,1]^d$ be a cube embedded in $T$ and $V$ an open neighborhood of $C$ in $T$ diffeomorphic to the open cube
$(-\epsilon,1+\epsilon)^d$ (for $\epsilon>0$). Suppose that we have cuts $H_+$ and $H_-$ for $f$ over $V$ of the same amplitude
$\delta$, such that $H_-\prec H_+$. Along this section we will shrink $\delta$ when it is necessary without explicitly mentioning it.
Define $X:=X_{int}(H_+,V,\delta)\setminus\dot{X}_{int}(H_-,V,\delta)$.

By Definition~\ref{1cuts}, an easy argument with cobordisms, and Ehresmann Fibration Theorem,
the restriction
\begin{equation}
\tau|_{\Sigma\cap X}:\Sigma\cap X\to V
\end{equation}
is a locally trivial fibration (trivial in fact, since the base is contractible), with fibre a trivial cobordism diffeomorphic to a 
disjoin union of cylinders, one for each connected component of $\Sigma\setminus V$. We will assume for notational simplicity that
$\Sigma\setminus V$ is connected, being the treatment of the general case completely analogous. We consider a diffeomorphism  
\begin{equation}
\label{productstructure}
\nu:\Sigma\cap X\to \SSS^1\times [0,1]\times V
\end{equation}
such that $\tau|_{\Sigma\cap X}=p_3\comp\nu$ (being $p_i$ the projection of $\SSS^1\times [0,1]\times V$
 to the $i$-th factor for $i=1,2,3$), 
$\nu(\Sigma\cap H_-)=\SSS^1\times\{0\}\times V$ and $\nu(\Sigma\cap H_+)=\SSS^1\times\{1\}\times V$.

By Definition~\ref{1cuts} there exist a normal bundle $\pi_+:N_+\to \Sigma\cap H_+$ and
a neighborhood of $\Sigma\cap H_+$ in $H_+$ which is equal to $T_{N_+}(\eta)$ for a certain positive
$\eta$.
Similarly there exist a normal bundle $\pi_-:N_-\to \Sigma\cap H_-$ and
a neighborhood of $\Sigma\cap H_-$ in $H_-$ coinciding with $T_{N_-}(\eta)$.

For any $t\in V$, if $z_1$ is a coordinate of $\CC^n$ not vanishing at $\Sigma_t$, the $z_1$-section 
(that is, the hyperplane $V(z_1-z_1(x))$)
is transverse to $\Sigma_t$ at any point $x\in\Sigma_t$, with finitely many exceptions. 
This allows to find an annulus $A$ embedded in $\Sigma\cap X_t$ not meeting the boundary, whose embedding in $\Sigma\cap X_t$
is a homotopy equivalence, such that the $z_1$-section is transverse to $\Sigma$ at any point of $A$. Consequently the vector subbundle
$\pi_A:N_A\to A$ of $\calT(T\times\CC^n,T)_A$ 
whose fibre over $a\in A$ is the $z_1$-section of the fibre $\calT(T\times\CC^n,T)_a$ is a
constant (and hence holomorphic) normal bundle to $\Sigma$ over $A$.

Using Lemma~\ref{normalbundle} we extend the bundles $\pi_+$, $\pi_-$ and $\pi_A$ to a normal bundle
\begin{equation}
\label{fibradonormal}
\pi:N\to\Sigma\cap X.
\end{equation}
We choose $\eta$ small enough so that $T_N(\eta)$ embeds as neighborhood of $\Sigma\cap X$ in $X$.
Observe that we have automatically the compatibility 
\[\tau\comp\pi|_{T_N(\eta)}=\tau|_{T_N(\eta)}.\]
By the last condition of
Definition~\ref{1cuts} and by Lemma~\ref{gentranstype}
the Milnor number at the origin of the restriction of $f$ to the fibre of $\pi$ over $x\in \Sigma\cap X$
is independent of $x$. Observe that $\Sigma\cap X_C$ is diffeomorphic to the product of $C$ and a cylinder, which is the same than 
the product of a circle and a $(d+1)$-dimensional cube. Since the restriction $\pi_A$ is holomorphic, an application of 
Corollary~\ref{assumptionA} enables us to use the results of Section~\ref{globalisationLR} to the restriction of the function $f$ to
$T_N(\eta)$, viewed as a $\mu$-constant family parametrised over $\Sigma\cap X$.

Thus, by Lemma~\ref{cutcircle} and Proposition~\ref{usefulconstr}, there exists a cut $S$ for $f$ over $\Sigma\cap X$ with amplitude
$\delta$ such that $Y:=Y_{int}(S,\Sigma\cap X,\delta)$ is contained in $T_N(\eta)$. We consider the composition 
$\nu\comp\pi:Y\to \SSS^1\times [0,1]\times C$. Given any $t\in C$, consider the circle 
\[K_t:=\SSS^1\times\{0\}\times\{t\}\] 
inside $\SSS^1\times [0,1]\times C$, and define 
\[Y_{K_t}:=(\nu\comp\pi)^{-1}(K_t)=\pi^{-1}(\Sigma_t\cap H_-)=Y_t\cap H_-.\] 
By Theorem~\ref{globalisation} there exist a homeomorphism
\begin{equation}
\label{1homeoisolated}
\Psi:Y\to Y_{K_t}\times [0,1]\times C
\end{equation}
such that
\begin{enumerate}
\item We have the equality $\nu\comp\pi=(p_1\comp\nu\comp\pi|_{Y_{K_t}},q_2,q_3)\comp\Psi$,
where $q_i$ is the projection of
$Y_{K_t}\times [0,1]\times C$ to the $i$-th factor.
\item We have $f|_{Y_{K_t}}\comp q_1\comp\Psi=f$.
\item the restriction of $\Psi$ to $\Sigma\cap X$ coincides with $\nu$, and the restriction of $\Psi$ to
$Y\setminus (\Sigma\cap X)$ is smooth.
\end{enumerate}

We will use the notation $Z=\overline{X\setminus Y}$. The subspaces $Y$ and $Z$ meet at a common boundary which coincides with the 
cut $S$. As for any $t\in C$ and $s\in D_\delta$ the fibre $f_t^{-1}(s)$ meets $S$ transversely, by Ehresmann Fibration Theorem 
and the contractibility of the base the restriction
\begin{equation}
\label{outerfib}
\psi|_Z:Z\to C\times D_\delta
\end{equation}
is a trivial fibration, which we call the {\em outer fibration}.

\subsection{Topological trivialisation of the space between two cuts}
We start studying the outer fibration using cobordism theory.

\begin{lema}
\label{hcobordisms}
For any $s\in D_\delta^*$, the cobordisms
\begin{equation}
\label{total}
(f_t^{-1}(s)\cap X,f_t^{-1}(s)\cap H_+,f_t^{-1}(s)\cap H_-)
\end{equation}
\begin{equation}
\label{partialouter}
(f_t^{-1}(s)\cap Z,f_t^{-1}(s)\cap Z\cap H_+,f_t^{-1}(s)\cap Z\cap H_-)
\end{equation}
are simply connected $h$-cobordisms, the second one with boundary.
\end{lema}
\begin{proof}
By Property~(4) of Definition~\ref{1cuts}, we may find a neighborhood $W$ of the origin adapted to $f_t$ such that
$W$ is contained in $\dot{X}_{int}(H_+,V,\delta)$, $\dot{X}_{int}(H_-,V,\delta)$ and
\begin{equation}
\label{joput1}
(f_t^{-1}(s)\cap X_{int}(H_+,V,\delta)\setminus\dot{W},f_t^{-1}(s)\cap H_+,f_t^{-1}(s)\cap \partial W)
\end{equation}
\begin{equation}
\label{joput2}
(f_t^{-1}(s)\cap X_{int}(H_-,V,\delta)\setminus\dot{W},f_t^{-1}(s)\cap H_-,f_t^{-1}(s)\cap \partial W)
\end{equation}
are homotopically trivial cobordisms for $s$ non-zero and sufficiently small. This easily implies that the 
cobordism~(\ref{total}) is also homotopically trivial, and that the three spaces involved in the cobordism 
have the homotopy type of $f_t^{-1}(s)\cap \partial W$, which is the boundary of the Milnor fibre of 
$f_t$ for the adapted neighborhood $W$.

Consider $\pi:Y\to\Sigma\cap X$ as in Construction~\ref{construccion}.
Define $\phi:Y\to(\Sigma\cap X)\times D_\delta$ by $\phi:=(\pi,f)$.
Consider the space 
\[B=\phi^{-1}((\Sigma_t\cap H_-)\times D_\delta).\] 
Observe that $B$ is included in $H_-$ and that we have the equality
\[f_t^{-1}(s)\cap B=\phi^{-1}((\Sigma_t\cap H_-)\times\{s\}).\] 
The mappings 
\begin{equation}
\label{engorde}
\pi|_B:B\to \Sigma_t\cap H_-
\end{equation}
\begin{equation}
\label{inters}
\pi|_{f_t^{-1}(s)\cap B}:f_t^{-1}(s)\cap B\to \Sigma_t\cap H_-
\end{equation}
are locally trivial fibrations over a circle with simply connected fibres: 
the fibres of the first mapping are contractible, and the fibres of the 
second are homeomorphic to the Milnor fibre of the transversal singularity at any point of
$\Sigma_t$, which is simply connected by the 
Kato-Matsumoto bound (recall that we have assumed $n\geq5$). Hence the inclusion $f_t^{-1}(s)\cap B\subset B$ induces an
isomorphism of (infinite cyclic) fundamental groups. Using this,
an easy application of Seifert-van Kampen Theorem shows that the fundamental 
group of $f_t^{-1}(s)\cap H_-$ is isomorphic to the fundamental group of $(f_t^{-1}(s)\cap H_-)\cup B$.

The restriction of the mapping~(\ref{outerfib}) to $(Z\cap H_-)_t$ yields a trivial fibration
\begin{equation}
\label{soloaqui}
\psi|_{(Z\cap H_-)_t}:(Z\cap H_-)_t\to \{t\}\times D_\delta.
\end{equation}
Using it we deduce that the space $(f_t^{-1}(s)\cap H_-)\cup B$ is homotopy equivalent to $f_{t}^{-1}(D_\delta)\cap H_-$. 
Using the fibration~(\ref{soloaqui}) we show that $f_t^{-1}(D_\delta)\cap H_-$ admits $(f_t^{-1}(0)\cap H_-)\cup B$
as a deformation retract. Working with the fibration~(\ref{engorde}) we prove that, in turn, the space 
$(f_t^{-1}(0)\cap H_-)\cup B$ admits $f_t^{-1}(0)\cap H_-$ as a deformation retract. 
Thus $(f_t^{-1}(s)\cap H_-)\cup B$ and $f_{t}^{-1}(0)\cap H_-$ 
have isomorphic fundamental groups. Consequently $f_t^{-1}(s)\cap H_-$ and $f_{t}^{-1}(0)\cap H_-$ have isomorphic fundamental 
groups for any $(s,t)\in C\times D_\delta)$.

We can show in the same way that the spaces 
$f_t^{-1}(s)\cap H_+$, and $f_t^{-1}(0)\cap H_+$ 
(respectively $f_t^{-1}(s)\cap X$, and $f_t^{-1}(0)\cap X$) have isomorphic 
fundamental groups for any $(t,s)\in C\times D_\delta$. 

With a similar construction we can show that 
$f_t^{-1}(s)\cap\partial W$ has the same fundamental group than $f_t^{-1}(0)\cap\partial W$, which is 
homotopy equivalent to the classical link of $f_t$. The later space is simply connected 
by~\cite{Mi}, Theorem~5.2 (use that $n\geq 5$). We have shown that~(\ref{total}) is a simply connected $h$-cobordism 

Consider the decomposition 
\begin{equation}
\label{decutil}
f_t^{-1}(s)\cap H_-=(f_t^{-1}(s)\cap H_-\cap Y)\bigcup (f_t^{-1}(s)\cap H_-\cap Z)
\end{equation}
The mapping $\pi$ fibres the first piece and the intersection of the two pieces over the circle $\Sigma_t\cap H_-$, with fibres
the Milnor fibre and the link of the transversal singularities respectively. 
The link of the transversal singularity is simply connected
due to Theorem~5.2 of~\cite{Mi} (we have assumed $n\geq 5$). Seifert-van Kampen Theorem implies now that 
$f_t^{-1}(s)\cap H_-$ and $f_t^{-1}(s)\cap H_-\cap Z$ have isomorphic fundamental groups. Therefore $f_t^{-1}(s)\cap H_-\cap Z$ is
simply connected. Similarly we
show that $f_t^{-1}(s)\cap H_+\cap Z$ and $f_t^{-1}(s)\cap Z$ are simply connected.

It remains to be shown the following vanishing of relative homology groups: 
\[H_*(f_{t}^{-1}(s)\cap Z,f_{t}^{-1}(s)\cap Z\cap H_-;\ZZ)=0\]
\[H_*(f_{t}^{-1}(s)\cap Z,f_{t}^{-1}(s)\cap Z\cap H_+;\ZZ)=0.\]

Using the ladder of long exact sequences formed by the Mayer-Vietoris sequences associated to the decomposition~(\ref{decutil}) and
the decomposition 
\begin{equation}
f_t^{-1}(s)=(f_t^{-1}(s)\cap Y)\bigcup (f_t^{-1}(s)\cap Z),
\end{equation}
and the fact that the inclusions 
\[f_t^{-1}(s)\cap H_-\cap Y\subset f_t^{-1}(s)\cap Y,\]
\[f_t^{-1}(s)\cap H_-\cap Y\cap Z\subset f_t^{-1}(s)\cap Y\cap Z\]
are homotopy equivalences (which is true due to the homeomorphism~(\ref{1homeoisolated})), we show that
 the first required vanishing follows from the vanishing
\[H_*(f_{t}^{-1}(s)\cap X,f_{t}^{-1}(s)\cap H_-;\ZZ)=0,\]  
which is true since the cobordism~(\ref{total}) is an homotopically trivial. The second vanishing is proved similarly.
\end{proof}

\begin{lema}
\label{boundarymf}
If $H$ is a cut for $f$ over $V$ with amplitude $\delta$ then $f_t^{-1}(s)\cap H$ is simply connected
for any $(t,s)\in V\times D_\delta$.
\end{lema}
\begin{proof}
We may assume that $H$ is equal to the cut $H_-$ of the previous Lemma. We have shown that $f_t^{-1}(s)\cap H_-$ is simply 
connected for $s\neq 0$, and that the fundamental groups of $f_t^{-1}(s)\cap H_-$ and $f_t^{-1}(0)\cap H_-$ are the same. 
\end{proof}

We have the ingredients for the proof one of our main technical propositions.

\begin{prop}
\label{pivote}
For any $t\in C$ there exists a homeomorphism
\begin{equation}
\label{comparacuts}
\Theta:X\to (H_-)_t\times [0,1]\times C
\end{equation}
such that
\begin{enumerate}
\item We have $\tau=p_3\comp\Theta$, where $p_i$ is the projection of $(H_-)_t\times [0,1]\times C$ to 
the $i$-th factor ($i=1,2,3$). 
\item We have $f|_{(H_-)_t}\comp p_1\comp\Theta=f$.
\item the restriction of $\Theta$ to $\Sigma\cap X$ is smooth with target 
$((H_-)_t\cap\Sigma)\times [0,1]\times C$ and the restriction of $\Theta$ to $X\setminus \Sigma$ is 
smooth.
\item The restriction $\Theta|_Y$ coincides with the homeomorphism~(\ref{1homeoisolated}).
\end{enumerate}
\end{prop}
\begin{proof}
We have the equalities
\[S_{K_t}:=S\cap(\nu\comp\pi)^{-1}(K_t)=S\cap\pi^{-1}(\Sigma_t\cap H_-)=(S\cap H_-)_t.\]

The restriction of the homeomorphism~(\ref{1homeoisolated}) to the boundary $S$ induces a diffeomorphism
\begin{equation}
\label{1homeobound}
\Psi|_S:S\to (S\cap H_-)_t\times [0,1]\times C
\end{equation}
preserving the fibres of $f$. 

The restrictions 
\begin{equation}
\psi|_S:S\to C\times D_\delta
\end{equation}
\begin{equation}
\psi|_{(S\cap H_-)_t}:(S\cap H_-)_t\to \{t\}\times D_\delta
\end{equation}
of the trivial fibration~(\ref{outerfib}) are trivial fibrations as well.
Choose a smooth trivialisation 
\[\varphi_1:=(\gamma,t,f_t):(S\cap H_-)_t\to \psi|_{(S\cap H_-)_t}^{-1}(t,0)\times\{t\}\times D_\delta\]
of $\psi|_{(S\cap H_-)_t}$.

Denote by $q_i$ the projection of $(S\cap H_-)_t\times [0,1]\times C$ to the $i$-th factor ($i=1,2,3$).
Since $\Psi|_S$ is a restriction of the homeomorphism~(\ref{1homeoisolated}), the mapping
\begin{equation}
\label{**}
\Upsilon:S\to\psi|_{(S\cap H_-)_t}^{-1}(t,0)\times[0,1]\times C\times D_\delta
\end{equation}
defined by $\Upsilon:=(\gamma\comp q_1,q_2,q_3,f_t\comp q_1)\comp\Psi|_S$ clearly satisfies
$\psi|_S=(p_3,p_4)\comp\Upsilon$, where $p_i$ is the projection of $\psi|_{(S\cap H_-)_t}^{-1}(t,0)\times [0,1]\times C\times D_\delta$
to the $i$-th factor.

Define 
\[\Xi:\psi|_{(S\cap H_-)_t}^{-1}(t,0)\times [0,1]\times C\times D_\delta\to\psi|_S^{-1}(t,0)\times C\times D_\delta\]
by $\Xi(y,u,t',s):=(\Psi|_S^{-1}(y,u,t),t',s)$. Then the composition
\[\varphi_2:=\Xi\comp\Upsilon:S\to \psi|_S^{-1}(t,0)\times C\times D_\delta\]
is a trivialisation of the fibration $\psi|_S$.

Consider an extension 
\[\varphi_3:Z\to\psi|_Z^{-1}(t,0)\times C\times D_\delta\]
of $\varphi_2$ to a trivialisation of $\psi|_Z$.

By Lemma~\ref{hcobordisms} the triple 
\begin{equation}
\label{***}
(\psi|_Z^{-1}(t,0),\psi|_Z^{-1}(t,0)\cap H_-,\psi|^{-1}_Z(t,0)\cap H_+)
\end{equation}
is a simply connected $h$-cobordism with boundary. The diffeomorphism~(\ref{**}) induces a trivialisation of the boundary cobordism
\[(\psi|_S^{-1}(t,0),\psi|^{-1}_S(t,0)\cap H_-,\psi|^{-1}_S(t,0)\cap H_+).\]
By $h$-cobordism Theorem the cobordism~(\ref{***}) admits a trivialisation 
\[\kappa:\psi^{-1}|_Z(t,0)\to (\psi^{-1}|_Z(t,0)\cap H_-)\times [0,1]\]
which extends the given trivialisation at the boundary.

Consider the diffeomorphism 
\[\Upsilon':Z\to (\psi^{-1}|_Z(t,0)\cap H_-)\times [0,1]\times C\times D_\delta\]
defined by $\Upsilon':=(\kappa\comp\sigma_1,\sigma_2,\sigma_3)\comp\varphi_3$, where $\sigma_i$ is the projection 
of $\psi|_Z^{-1}(t,0)\times C\times D_\delta$ to the $i$-th factor.
By definition $\Upsilon'$ extends $\Upsilon$.

The mapping $\varphi_3$ restricts to a trivialisation 
\[\varphi_3|_{(Z\cap H_-)_t}:(Z\cap H_-)_t\to\psi|_Z^{-1}(t,0)\cap H_-\times \{t\}\times D_\delta\]
of the fibration $\psi|_{(Z\cap H_-)_t}:(Z\cap H_-)_t\to\{t\}\times D_\delta$. Define
\[\Psi':Z\to (Z\cap H_-)_t\times [0,1]\times C\]
by $\Psi'(z)=(\varphi_3|^{-1}_{(Z\cap H_-)_t}(x,s),u,t')$, where $(x,u,t',s):=\Upsilon'(z)$.

As $\Upsilon$ extends $\Upsilon'$ we obtain that $\Psi$ and $\Psi'$ concide at the intersection of their domains.
Hence the desired homeomorphism $\Theta$ 
can be defined piecewise over $Y$ and $Z$ gluing $\Psi$ and $\Psi'$.  

In order to fulfill the differentiability of $\Theta$ at $S$ it is sufficient take a collar of $S$ and 
modify $\Psi$ and $\Psi'$ in each of the sides so that the gluing is smooth.
\end{proof}

\begin{cor}
\label{indeplinkmilnorfib}
Let $f:(\CC^n,O)\to\CC$ be a holomorphic germ with critical set of dimension $1$ at the origin. Given a system of neighborhoods
$\{U_\alpha\}_{\alpha\in A}$ of the origin adapted to $f$, the homeomorphism type of the embedded link and the diffeomorphism type of
the Milnor fibration of $f$ associated to $\{U_\alpha\}_{\alpha\in A}$ coincide with the classical ones.
\end{cor}
\begin{proof}
In Proposition~\ref{ejemplogoodneigh} we exhibit a system of neighborhoods adapted to $f$ for which the statement is true.
Suppose that we have two neighborhoods $W_1$, $W_2$ of the origin adapted to $f$ such that $W_1\subset\dot{W}_2$. By Remark~\ref{gnc}
there is $\delta>0$ such that $\partial W_i\cap f^{-1}(D_\delta)$ is a cut for $f$ of amplitude $\delta$ for $i=1,2$.
A straightforward application of Proposition~\ref{pivote} gives the result. 
\end{proof}

\subsection{The embedded topological type and the Milnor fibration}

We show that being a cut is an open property on the base: 

\begin{lema}
\label{opennature}
Let $f_t:(\CC^n,O)\to\CC$ be a smooth family of holomorphic germs depending on a parameter $t$ varying in a smooth manifold $T$,
and $V$ be a compact submanifold of $T$. Assume that $f$ is equisingular at the critical set.
Suppose that we have a positive $\delta>0$ and a closed smooth hypersurface $H$ of
$\psi^{-1}(T\times D_\delta)$ which satisfies Property (1) of Definition~\ref{1cuts} and which is such that the restriction $H|_V$ is a cut over $V$ of amplitude $\delta$.
Then there is an open neighborhood $U$ of $V$ in $T$ such that $H|_U$ is a cut over $U$ of
amplitude $\delta$.
\end{lema}
\begin{proof}
We have to check that $H_U$ satisfies Properties (2)-(6) of Definition~\ref{1cuts} for a certain open neighborhood $U$ of $V$ in
$T$. The second property follows from the fact that Property (1) holds and an argument like in the proof of 
Lemma~\ref{conormality}. Let $U$ be an open neighborhood over which Properties (1) and (2) hold. It is easy to check 
using Ehresmann fibration Theorem that 
Property (3) holds as well. The fact that Property (5) holds (after possibly shrinking $U$) 
is an easy argument involving manipulations with cobordisms, Ehresmann fibration Theorem, and the fact that
cobordism~(\ref{coborcomp}) is trivial. To show that Property (6) holds (after possibly shrinking $U$) we only have to use that
it holds over $V$, that $V$ is compact and that the family is equisingular at the critical set.

Proving Property (4) is slightly more involved. Let $U$ be an open neighborhood of $V$ in $T$ over which
Properties (1)-(3), (5) and (6) are satisfied, and such that $V$ meets all the connected components of 
$U$. Then $X:=X_{int}(H,U,\delta)$ is defined. The mapping
\begin{equation}
\label{fibtotspac}
\tau:X\to U
\end{equation}
is a locally trivial fibration with contractible fibres, and
the mapping 
\begin{equation}
\label{milnorextend}
\psi=(\tau,f):X\setminus f^{-1}(0)\to U\times D_\delta^*
\end{equation}    
is a locally trivial fibration whose fibre is diffeomorphic to the Milnor fibre of $f_t$ for any $t\in V$
(notice that $H$ is a cut over $V$).

Choose any $t\in U$. Let $W$ be a neighborhood of the origin adapted to $f_t$ satisfying that
$\{t\}\times W$ is contained in $\dot{X}_t$. We will prove that the cobordism 
\begin{equation}
\label{objetivo}
(f_t^{-1}(s)\cap X_t\setminus\dot{W},f_t^{-1}(s)\cap H,f_t^{-1}(s)\cap\partial W)
\end{equation}
is a simply connected $h$-cobordism for any $s\in D_\delta^*$.

Using the fibration~(\ref{milnorextend}) and Lemma~\ref{boundarymf} we deduce that 
$f_t^{-1}(s)\cap H$ is simply connected for any $(t,s)\in U\times D_\delta$. The space 
$f_t^{-1}(s)\cap \partial W$ is simply connected by Lemma~\ref{boundarymf} since $\partial W$ is a 
cut over $t$. By fibration~(\ref{milnorextend}) and the fact that $H$ is a cut over $V$ we deduce that
$f_t^{-1}(s)\cap X_t$ has the same homotopy type than the Milnor fibre of $f_t$, which is simply 
connected by Kato-Matsumoto bound. As $W$ is a neighborhood adapted to $f$ the space 
$f_t^{-1}(s)\cap W$ is also homotopic to the Milnor fibre of $f_t$. Seifert-van Kampen Theorem applied
to the decomposition
\[f_t^{-1}(s)\cap X_t=f_t^{-1}(s)\cap X_t\setminus\dot{W}\cup f_t^{-1}(s)\cap W\]
gives that $f_t^{-1}(s)\cap X_t\setminus\dot{W}$ is simply connected.

To finish the proof it is enough to show that all the relative homology groups 
\[H_*(f_{t}^{-1}(s)\cap X_t\setminus\dot{W},f_{t}^{-1}(s)\cap\partial W;\ZZ)\] 
vanish. By excision it is equivalent to show the vanishing 
\[H_*(f_{t}^{-1}(s)\cap X_t,f_{t}^{-1}(s)\cap W;\ZZ)=0.\] 
Using the inclusions $f_{t}^{-1}(s)\cap X_t\subset X_t$, 
$f_{t}^{-1}(s)\cap W\subset f_t^{-1}(D_\delta)\cap W$, the fact that $X_t$ and $f_t^{-1}(D_\delta)\cap W$
are contractible, and long exact sequences in homology, it is easy to show that the wanted vanishing is 
equivalent to the vanishing
\[H_*(X_t\setminus\dot{W},(f_{t}^{-1}(s)\cap X_t\setminus\dot{W})\cup (f_{t}^{-1}(D_\delta)\cap\partial W);\ZZ)=0.\]

As $\partial W$ is a cut for $f_t$ over $t$ with amplitude $\delta$ and 
$X_{int}(\partial W,\{t\},\delta)=f_t^{-1}(D_\delta)\cap W$ is contained in $\dot{X}_t$, performing 
Construction~\ref{construccion} for $f_t$, $T=C=\{t\}$ we obtain an splitting
$X_t\setminus\dot{W}=Y\cup Z$, such that the restriction
\[f_t|_Z:Z\to D_\delta\] 
is a trivial fibration and the space $Y$ is equal to $f_t^{-1}(D_\delta)\cap T_N(\eta)$ for a
a normal bundle 
$p:N\to\Sigma_t\cap (X_t\setminus\dot{W})$ to $\Sigma_t$ and a positive $\eta$.

Using the triviality of $f_t|_Z$ and homology excision we reduce our problem to prove the vanishing of 
\[H_*(Y,Y\cap (f_t^{-1}(s)\cup\partial W));\ZZ)=0,\]
 which is easy using the homeomorphism~(\ref{1homeoisolated}) and the homotopy invariance of homology.
\end{proof}

\begin{theo}
\label{premasseyemblink}
Assume $n\geq 5$.
Let $f_t$ be a smooth family of holomorphic germs parametrised over a connected family $T$, such that 
$\Sigma_t$ is $1$-dimensional at the origin for any $t$. If $f_t$ is equisingular at the critical set 
then the diffeomorphism type of the Milnor fibration and the homeomorphism type of the embedded
link of $f_t$ is independent of $t$. 
\end{theo}
\begin{proof}
As $T$ is connected we are reduced to prove a local statement in the base. Given any $t\in T$ we consider
an adapted neighborhood $W$ to $f_t$. By Corollary~\ref{indeplinkmilnorfib}, for a certain positive $\delta$, the 
Milnor fibration of $f_t$ is $C^\infty$-equivalent to 
\[f_t:W\cap f_t^{-1}(\partial D_\delta)\to\partial D_\delta\]
and the embedded link is homeomorphic to
\[(\partial (W\cap f_t^{-1}(D_\delta)),f_t^{-1}(0)\cap\partial W).\]

Consider another neighborhood of the origin $W'$ adapted to $f_t$, such that we have the inclusion $W'\subset\dot{W}$. 
By Remark~\ref{gnc} there exists $\delta>0$ such that 
\[(\{t\}\times\partial W)\cap f_t^{-1}(D_\delta)\]
and
\[(\{t\}\times\partial W')\cap f_t^{-1}(D_\delta)\] 
are cuts over $t$ of amplitude $\delta$. By 
Lemma~\ref{conormality}~(1) and 
Lemma~\ref{opennature} there is a neighborhood $C$ of $t$ in $T$ (which can be taken take cubical) such 
that 
\[H_+:=(C\times\partial W)\cap\psi^{-1}(C\times D_\delta)\]
and
\[H_-:=(C\times\partial W')\cap\psi^{-1}(C\times D_\delta)\] 
are cuts over $C$ of amplitude $\delta$. 

Then the fibration 
\[\psi:X_{int}(H_+,C,\delta)\setminus f^{-1}(0)\to C\times D_\delta^*\]
is locally trivial. Consequently the diffeomorphism type of the fibration
\begin{equation}
\label{modelmilfib}
f_{t'}:\{t'\}\times W\cap f_{t'}^{-1}(\partial D_\delta)\to\partial D_\delta
\end{equation}
is independent of $t'\in C$. For $t'=t$ we obtain the diffeomorphism type of the Milnor fibration of $t$.

An easy argument using Proposition~\ref{pivote} applied to $H_-$ and $H_+$ easily implies that the homeomorphism type of the pair
\begin{equation}
\label{modellink}
((\partial X_{int}(H_+,C,\delta)_{t'},\partial X_{int}(H_+,C,\delta)_t\cap f_{t'}^{-1}(0))
\end{equation}
is independent of $t'\in C$. For $t'=t$ we obtain the homeomorphism type of the embedded link of $f_t$.

Fix any $t'\in C$. Consider any neighborhood $W''$ adapted to $f_{t'}$ such that $W''\subset\dot{W}$.
By Corollary~\ref{indeplinkmilnorfib} the Milnor fibration of $f_{t'}$ is diffeomorphic to 
\begin{equation}
\label{modelmilfib'}
f_{t'}:W''\cap f_{t'}^{-1}(\partial D_\delta)\to\partial D_\delta
\end{equation} 
and the embedded link of $f_{t'}$ is homeomorphic to  
\begin{equation}
\label{modellink'}
(\partial (W''\cap f_{t'}^{-1}(D_\delta)),f_{t'}^{-1}(0)\cap\partial W)
\end{equation}
for $\delta$ small enough.

By Remark~\ref{gnc} if we shrink $\delta$ enough we have that  
$(\{t'\}\times\partial W''\{t'\})\cap f_{t'}^{-1}(D_\delta)$ is a cut over $t'$ of amplitude $\delta$.
Another application of Proposition~\ref{pivote}, now for $C=\{t'\}$, 
$H_+:=(\{t'\}\times\partial W)\cap f_{t'}^{-1}(D_\delta)$ and 
$H_-:=(\{t'\}\times\partial W'')\cap f_{t'}^{-1}(D_\delta)$ implies that the diffeomorphism types of 
the fibrations~(\ref{modelmilfib}) and~(\ref{modelmilfib'}) are the same and that the 
pairs~(\ref{modellink}) and~(\ref{modellink'}) are homeomorphic.      
\end{proof}

\section{Topological $R$-equivalence} Our aim is to prove that any family $f_t$ of analytic germs, parametrised over 
a  cube $C$, with critical set of dimension 
at most $1$, which is topologically equisingular at the critical set, is, in fact, topologically equisingular with respect to 
$R$-equivalence. Once we have the results of the previous section we can follow closely the strategy that we used for $\mu$-constant
families in Section~\ref{globalisationLR}. We do it now pointing specially the aspects in which the proofs are different.

\subsection{Extension of cuts: the non-isolated case}

We define extension of cuts as in Definition~\ref{extension}. Consider the closed and open cubes $C$ and $U$ as in~\ref{segunda}.
Consider two cuts $H_+$ and $H_-$ for $f$ defined over $U$ of the same amplitude $\delta$, pairs 
\[(B_1,\mu_1),...(B_{k_1},\mu_{k_1})\quad\quad (C_1,\nu_1),...(C_{k_2},\nu_{k_2}),\]
as in ~\ref{segunda}, and a cut $H_0$ of amplitude $\delta$ over a neighbourhood $V$ of a contractible union of faces $A$.

We show that Lemma~\ref{extcut} still holds in the new setting:

\begin{lema}
\label{1extcut}
In the setting above, after possibly shrinking $U$ and $V$ to smaller neighborhoods of the corresponding sets, there exists a cut 
$H'_0$ over $U$ which extends $H_0$ and satisfies $H_-\prec H_0$, $\mu_i\prec H_0$ for any $i$, $H_0\prec \nu_i$ for 
any $i$, and $H_0\prec H_+$.
\end{lema}
\begin{proof}
The structure of the proof is the same than the one of Lemma~\ref{extcut}. However we need to adapt several arguments.

Define 
\[X:=\overline{X_{int}(H_+,U,\delta)\setminus X_{int}(H_-,U,\delta)}.\]

\textbf{Step 1}: we shall inductively reduce to the case in which $k_1+k_2=0$. Suppose that $k_1+k_2>0$. Assume that $k_2>0$.

Observe that by the definition of a cut the projections 
\[\tau:\Sigma\cap X\to U,\] 
\[\tau:\Sigma\cap (X_{int}(H_+,C_1',\delta)\setminus (\dot{X}_{int}(\nu_1,U,\delta))\to C'_1\]
\[\tau:\Sigma\cap ({X}_{int}(\nu_1,U,\delta)\setminus (\dot{X}_{int}(H_-,C'_1,\delta))\to C'_1\]
are trivial fibrations with fibres trivial cobordisms diffeomorphic to a disjoin union of cylinders, one for each connected
component of $\Sigma\cap H_+$. We will assume for simplicity that there is a unique connected component, being the general case 
analogous. We denote $\Sigma\cap X$ by $\Sigma'$.
Taking into account that the "$\prec$" relations between the cuts are preserved by intersecting with
$\Sigma$, an argument like Step 1 of Lemma~\ref{extcut} shows that (after possibly shrinking $U$ and
$C'_1$ to smaller neighborhoods) there exists a smooth closed 
hypersurface $K\subset\Sigma'$ with the following properties:
\begin{enumerate}
\item (extension) we have $K_{C'_1}=\nu_1\cap\Sigma_{C'_1}$.
\item The space $\Sigma'\setminus K$ has two connected components
$A_+$ and $A_-$ each of them containing respectively $H_+\cap\Sigma$ and $H_-\cap\Sigma$. The component $A_-$ contains moreover the 
intersections $H_0\cap\Sigma$, $\mu_i\cap\Sigma$ and $\nu_j\cap\Sigma$ for any $i$ and $j\neq 1$.
\item The restrictions $\tau:\overline{A}_-\to U$ and $\pi:\overline{A}_+\to U$ are locally trivial fibrations with fibre a trivial cobordism (a cylinder).
\end{enumerate}
In this situation there exists a diffeomorphism 
\[\varphi:\Sigma'\to (\Sigma'\cap H_-)\times [0,1]\]
such that 
\[\varphi(H_-\cap \Sigma')=(H_-\cap \Sigma')\times\{0\},\] 
\[\varphi(H_+\cap \Sigma')=(H_-\cap \Sigma')\times\{1\},\]
\[\varphi(K\cap \Sigma')=(H_-\cap \Sigma')\times\{1/2\}\] 
and 
\[\tau=\tau\comp q_1 \comp\varphi,\]
with $q_i$ the projection of $(\Sigma'\cap H_-)\times [0,1]$ to the $i$-th factor $i=1,2$.

Let $\pi_+:M_+\to\Sigma'\cap H_+$, $\pi_-:M_-\to\Sigma'\cap H_+$ and $\pi_{\nu_1}:M_{\nu_1}\to \nu_1\cap\Sigma'|_{C'_1}$ be normal
bundles such that a neighborhood of $\Sigma'\cap H_+$, $\Sigma'\cap H_-$ and $\Sigma'\cap \nu_1$ in $H_+$, $H_-$ and $\nu_1$ coincide
with a tubular neighborhood of the zero section of $\pi_+$, $\pi_-$, and $\pi_{\nu_1}$ respectively.
As in Construction~\ref{construccion} we construct a normal bundle $\pi:N\to\Sigma'$
which extends $\pi_+$, $\pi_-$ and $\pi_{\nu_1}$ and is holomorphic over an annulus contained in $\Sigma'_t$ for a certain 
$t\in C$, whose inclusion into $\Sigma'_t$ is a homotopy equivalence.
Consider $\eta$ such that $T_N(\eta)$ is naturally embedded in $T\times\CC^n$ and observe that 
the Milnor number of the restriction of $f$ to the fibres of $\pi$ at the origin of the fibres is constant. 
We have a $\mu$-constant family over $\Sigma'$ for which Assumption A is satisfied (see Corollary~\ref{assumptionA}).

Define 
\[\phi:=(\pi,f):T_N(\eta)\to\Sigma'\times\CC.\]
Consider the restriction
\[\psi|_X=(\tau,f)|_X:X\to U\times D_\delta.\]

As in Construction~\ref{construccion} we may find a smooth hypersurface $S\subset X$, which splits $X$ in two submanifolds with
boundary, whose closures are denoted by $Y$ and $Z$, (being $Y$ the one containing $\Sigma'$) and a homeomorphism
\[\Psi:Y\to (Y\cap H_-)\times [0,1]\] 
whose restrictions to $\Sigma'$ and $Y\setminus\Sigma'$ are smooth, and satisfying $\psi|_Y=\psi|_{Y\cap H_-}\comp q_1\comp\Psi$, 
and $\varphi\comp \pi|_Y=((\pi|_{Y\cap H_-},id)\comp \Psi$, being $q_i$ the projection of $(Y\cap H_-)\times [0,1]$ to the 
$i$-th factor ($i=1,2$),
The inverse $\Psi^{-1}$ may be seen as a flow that integrates a (not necessarily continuous) vector field
$\calY$ defined over $Y$, which is smooth over $\Sigma'$ and $Y\setminus \Sigma'$, is tangent to the fibres of $\psi$, and whose 
integral $\Psi^{-1}$ extends $\varphi^{-1}$ and takes fibres of $\pi$ to fibres of $\pi$.

Our aim now is to extend $\calY$ to a vector field defined over $X$. We work first in $X|_{C'_1}$. Define
\[B_+:=Z\cap X_{int}(H_+,C_1',\delta)\setminus\dot{X}_{int}(\nu_1,C'_1,\delta),\]
\[B_-:=Z\cap {X}_{int}(\nu_1,C'_1,\delta)\setminus (X_{int}(H_+,C_1',\delta).\]
By Ehresmann fibration Theorem and contractibility of the base the restrictions $\psi:B_+\to C'_1\times D_\delta$ and
 $\psi:B_-\to C'_1\times D_\delta$ are trivial fibrations, with fibres cobordisms with boundary, being the boundary
cobordisms the fibres of the restrictions $\psi|_{B_+\cap S}$ and $\psi|_{B_-\cap S}$. All
the boundary cobordisms are simultaneously trivialised by the restrictions $\Psi|_{B_+\cap S}$ and $\Psi|_{B_-\cap S}$.

Applying the procedures of the  proof of Proposition~\ref{pivote} to $B_-$ and $B_+$, we obtain that there is a diffeomorphism
\[\Theta:Z_{C'_1}\to (Z\cap H_-)_{C'_1}\times [0,1]\] 
satisfying 
\begin{enumerate}[(i)]
\item $\Theta|_{S_{C'_1}}=\Psi|_{S_{C'_1}}$, 
\item $\psi|_{Z_{C'_1}}=\psi|_{(Z\cap H_-)_{C'_1}}\comp p_1\comp\Theta$ (being $p_1$ the projection of 
$(Z\cap H_-)_{C'_1}\times [0,1]$ to the first factor),
\item $\Theta(\nu_1|_{C'_1}\cap Z)=(H_-|_{C'_1})\cap Z)\times \{1/2\}$.
\end{enumerate}
The first compatibility implies that $\Theta$ and 
$\Psi$ glue to a homeomorphism $\Xi$ from the union of their domains ($Y\cup Z_{C'_1}$) to the union of the images.
Using a collaring we may assume that $\Xi$ is smooth at $S$. Observe that the cut $\nu_1|_{C'_1}$ is contained in the domain of
definition of $\Xi$ and that 
\[\Xi(\nu_1|_{C'_1})=H_-|_{C'_1}\times \{1/2\}.\]
The inverse $\Xi^{-1}$ may be seen as a flow of a vector field $\calZ$ defined on $Y\cup Z_{C'_1}$.

Observe that the restriction
\[\psi:Z\to U\times D_\delta\]
is a trivial fibration with fibre a trivial cobordism with boundary. Following the procedure used in the proof of Lemma~\ref{extcut}
to extend the vector field $\calY$ to $\calX$,
we can construct a vector field $\calZ'$ in $(Z\setminus\dot{Z})_{C'_1}$ which coincides 
with $\calZ$ at the intersection $S\cup Z_{\partial C'_1}$ of the domains of $\calZ$ and $\calZ'$, 
and whose flow induces a diffeomorphism
\[\Xi':(Z\setminus\dot{Z}|_{C'_1})\to (H_-\cap ((Z\setminus\dot{Z})_{C'_1})\times [0,1].\]
satisfying the analog of property (ii) above. A careful construction of $\calZ'$ (using a collaring) 
yields that $\calZ$ and $\calZ'$ glue to a vector field $\calX$ defined over $X$, which is the desired extension of $\calY$ to $X$.

Step 1 finishes using the flow of $\calX$ as we use the flow $\varphi$ in the proof of Lemma~\ref{extcut}. Step 2 is also as
in Lemma~\ref{extcut}.
\end{proof}

\subsection{Existence of cuts: the non-isolated case}
\label{1extcutii}

Let $f:U\times\CC^n\to\CC$ be a smooth family of holomorphic functions (with $U=(-\epsilon,1+\epsilon)^d$). Consider 
$C=[0,1]^d\subset U$. Let $\rho$ be a distance function associated to a hermitian metric in the trivial vector bundle
$\tau:U\times\CC^n\to U$. Recall that $\psi=(\tau,f)$. Given any subset $B\subset U\times\CC^n$ we denote by $\partial_\tau B$ the 
union $\partial_\tau B:=\cup_{u\in U}\partial B_u$.

\begin{prop}
\label{1usefulconstr}
Let $\theta:C\to (0,\infty)$ be any continuous function.
There exist a positive $\delta$ and a cut $H$ over $C$ with
amplitude $\delta$ such that $X_{int}(H,C,\delta)$ is contained in $B(C,\theta)$.
\end{prop}

\begin{proof}
The proof is completely analogous to the proof of Proposition~\ref{usefulconstr}, using
the fact that for any neighborhood $W$ of the origin adapted to $f_t$ there is a neighborhood $V$ of $t$ in $U$ and a 
positive $\delta$ such that $\partial W\cap \psi^{-1}(V\times D_\delta)$ is a cut for $f$ over $V$ of amplitude $\delta$.
\end{proof}

\subsection{Topological equisingularity} 
Let $C$, $\tau:E\to U$, $f$ and $\psi$ as in Section~\ref{1extcutii}. Let $\theta_1:U\to\RR_+$ be a positive continuous function.
By Proposition~\ref{1usefulconstr} there exists $\delta_1>0$ and a cut $H_1$ for $f$ over $U$ of 
amplitude $\delta_1$ such that $X:=X_{int}(H_1,C,\delta_1)$ is contained in $B(C,\theta_1)$. Define $X^*:=X\setminus f^{-1}(0)$.\\

\textbf{Construction} ($\dagger$). Let $H_2$ be a cut for $f$ over $C$ of amplitude $\delta_2<\delta_1$ such that $H_2\prec H_1$.
Consider $Z_1:=X_{int}(H_1,C,\delta_2)\setminus\dot{X}_{int}(H_2,C,\delta_2)$.
As in Construction~\ref{construccion} we obtain a normal
bundle $\pi_1:M_1\to Z_1\cap\Sigma$ such that $T_\eta(M_1)$ embeds in $Z_1$ as a tubular neighborhood of $\Sigma\cap Z_1$ 
for a positive $\eta$. Using Proposition~\ref{pivote} it is easy to construct a vector field $\calZ_1$ in $Z_1$ which is smooth 
outside $\Sigma$, smooth and tangent to $\Sigma$ at $\Sigma\cap Z$, and such that its flow
induces a homeomorphism
\[\Theta:Z_1\to (Z_1\cap H_1)\times [0,1]\]
satisfying
\begin{enumerate}[(a)]
\item we have $\psi=\psi|_{Z_1\cap H_1}\comp p_1\comp\Theta$, where $p_i$ is the projection of $(Z_1\cap H_1)\times [0,1]$ to the 
$i$-th factor,
\item we have the equalities
\[\Theta(Z_1\cap H_1)=(Z_1\cap H_1)\times\{0\}\]
\[\Theta(Z_1\cap H_2)=(Z_1\cap H_1)\times\{1\},\]
\item for $\eta$ sufficiently small,
the diffeomorphism given by the restriction of the flow to $T_\eta(M_1)$ for a fixed time $s$ takes the fibres 
of $\pi_1$ to fibres of $\pi_1$.
\end{enumerate}

\begin{lema}
\label{1vectorfield}
There exists a (not necessarily continuous) vector field $\calX$ in 
$X\setminus C$ with the following properties:
\begin{enumerate}
\item its restrictions to $X\setminus\Sigma$ and to $\Sigma\setminus C$ are smooth
\item it is tangent to the fibres of $\tau$, 
\item there exists a vector field $\calW$ in $D_\delta^*$, which is radial, pointing to the origin and of modulus 
$|\calW(z)|\leq |z|^2$, such that $df(\calX)(x)=\calW(f(x))$ for any $x\in X\setminus f^{-1}(0)$,
\item the vector field $\calX$ is tangent to $f^{-1}(0)$ outside $\Sigma$ and it is tangent to $\Sigma$ in $\Sigma$,
\item any integral curve converges to the
origin of the fibre of $\tau$ in which it lies in positive infinite time,
\item the flow of $\calX$ is continuous.
\end{enumerate}
\end{lema}
\begin{proof}
As in the proof of Lemma~\ref{vectorfield} we construct $\calX$ as the amalgamation two vector fields $\calY$ and $\calZ$. 

There is a continuous function $\theta_2:C\to\RR$ such that
$B(C,\theta_2)$ is included in $X_{int}(H_1,C,\delta_1)$. By Proposition~\ref{1usefulconstr}
there exists $\delta_2$ satisfying $\delta_1>\delta_2>0$ and a cut $H_2$ for $f$ over $C$ of 
amplitude $\delta_2$ such that $X_2:=X_{int}(H_1,C,\delta_2)$ is contained in $B(C,1/2\theta_2)$.
We iterate this procedure to obtain an infinite sequence of constants $\delta_i$ and cuts $H_i$ over $C$ such that the sets $X_i:=X_{int}(H_i,C,\delta_i)$ form a nested sequence
\begin{equation}
\label{1nestedset}
X=X_1\supset X_2\supset ...\supset X_i\supset ...
\end{equation}
of closed neighborhoods of $C$ such that $\cap_{i=1}^\infty X_i$ is equal to the zero section $C$.

Let $Z_i:=(X_i\setminus\dot{X}_{i+1})\cap f^{-1}(D_{\delta_{i+1}})$. Let $\pi_i:M_i\to\Sigma\cap Z_i$ be the normal bundle
appearing in Construction~($\dagger$). If the normal bundles $\pi_i$ are chosen carefully enough they glue to a normal bundle
$\pi:M\to\Sigma\setminus C$. Let $\xi:\Sigma\setminus C\to\RR_+$ be a smooth function such that 
$T_M(\xi)$ embeds as a neighborhood of $\Sigma\setminus C$ in $X\setminus C$.
Construction ($\dagger$) produces a vector field $\calZ_i$ over $Z_i$. A partition of unity argument glues the vector
fields $\calZ_i$ to a vector field $\calZ$ defined on $Z:=\cup_{i=1}^\infty Z_i$, smooth in $Z\setminus\Sigma$, tangent to the fibres
of $\psi$ at their smooth points, and such that $\calZ|_\Sigma$ is smooth and tangent to $\Sigma$. 
Observe that we have $df(\calZ)=0$. Let 
\[\Psi:W\to T_M(\xi)\]
be the flow of $\calZ|_{T_M(\xi)}$, where $W\subset T_M(\xi)\times\RR$ is the maximal domain where it is defined. By the property
$(c)$ of the flow of $\calZ_i$ (see Construction ($\dagger$)), for any $(x,u)\in W$ the conmuting relation 
\begin{equation}
\label{conmuting}
\pi(\Psi(x,u))=\Psi(\pi(x),u)
\end{equation}
holds.

Consider in $C\times D_\delta^*$ the vector field $\calV'$ characterised by being tangent to the fibres of the projection from 
$C\times D_\delta^*$ to the first factor, and a lift by the projection from $C\times D_\delta^*$ to the second factor of the vector
field $\calV$ in $D_\delta^*$ which is radial, pointing at the origin and of modulus $||\calV(z)||=||z||^{2}$.
As $\psi:X^*\to U\times D_\delta^*$ is submersive we can define the vector field $\calY$ in 
$X^*$ to be a lifting of $\calV'$ by the mapping $\psi$. Moreover, as the restriction of $f$ to the fibres of $p:T_M(\xi)\to\Sigma$
only has critical points at $\Sigma$, (perhaps having to shrink $\xi$) we may construct $\calY$ such that its restriction to 
$T_M(\xi)$ is tangent to the fibres of $\pi:T_M(\xi)\to\Sigma$. 

Let $\rho_1:D_\delta\to [0,1)$ be a smooth function vanishing at $0$ and positive in $D_\delta^*$. We choose $\rho_1$ small enough
that the modulus $||\rho_1(f(z))\calY(z)||$ converges to $0$ as $||f(z)||$ approaches $0$.
Let $\rho_2:U\to\RR$ be an smooth function with support contained in the interior of 
$Z$ and which is identically $1$ in a neighborhood of $f^{-1}(0)\cap Z$ in $Z$. Define the vector field $\calX:=(\rho_1\comp f)\calY+\rho_2\calZ$ on $X$.

As $\calX|_{f^{-1}(0)}$ coincides with $\calZ|_{f^{-1}(0)}$ it is clear that any 
integral curve in $f^{-1}(0)$ converges to the origin in positive infinite time. It is also clear that the restrictions of 
$\calX$ to $Y\setminus\Sigma$ and to $\Sigma\setminus C$ are smooth, and that the later is tangent to $\Sigma\setminus C$. All the properties required to $\calX$ are clear except (3), (5) and (6).

Property (3) is true taking $\calW:=\rho_1\calV$.

Observe that the restriction $\calY|_{H_1\cap Y}$ is tangent to $H_1$, and that $\calZ$ points into $X$ at any $y\in Z\cap H_1$. 
On the other hand $df(\calX)$ is radial and pointing to the origin. This shows that no integral curve of $\calX$ can go out 
of the domain $X$ in positive time. As $df(\calX)=\rho_1\calV$ and $\calV$ is radial, pointing to the origin, and of
module small enough that any of its integral curves converges to the origin in positive infinite time, we deduce that any 
integral curve of $\calX$ not lying in $f^{-1}(0)$ is defined in positive infinite time. As $\calX$ coincides with $\calZ$ in 
$f^{-1}(0)$ it is also clear that the integral curves contained in $f^{-1}(0)$ are defined in positive infinite time.

The continuity of the flow of $\calX$ 
\[\Phi:(X\setminus C)\times [0,\infty)\to X\]
is clear outside $(\Sigma\setminus C)\times [0,\infty)$, since $\calX$ is smooth outside $\Sigma$. 
Observe that from the Relation~(\ref{conmuting}) and from the tangency of $\calY|_{T_M(\xi)}$ to the fibres 
of $\pi:T_M(\xi)\to\Sigma$ we obtain that if $\Phi:W'\to T_M(\xi)$ is the flow of $\calX|_{T_M(\xi)}$ 
(where $W'$ is its maximal domain of definition), then we have the conmutation relation
\begin{equation}
\label{conmuting2}
\pi(\Phi(x,u))=\Phi(\pi(x),u)
\end{equation}
for any $(x,u)\in W'$.
The continuity of $\Phi$ at any point of $(\Sigma\setminus C)\times [0,\infty)$ follows easily using this relation and the 
fact that the modulus $||\rho_1(f(z))\calY(z)||$ converges to $0$ as $||f(z)||$ approaches $0$. This shows (6).

As the vector field $\calX$ is smooth outside $\Sigma$ the only accumulation points of an integral curve of $\calX$ as time 
tends to $+\infty$ can be in $\Sigma$. As any integral curve of $\calX$ inside $\Sigma$ converges to the origin in positive infinite
time, from the continuity of the flow of $\calX$ we deduce that the 
only accumulation point of an integral curve of $\calZ$ can be the origin of the fibre of $\tau$ where it is contained.
This proves Property (5).
\end{proof}

Here is our main topological equisingularity theorem

\begin{theo}
\label{porfin}
Let $f:C\times \CC^n\to\CC$ be a family of holomorphic germs at the origin, with $1$-dimensional critical set and smoothly parametrised
over a cube $C$.
Suppose that it is equisingular at the critical set. 
Let $t$ be any point of $C$.
Use the notations considered in this section.
There exists a homeomorphism
\begin{equation}
\label{1homeoisol}
\Psi:X\to C\times X_t
\end{equation}
such that
\begin{enumerate}
\item we have $\tau=p_1\comp\Psi$, where $p_i$ is the projection of $C\times X_t$ to the 
$i$-th factor ($i=1,2$),
\item we have $f|_{X_t}\comp p_2\comp\Psi=f$, 
\item the restriction of $\Psi$ to $X\setminus\Sigma$, and to $\Sigma\setminus C$ is smooth.
\end{enumerate}
\end{theo}
\begin{proof}
The first step is to construct the restriction of $\Psi$ to $\partial_\tau X$. We have the decomposition
$\partial_\tau X=(X\cap f^{-1}(\partial D_\delta))\cup (H_1\cap f^{-1}(D_\delta)$. The restriction of $\Psi$ to 
$H_1\cap f^{-1}(D_\delta)$ can be 
constructed applying Proposition~\ref{pivote} to the cuts $H_2\prec H_1$ (see Construction~($\dagger$). The extension to 
$X\cap f^{-1}(\partial D_\delta)$ is easy to obtain using that
\[\psi:X\cap f^{-1}(\partial D_\delta)\to C\times\partial D_{\delta}\] 
is a locally trivial fibration.

After this the proof is completely analogous to the proof 
of Theorem~\ref{globalisationLR} replacing the reference to Lemma~\ref{vectorfield} by a reference to Lemma~\ref{1vectorfield}.
\end{proof}

\section{Families with constant L\^e numbers}

Let $f:T\times\CC^n\to\CC$ a holomorphic function, where $T$ is a connected complex manifold.
Unless we state the contrary we suppose that the critical set of $f_t$ has dimension $1$ at the origin for any $t\in T$. 
Let $\mathbf{Z}:=\{z_1,...,z_n\}$ be a coordinate system of $\CC^n$. 
The symbol $\lambda^i_{f,\mathbf{Z}}(x)$ denotes the $i-th$ L\^e number of a holomorphic function $f$ at the point $x$ with respect
to $\mathbf{Z}$.
Our aim in this section is to prove the following theorem:

\begin{theo}\label{wanted}
If any of the following conditions hold:
\begin{enumerate}
\item there is a coordinate system $\mathbf{Z}$ such that the L\^e numbers at the origin of $f_t$ 
with respect to $\mathbf{Z}$ are defined and independent of $t$,
\item the generic L\^e numbers at the origin of $f_t$ are independent of $t$,
\end{enumerate}
then the embedded topological type of $f_t$ at the origin and the the diffeomorphism type of the 
Milnor fibration is independent of $t\in T$. Moreover, if $T$ is (diffeomorphic to) a cube, 
the restriction of the family over $T$ is topologically $R$-equisingular.
\end{theo}

Because of the results of the previous sections to prove the Theorem it is enough to show that the family
is equisingular at the singular set. Thus we are reduced to prove:

\begin{theo}\label{leconsteqsing}
If any of the conditions of the previous theorem hold then $f_t$ is equisingular at the critical set.
\end{theo}

Since equisingularity at the critical set is local in the base, we can assume that $T$ is an open 
neighborhood of the origin in an affine space, and it is enough to prove that $f$ is equisingular at the 
critical set in a small neighborhood $V$ of the origin in $T$. 

A generic coordinate system $\Zb$ for $f_{t_0}$ satisfies that the L\^e 
numbers at the origin of $f_t$ with respect to it are defined for any $t$ close enough to $t_0$. 
Moreover, by the lexicographical upper semicontinuity of the L\^e numbers (\cite{Ma3}~Corollary~4.16)
the constancy of the generic L\^e numbers
implies the constancy of the L\^e numbers of $f_t$ with respect to $\Zb$, for $t$ close enough to $t_0$.
Therefore Condition~(2) implies Condition~(1), at least locally in the base.
  
By the definition of equisingularity at the critical set it is clear that
(by slicing), it is enough to assume that $T$ has complex dimension $1$. Thus we make the assumption 
$T=\CC$. To any coordinate system $\mathbf{Z}=\{z_1,...,z_n\}$ for $(\CC^n,O)$ we associate the 
coordinate system $\mathbf{Z'}:=\{t,z_1,...,z_n\}$ of $\CC\times\CC^n$. 

Given $A\subset\CC\times\CC^n$ we consider the notation $A_t:=A\cap (\{t\}\times\CC^n$).

Let $\Sigma_f:=V(\partial f/\partial t, \partial f/\partial z_1,...,\partial f/\partial z_n)$ denote the critical set of $f$ and let
$\Sigma:=V(\partial f/\partial z_1,...,\partial f/\partial z_n)$ be the union of the critical sets of the $f_t$'s.

\begin{lema}
\label{thom}
Let $f_t$ be a family with critical set not necessarily $1$-dimensional: let $s:=dim_O(\Sigma_0)$. 
Suppose that for small $t$ the L\^e numbers
$\lambda^i_{f_t,\mathbf{Z}}(O)$ are defined and independent on $t$ for any $0\leq i\leq s$. Then the sets $\Sigma_f$ and $\Sigma$ are equal
in a neighborhood of the origin. 
\end{lema}
\begin{proof}
We only have to show $\Sigma\subset\Sigma_f$. For this we show that for any 
$t$ small enough the set $(\Sigma_f)_t$ contains $\Sigma_t$.
If the inclusion does not hold there exists a sequence $y_n=(t_n,x_n)$ converging to $(O,0)$ such that $x_n$ belongs to $\Sigma_{t_n}\setminus\Sigma_f$.
Hence $V(f-f(y_n))$ is smooth at $y_n$ and the tangent 
space $T_{y_n}V(f)$ is equal to $V(t-t_n)$. Consequently the limit of the tangent spaces is
$V(t)$.

On the other hand, by \cite{Ma3}, Theorem~6.5, $\CC\times\{O\}$ satisfies Thom's $A_f$ condition
with respect to the open stratum at the origin, which contradicts the fact that the limit of 
tangent spaces is $V(t)$.
\end{proof}

A consequence is that the first polar variety of $f$ is empty, since, it is defined by
\begin{equation}
\label{1pF}
\Gamma^1_{f,\Zc}=\overline{V(\partial f/\partial z_1,...,\partial f/\partial z_n)\setminus\Sigma_f}.
\end{equation}
We deduce the following vanishings for $t$ small enough 
(the L\^e varieties and L\^e cycles computations that follow are done using \cite{Ma3}, Chapter~1):
\begin{equation}
\label{0lF}
\gamma^1_{f,\Zc}(t,O)=0\quad\quad\Lambda^0_{f,\Zc}=0\quad\quad\lambda^0_{f,\Zc}(t,O)=0.
\end{equation}

Let $\Sigma=\cup_{i=1}^k\Sigma_{i}$ and $\Sigma_{0}=\cup_{i=1}^r\Sigma_{0,i}$ be
decompositions of
$\Sigma$ and $\Sigma_{0}$ in irreducible components. Reordering conveniently we can find 
numbers $1\leq k_1\leq k_2\leq k$ such that 
\begin{enumerate}
\item The component $\Sigma_{i}$ is $2$-dimensional for $1\leq i\leq k_1$.
\item The component $\Sigma_{i}$ is $1$-dimensional and contained in $V(t)$ for 
$k_1+1\leq i\leq k_2$.
\item The component $\Sigma_{i}$ is $1$-dimensional and not contained in $V(t)$ for 
$k_2+1\leq i\leq k$.
\item There exist numbers $0=r_0\leq r_1\leq ...\leq r_{k_2}= r$ such that
\begin{equation}
\Sigma_{j}\cap V(t)=\cup_{i=r_{j-1}+1}^{r_j}\Sigma_{0,i}
\end{equation}
for any $1\leq j\leq k_2$.
\end{enumerate}


Now we assume that all the L\^e numbers of $f_t$ at the origin with respect to $\Zb$ are defined for 
any $t$.

In the following discussion we draw further consequences from the fact that the L\^e numbers $\lambda^i_{f_t,\mathbf{Z}}(O)$ are
independent on $t$ for any $0\leq i\leq s=1$.
Let $X, Y\subset U$ be closed analytic subspaces of an open subset $U$ of $\CC^n$. We denote by $\overline{X\setminus Y}$ 
the scheme theoretical closure of $X\setminus Y$ in $U$, and by $[X]$ the cycle associated with the scheme $X$
(see \cite{Ma3}, Chapter~1).

Since the first L\^e number at the origin of $f_t$ with respect to $\Zb$ is defined, and $\Sigma_{t}$ is
$1$-dimensional, all the irreducible components of $\Gamma^2_{f_t,\Zb}$ are of 
the expected dimension $2$ for any $t$ (see~\cite{Ma3}, pages 11-14). As $\Sigma_{f_t}$ is of dimension $1$ and, by definition,
\begin{equation}
\Gamma^2_{f_t,\Zb}=\overline{V(\partial f_t/\partial z_3,...,\partial f_t/\partial z_n)\setminus\Sigma_{f_t}},
\end{equation}
we conclude that all the irreducible components of
$V(\partial f_t/\partial z_3,...,\partial f_t/\partial z_n)$ are $2$-dimensional. Hence 
$(\partial f_t/\partial z_3,...,\partial f_t/\partial z_n)$ is a regular sequence and the analytic subscheme 
$V(\partial f_t/\partial z_3,...,\partial f_t/\partial z_n)$ has no embedded components. Thus we have
\begin{equation}
\label{2pf}
\Gamma^2_{f_t,\Zb}=V(\partial f_t/\partial z_3,...,\partial f_t/\partial z_n).
\end{equation}

The third polar variety of $f$ is 
\begin{equation}
\label{3pF}
\Gamma^3_{f,\Zc}=\overline{V(\partial f/\partial z_3,...,\partial f/\partial z_n)\setminus\Sigma_f}=V(\partial f/\partial z_3,...,\partial f/\partial z_n).
\end{equation}
The first equality is by definition. From the fact that all the irreducible components of 
$V(\partial f_t/\partial z_3,...,\partial f_t/\partial z_n)$ are $2$-dimensional and the Theorem on the dimension of the fibres of a morphism we deduce that all the irreducible components of
$V(\partial f/\partial z_3,...,\partial f/\partial z_n)$ have at most dimension $3$, which is on the other hand its minimal possible 
dimension given the number of equations defining the set. Hence 
$(\partial f/\partial z_3,...,\partial f/\partial z_n)$ is a regular sequence and the analytic subscheme 
$V(\partial f/\partial z_3,...,\partial f/\partial z_n)$ has no embedded components. This, together with the fact that $\Sigma$
is $2$-dimensional implies the second equality.


Consider a decomposition in irreducible components of 
$\Gamma^3_{f,\Zc}\cap V(\partial f/\partial z_2)=V(\partial f/\partial z_2,...,\partial f/\partial z_n)$:
\begin{equation}
\label{decom}
V(\partial f/\partial z_2,...,\partial f/\partial z_n)=(\bigcup_{i=1}^{k_1}X^i)\bigcup(\bigcup_{i=1}^dY^i),
\end{equation}
where $X^i$ coincides as a set with the component $\Sigma_{i}$ of $\Sigma$. As any $Y_i$ is (at least) $2$-dimensional and 
the components $\Sigma_i$ are $1$-dimensional for $i>k_1$, no $Y^i$ is a component
of $\Sigma$. We have 
\begin{equation}
\label{2pF}
\Gamma^2_{f,\Zc}=\overline{\Gamma^3_{f,\Zc}\cap V(\partial f/\partial z_2)\setminus\Sigma}=\bigcup_{i=1}^dY^i,
\end{equation}
\begin{equation}
\label{2lF}
\Lambda^2_{f,\Zc}=[\Gamma^3_{f,\Zc}\cap V(\partial f/\partial z_2)]-[\Gamma^2_{f,\Zc}]=\sum_{i=1}^{k_1}[X^i].
\end{equation}

Now we want to compute the first L\^e cycle of $f_t$ for a certain value of the parameter $t$. By equation~(\ref{2pf}) we have
\begin{equation}
\Gamma^2_{f_t,\Zb}\cap V(\partial f_t/\partial z_2)=V(\partial f_t/\partial z_2,...,\partial f_t/\partial z_n).
\end{equation}
Hence, taking into account Equation~(\ref{decom}), we have
\begin{equation}
\label{decomt}
[\Gamma^2_{f_t,\Zb}\cap V(\partial f_t/\partial x_2)]=(\sum_{i=1}^{k_1}[X^i_t])+(\sum_{i=1}^d[Y^i_t]).
\end{equation}
Clearly $(\cup_{i=1}^{k_1}X^i_t)$ is included in $\Sigma_{f_t}$, and therefore
\begin{equation}
\label{desigualdad1}
\Lambda^1_{f_{t_0},\Zb}\geq\sum_{i=1}^{k_1}[X^i_{t_0}]=(\sum_{i=1}^{k_1}[X^i])\centerdot [V(t-t_0)].
\end{equation}

\begin{lema}
\label{firstle}
The previous inequality is actually an equality for any $t_0\in T$:
\begin{equation}
\label{igualdad1}
\Lambda^1_{f_{t_0},\Zb}=\sum_{i=1}^{k_1}[X^i_{t_0}]=(\sum_{i=1}^{k_1}[X^i])\centerdot [V(t-t_0)].
\end{equation}
Moreover, the components $\Sigma_{k_1+1}$,...,$\Sigma_{k_2}$ do not appear. In other words, we have $k_1=k_2$.
\end{lema}
\begin{proof}
We claim that no component of $Y^i_{t_0}$ is contained in $\Sigma_{f_{t_0}}$ for any $i\leq d$ and any $t_0$. 
The claim implies easily both assertions of the Lemma.

By Lemma~\ref{thom} it is enough to show that no component of $Y^i_{t_0}$ is contained in $\Sigma$.
This is clear (after possibly shrinking the base of the family) for $t_0\neq 0$. Thus equality~(\ref{igualdad1}) holds in this case.

Inequality~(\ref{desigualdad1}) implies that the first L\^e number of $f_{t_0}$ satisfies
\begin{equation}
\label{desigualdad2}
\lambda^1_{f_{t_0},\Zb}\geq\mult_{[O]}((\sum_{i=1}^{k_1}[X^i])\centerdot [V(t-t_0)]\centerdot [V[z_1]]
=\mult_{[O]}((\sum_{i=1}^{k_1}[X^i\cap V(z_1)])\centerdot [V(t-t_0)].
\end{equation}
This inequality becomes an equality if and only if our claim is true. This is the case if $t_0\neq 0$.

As $\mult_{[O]}((\sum_{i=1}^{k_1}[X^i\cap V(z_1)])\centerdot [V(t)]$ is upper semicontinuous in $t$ and the L\^e number 
$\lambda^1_{f_t,\Zb}$ is independent on $t$ we conclude that Inequality~(\ref{desigualdad2}) must be an equality also for $t_0=0$. \end{proof}

\begin{lema}
\label{leponumF}
The L\^e numbers $\lambda^i_{f,\Zc}(t,O)$ are defined for $i\leq 2$ and any $t$.
\end{lema}
\begin{proof}
We work for $t=0$. The case $t\neq 0$ is analogous.
We have already checked that $\lambda^0_{f,\Zc}(O)$ is defined (and zero).

Two consequences of Lemma~\ref{firstle} are the equalities 
\begin{equation}
\label{loc0}
\Lambda^1_{f_0,\Zb}=\sum_{i=1}^{k_1}[X^i\cap V(t)],
\end{equation}
\begin{equation}
\label{loc1}
\Gamma^1_{f_0,\Zb}=\sum_{i=1}^{d}[Y^i\cap V(t)].
\end{equation}

The definedness of $\lambda^1_{f_0,\Zb}(O)$ implies that  
$X^i\cap V(t,x_1)$ is $0$-dimensional at the origin for any $i$. Hence, taking into account~(\ref{2lF}), 
we obtain that $\lambda^2_{f,\Zc}(0,O)$ is defined.

It only remains to be shown that $\lambda^1_{f,\Zc}(O)$ is defined. The first L\^e cycle is 
\[\Lambda^1_{f,\Zc}=[\Gamma^2_{f,\Zc}\cap V(\partial f/\partial z_1)]-[\Gamma^1_{f,\Zc}]\]
\[=\sum_{i=1}^{d}[Y^i\cap V(\partial f/\partial z_1)].\]
We only need to prove that $Y^i\cap V(\partial f/\partial z_1)\cap V(t)$ is $0$-dimensional at the origin. By Identity~(\ref{loc1})
the last set is equal to $\Gamma^1_{f_0,\Zb}\cap V(\partial f_0/\partial z_1)$. We have to show that
$\Gamma^1_{f_0,\Zb}\cap V(\partial f_0/\partial z_1)$ is $0$-dimensional. If this is not the case then
$\Gamma^1_{f_0,\Zb}$ is contained in $V(\partial f_0/\partial z_1,...,\partial f_0/\partial z_n)$ and meets the origin.
Hence it is contained in $V(f_0)$, but this is impossible by definition of the relative polar varieties.
\end{proof}

\begin{lema}
\label{constconst}
The components $\Sigma_{k_2+1}$,...,$\Sigma_{k}$ do not appear. In other words, we have $k_2=k$
\end{lema}
\begin{proof}
Applying Proposition~1.21 of \cite{Ma3}, the vanishing~(\ref{0lF}) and the constancy of the L\^e numbers of $f_t$ 
we obtain that for any $t$ small enough $\lambda^1_{f,\Zc}(t,O)$ is constant.

Suppose $k_2\neq k$. As $\Sigma_{k_2+1}$ is a $1$-dimensional irreducible component of the singular locus the cycle $[\Sigma_{k_2+1}]$ is a positive summand of $\Lambda^1_{f,\Zc}$. 
As it goes through $(0,O)$ but not through $(t,O)$ for $t$ non-zero,
we have that $\lambda^1_{f,\Zc}(0,O)$ is strictly bigger than $\lambda^1_{f,\Zc}(t,O)$, which is a contradiction.
\end{proof}

Let $\Sigma_{red}$ denote the analytic subset $\Sigma$ with its reduced structure. We have
\begin{cor}
\label{flatness}
The analytic subset $\Sigma$ admits a decomposition in irreducible components at the origin of the form 
$\Sigma:=\cup_{i=1}^k\Sigma_{i}$ such that each component is $2$-dimensional at the origin. Moreover the restriction 

\begin{equation}
\label{proyres}
\pi|_{\Sigma_{red}}:\Sigma_{red}\to\CC
\end{equation}
is flat.
\end{cor}
\begin{proof}
We have proved everything except flatness. 
Each of the components of the primary decomposition of $\Sigma_{red}$ is an irreducible component, and all of them map dominant 
by $\pi|_{\Sigma_{red}}$ to the $1$-dimensional target $\CC$. It is well known that this implies that $\pi|_{\Sigma_{red}}$ is
flat.
\end{proof}

Now we study the geometry of the fibres of the mapping~(\ref{proyres}) in a neighborhood of 
the origin. To obtain a lighter notation we denote $\Sigma_{red}$ by $S$.

Consider the mapping
\begin{equation}
\label{alpha1}
\alpha:=(z_1,t):S\to\CC^2.
\end{equation}
As the origin $(0,O)$ is an isolated point of $\alpha^{-1}(O)$, there is a neighborhood
$U$ of $(0,O)$ in $S$ such that, for any $\epsilon$, $\eta$ positive and small enough, the restriction
\begin{equation}
\label{alpha}
\alpha|_{E}:E:=\alpha^{-1}(D_\epsilon\times D_\eta)\cap U\to D_\epsilon\times D_\eta
\end{equation}
is a finite mapping, the space $E$ is a neighborhood of the origin in $S$ and $\alpha|_E^{-1}(0,0)=\{(0,O)\}$.
We will denote $\alpha|_E$ simply by $\alpha$.
In the sequel $\epsilon$ and $\eta$ are chosen always so that $\alpha$ is finite.

We choose $\epsilon$ and $\eta$ so small 
that the decompositions $E=\cup_{i=1}^k E_i$ and $E_0=\cup_{i=1}^{r}E_{0,i}$ in irreducible components coincide at the origin
with the previously given decompositions; in other words, as germs of reduced sets at the origin we have
$E_{i}=\Sigma_i$ and $E_{0,i}=\Sigma_{0,i}$. We will always choose $\epsilon$ small enough that any two different components 
of $E_0$ only intersect at the origin.

It is well known that for $\epsilon$ and $\eta$ sufficiently small the mapping 
\begin{equation}
\label{fibrsing}
\pi|_{E}:E\to D_\eta
\end{equation}
is a topological locally trivial fibration over $D_\eta^*$.

We consider a decomposition of $E_{t}$ in irreducible components $E_{t,1},...,E_{t,r_t}$ for any $t$. 
The topological local triviality of $\pi|_E$ over $D_\eta^*$ implies that the number $r_{gen}:=r_t$ is independent on $t$ as long as
$t\in D_\eta^*$. For any $t\in D_\eta^*$ we consider a numbering $E_{t,1},...,E_{t,r_{gen}}$ of the irreducible components of $E_t$.
There is a mapping $\beta_t:\{1,...,r_{gen}\}\to\{1,...,k\}$ defined by the relation 
$E_{t,j}\subset E_{\beta_t(j)}$. We also define a mapping $\beta_{0}:\{1,...,r\}\to\{1,...k\}$ 
by the relation $E_{0,j}\subset E_{\beta_{0}(j)}$.

Let $\mu_{t,j}$ denote the generic transversal Milnor number of $f_t$ at a generic point of $E_{t,j}\setminus\{O\}$. 
Observe that there is a finite set $I_t\subset D_\epsilon$ such that, for any
$s\in D_\epsilon\setminus I_t$, the hyperplane $V(z_1-s)$ meets $E_{t}$ transversely at points with transversal Milnor number 
$\mu_{t,j}$. By Lemma~\ref{gentranstype}, for any $s\in D_\epsilon\setminus I_t$,
the Milnor number of any of the isolated singularities of $f_t|_{V(z_1-s)}$ contained in $E_{t,j}$ 
is equal to $\mu_{t,j}$.
Moreover, for any $x\in V(z_1-s)\cap E_t$, we have the equality
\begin{equation}
\label{legen}
I_x(\Lambda^1_{f_t,\Zb},V(z_1-s))=\mu_{t,i};
\end{equation}
indeed, since $f_t$ has smooth critical set at $s$ and $V(z_1-s)$ meets $\Sigma_t$ transversely in $x$, the first L\^e number
$\lambda^1_{f_t,\Zb}(x)$ is equal to the generic transversal Milnor number.
On the other hand, by definition, the L\^e number is equal to 
the intersection multiplicity $I_x(\Lambda^1_{f_t,\Zb},V(z_1-s))$.

Given any component $E_i$ we define a family of germs with isolated singularities parametrised over 
it, by imposing that, for any $x\in E_i$, we assign the singularity at $x$ defined by the 
restriction $f_{\pi(x)}|_{V(z_1-z_1(x))}$. Let $\mu_i$ be the generic Milnor number of the family. Observe that, by the upper
semicontinuity of the Milnor number, if $E_{t,j}$ is contained in $\Sigma_i$ then $\mu_i\leq \mu_{t,j}$. As $\mu_i$ and 
$\mu_{t,j}$ are the generic values of the Milnor number of the hyperplane sections $V(z_1-s)$ at $E_i$ and $E_{t,j}$ respectively,
we obtain that for any generic $t$ we have the equality $\mu_{i}=\mu_{t,j}$: the subset of 
$E_i$ in which the hyperplane section has Milnor number strictly bigger than $\mu_i$ is closed in the Zariski topology. 
We may express the last equality as 
\[\mu_{t,j}=\mu_{\beta_t(j)}\]
for $t$ generic.

By Sard's Theorem the mapping $\pi|_E$ is generically submersive at the regular locus of $E$.

\begin{prop}
\label{reduced}
The following assertions hold:
\begin{enumerate}
\item If $\epsilon$ is small enough, the analytic set $E_0\setminus\{O\}$ is smooth and $\pi|_E$ is submersive at $E_{0}\setminus\{O\}$.
\item For any $j\in\{1,...,r_{0}\}$ we have $\mu_{0,j}=\mu_{\beta_0(j)}$. 
\end{enumerate}
\end{prop}

\begin{proof}
Clearly, if $\epsilon$ is small enough, the only singularity of $E_0$ is at the origin.

The first L\^e number $\lambda^1_{f_{t_0},\Zb}(O)$ is $I_O(\Lambda^1_{f_{t_0},\Zb},V(z_1))$ by definition. In Lemma~\ref{firstle}
we have proved that
$\Lambda^1_{f_{t_0},\Zb}$ is equal to $\sum_{i=1}^{k_1}[X^i_{t_0}]=\sum_{i=1}^{k_1}[X^i]\centerdot [V(t-t_0)]$ for any $t_0$. 
By conservation of number in 
intersection theory we have that for a certain $\epsilon$ small enough, there exist positive $\xi$ and $\eta$
such that if $s\in D_{\xi}$ and $t\in D_{\eta}$, then
\begin{equation}
\label{intmullt1}
\lambda^1_{f_0,\Zb}(O)=\sum_{x\in X_t\cap V(z_1-s)}I_x([X_t],V(z_1-s)),
\end{equation}
By Equality~(\ref{legen}), if $s\in D_\epsilon\setminus I_t$, this quantity is equal to the sum 
\begin{equation}
\label{sumaaux}
\sum_{x\in V(f_t,z_1-s)}\mu(f_t|_{V(z_1-s)},x)
\end{equation}
of the Milnor numbers of the singularities of the restriction of $f_t$ to $V(z_1-s)$ which are contained in the zero set 
$V(f_t,z_1-s)$. 

If either a component of $E_0$ is contained in the singular locus of $E$, or $\pi|_E$ is not generically submersive at it,
then $E_0$ is contained in the ramification locus of $\alpha$, and therefore $D_\epsilon\times\{0\}$ is a component
of the branching locus of $\alpha$. Hence, if $t\in D_\eta^*$ and $\eta$ is small enough the set $D_\epsilon\times\{t\}$
is not included in the branching locus of $\alpha$. Then, if $s$ does not belong to the finite set $I_0\cup I_t$, the 
cardinality of $X_0\cap V(z_1-s)$ is strictly smaller than the cardinality of $X_t\cap V(z_1-s)$ for $t\neq 0$. As the cardinality
of $X_t\cap V(z_1-s)$ is the number of singularities of $f_t|_{V(z_1-s)}$ contained in $f_t|_{V(z_1-s)}^{-1}(0)$, 
by the non-splitting result
of L\^e-Lazzery we have the strict inequality
\begin{equation}
\label{ineqaux}
\sum_{x\in V(f_t,z_1-s)}\mu(f_t|_{V(z_1-s)},x)<\sum_{x\in V(f_0,z_1-s)}\mu(f_0|_{V(z_1-s)},x),.
\end{equation}
which is impossible since both quantities are equal to $\lambda^1_{f_0,\Zb}(O)$ by Equalities~(\ref{intmullt1}),~(\ref{sumaaux}).
This proves that $\pi$ is generically submersive at any component of $E_0$. Shrinking $\epsilon$ we obtain assertion (1).

Fix $s_0\neq 0$ small enough so that $V(z_1-s_0)$ meets $E_0$
transversely in finitely many points. Then there exits $\xi$ small enough that $V(z_1-s_0)$ meets $E_t$ transversely 
for any $t\in D_\xi$, and hence 
\[\alpha|_{\alpha^{-1}(\{s_0\}\times D_\xi)}:\alpha^{-1}(\{s_0\}\times D_\xi)\to\{s_0\}\times D_\xi\]
is finite and etale. Thus, the space $\alpha^{-1}(\{s_0\}\times D_\xi)$ splits in a finite number of disks and 
we have a natural bijection between $\Sigma_{f_t}\cap V(z_1-s_0)$ and 
$\Sigma_{f_0}\cap V(z_1-s_0)$ induced by incussion in $\alpha^{-1}(\{s_0\}\times D_\xi)$.

Choose $x\in \Sigma_{f_t}\cap V(z_1-s_0)$ for a certain $t\in D_\xi^*$, 
and let $y$ be the corresponding point in 
$\Sigma_{f_0}\cap V(z_1-s_0)$ by the bijection.
Suppose that $x$ belongs to $E_{t,j_1}$ and that $y$ belongs to $E_{0,j_2}$. Clearly 
$E_{t,j_1}$ and $E_{0,j_2}$ belong to the same irreducible component of $E$, in other words $\beta_t(j_1)=\beta_0(j_2)$.
The summand of the left hand side of Inequality~(\ref{ineqaux}) corresponding to $x$ is equal to 
$\mu_{t,j_1}$ and the summand of the right hand
side of Inequality~(\ref{ineqaux}) corresponding to $y$ is equal to $\mu_{0,j_2}$. 
Observe that $\mu_{t,j_1}$ is equal to $\mu_{\beta_t(j_1)}$ for $t\in D_\eta^*$ generic, and
that, in general, we have the inequality $\mu_{0,j_2}\geq\mu_{\beta_0(j_2)}=\mu_{\beta_t(j_1)}$.
Hence the contribution of $x$ to the LHS of Inequality~(\ref{ineqaux}) is smaller 
or equal than the contribution of the corresponding point $y$ to the RHS, and the contribution is strictly smaller if and only if 
$\mu_{0,j_2}<\mu_{\beta_0(j_2)}$. 
Therefore, if there exists $j\in\{1,...,r_0\}$ for which the inequality $\mu_{0,j}\leq\mu_{\beta_0(j)}$
is strict, then the Inequality~(\ref{ineqaux}) is strict and we get a contradiction again.
\end{proof}

\begin{prop}
\label{smooth}
There exist positive $\epsilon$ and $\eta$ such that
\begin{enumerate}
\item the ramification locus of $\alpha|_E$ is $\CC\times\{O\}$,
\item for any $t\in D_\eta$ and any $x\in E_{t,j}\setminus\{O\}$ 
we have $\mu(f_t|_{V(z_1-z_1(x))},x)=\mu_{\beta_t(j)}$ if $t\neq 0$ and $\mu(f_t|_{V(z_1-z_1(x))},x)=\mu_{\beta_{0}(j)}$ if $t=0$.
\end{enumerate}
\end{prop}
\begin{proof}
We show that if any of the conclusions is false then the L\^e number $\lambda^1_{f,\Zc}(t,O)$ is not 
independent on $t$, which is a contradiction, as we have seen in the proof of Lemma~\ref{constconst}.

Since the first polar variety of $f$ is empty, the first L\^e cycle $\Lambda^1_{f,\Zc}$ is the cycle associated to the scheme 
\[W:=\overline{V(\partial f/\partial z_2,...,\partial f/\partial z_n)\setminus \Sigma}\cap V(\partial f/\partial z_1).\]
The first L\^e number of $f$ at $(t_0,O)$ is the intersection multiplicity
\[\lambda^1_{f,\Zc}(t_0,O)=I_{(t_0,O)}([W],[V(t-t_0)]).\]
By conservation of number in intersection theory, if $t_0$ is small enough, we have
\[I_{(0,O)}([W],[V(t-t_0)])=\sum_{(t_0,x)\in W\cap V(t-t_0)} I_{(t_0,x)}([W],[V(t-t_0)]).\]
Hence, if there is an irreducible component of $W$ different from $\CC\times\{O\}$ and passing through $(0,O)$ then
$\lambda^1_{f,\Zc}(0,O)$ is strictly 
bigger than $\lambda^1_{f,\Zc}(t,O)$ for $t\neq 0$.

Suppose that we have $(t,x)\in E$ belonging to the ramification locus of $\alpha|_E$. We may assume 
(if $\eta$ is small enough) that in the discriminant of $\alpha$ there are no components of the form 
$\CC\times\{t\}$ for $t\neq 0$. Hence for any $s$ close enough to $z_1(x)$ the singularity that $f_t|_{V(z_1-z_1(x))}$ has at $x$ 
splits in several critical points of $f_t|_{V(z_1-s)}$. By the L\^e-Lazzery non-splitting result we deduce that there is at least one of them, called $x_s$, such that $f_t(x_s)$ is not zero. As $\Sigma$ is contained in $f^{-1}(0)$ we have that $(t,x_s)$ belongs to
$V(\partial f/\partial z_2,...,\partial f/\partial z_n)\setminus \Sigma$. As $x_s$ converges to $x$ as $s$ converges to $z_1(x)$ we conclude that $(t,x)$ belongs to $\overline{V(\partial f/\partial z_2,...,\partial f/\partial z_n)\setminus \Sigma}$. As 
$\partial f/\partial z_1(t,x)=0$ (for being $(t,x)\in\Sigma$) we deduce that $(t,x)$ belongs to $W$. Therefore, if the first assertion
is false there is an irreducible component of $W$ which is different from $\CC\times\{O\}$ and we get a contradiction.

If $t=0$, Assertion (2) follows easily (by shrinking $\epsilon$) from the second assertion of Proposition~(\ref{reduced}).
Suppose now that we have $(t,x)\in E$, with $t\neq 0$,  such that $x$ belongs to $E_{t,j}$ and 
\begin{equation}
\label{aqui2}
\mu(f_t|_{V(z_1-z_1(x))},x)>\mu_{\beta_t(j)}.
\end{equation}
Choose a neighborhood $U$ of $x$ in $\CC^n$ such that $x$ is the only critical point of $f_t|_{V(z_1-z_1(x))}$ in 
$U\cap V(z_1-z_1(x))$.
As, for {\em generic} $s$, the Milnor number of any of the singularities of $f_t|_{V(z_1-s)}$ at $\Sigma_{t,j}$ is equal to 
$\mu_{\beta_t(j)}$ (use previous Proposition applied to $E_t$ instead of $E_0$), by the strict inequality~(\ref{aqui2}), 
we can approximate $z_1(x)$ by a sequence $\{s_n\}_{n\in\NN}$ such that $f_t|_{V(z_1-s_n)}$ 
has several critical points in $U$. By the L\^e-Lazzery non-splitting result we deduce that there is at least 
one of them, called $x_n$, such that $f_t(x_n)$ is not zero.
From this point the proof proceeds like the proof for the first assertion.
\end{proof}

\begin{proof}[Proof of Theorem~\ref{leconsteqsing}]
By the remarks after the statement of the theorem we know that it is sufficient to prove that the first 
condition implies the conclussion.

The dimensional requirement of the first condition for equisingularity at the critical set follows from 
Lemma~\ref{thom} , Lemma~\ref{firstle} and Lemma~\ref{constconst}.
The first assertion of Proposition~\ref{smooth} implies that $E_t$ is smooth outside the origin, and that
\[\pi|_{E\setminus (\CC\times\{O\})}:E\setminus (\CC\times\{O\})\to D_\eta\] 
is a submersion. This gives the first condition for equisingularity at the critical set.

The first assertion of Proposition~\ref{smooth} also implies that for any $t\in D_\eta$ the mapping $z_1:E_t\setminus\{O\}\to\CC$ has 
no critical points. Therefore the mapping $|z_1|^2:E_t\to\RR$ has the origin as its only minimum and has no critical points 
in $E_t\setminus\{O\}$. Thus the restriction $|z_1|^2:E_t\setminus\{O\}\to (0,\epsilon]$ is a proper function without 
critical points for any $t\in D_\eta$. This function trivialises the desired cobordism and implies the second condition for
equisingularity at the critical set.

The third condition of equisingularity at the critical set follows from the second assertion of Proposition~\ref{smooth}.
\end{proof}

\section{Topological stems}

Here we modify Pellikaan's inductive definition of stem of degree $d$ in order to define topological 
stems, and prove Theorem~D. We refer the reader to the introduction and to the papers cited there to
find motivation for the definition and applications of stems, and for an explanation of why our 
modification of the definition is a reasonable one.

Let $\mm$ denote the maximal ideal of $\calO_{(\CC^n,O)}$.

\begin{definition}
\label{topstem}
Let $d$ be a positive integer. We define topological stems of degree $d$ inductively as follows:
\begin{itemize}
\item A holomorphic function germ $f:(\CC^n,O)\to\CC$ is a {\em topological stem of degree} $1$ if there
exists a positive integer $N$ such that for any $g\in\mm^N$, and $t\in\CC$ sufficiently small, 
either $f+tg$ has an isolated singularity 
at the origin, or it is topologically $R$-equisingular to $f$.
\item A holomorphic function germ $f:(\CC^n,O)\to\CC$ is a {\em topological stem of degree} $d$ if there
exists a positive integer $N$ such that for any $g\in\mm^N$, and $t\in\CC$ sufficiently small, either 
$f+tg$ is a stem of degree strictly 
smaller than $d$, or it is topologically $R$-equisingular to $f$. 
\end{itemize}
\end{definition}

\begin{theo}
\label{charstem}
A holomorphic function germ $f:(\CC^n,O)\to\CC$ is a topological stem of positive finite degree if and 
only if its critical set is $1$-dimensional at the origin.
\end{theo}
\begin{proof}
Suppose that the critical set of $f$ is $1$-dimensional at the origin. We will show that $f$ is a 
topological stem 
of degree bounded by its first L\^e number with respect to a generic coordinate system. We work by 
induction on the first L\^e number. 
Let $\mathbf{Z}$ be generic coordinate system. Then $f|_{V(z_1)}$ has an isolated singularity at the 
origin, with Milnor number $\mu$. As $f|_{V(z_1)}$ is $\mu+1$-determined, for any $g\in\mm^{\mu+2}$ and
any $t\in\CC$ the restriction of the function $h_t:=f+tg$ to $V(z_1)$ has an isolated singularity at the
origin (with the same Milnor number). In this situation the L\^e numbers of $h_t$ at the origin with 
respect to $\mathbf{Z}$ are well defined for any $t\in\CC$.

Suppose that $\lambda^1_{h_t,\Zb}(O)$ is constant for $t$ sufficiently small. 
In view of~Definition~4.1~of~\cite{Ma3}, it is easy to show that any polar ratio of $h_t$ is bounded 
above by $M_t:=\lambda^0_{h_t,\Zb}(O)+1$. By the lexicographical upper semicontinuity of the L\^e numbers
we have that
$M:=M_0\geq M_t$ for $t$ small. The function 
$h_0+z_1^M$ has an isolated singularity at the 
origin (see~Lemma~4.3,~\cite{Ma3}). Denote by $\mu_0$ its Milnor number. Consider $N:=\mu_0+2$. 
Notice that, as $h_0+z_1^M$ is 
$N-1$-determined, if $g$ belongs to $\mm^N$ we have the equality of Milnor numbers
\begin{equation}
\label{eqmu}
\mu(h_0+z_1^M)=\mu(h_t+z_1^M)
\end{equation}
for any $t$. 
By the L\^e-Iomdine formula, if the first L\^e number of $h_t$ is
constant then the $0$-th L\^e number is constant as well, and applying Theorem~\ref{wanted} we find that
$h_t$ is topologically $R$-equisingular for $t$ small enough. 

Otherwise $\lambda^1_{h_t,\Zb}(O)$ ($t\neq 0$ sufficiently small) is 
strictly smaller than $\lambda^1_{h_0,\Zb}(O)$. In this case, by induction hypothesis, the germ 
$h_t$ is a topological 
stem of finite degree bounded by its generic first L\^e number at the origin, wich is in turn 
bounded by $\lambda^1_{h_t,\Zb}(O)$. We conclude that $f$ is a topological stem 
of degree bounded by $\lambda^1_{f,\Zb}(O)$.
\end{proof}

\section{The underlying deformation of the critical set}
\label{greuel}

Notice that the condition of equisingularity at the critical set imposes no condition at the origin. 
This observation is probably the explanation of the phenomena shown in this section. Below we will show
how to produce examples of families which are equisingular at the critical set, but such that the critical set experiences 
drastic changes from the analytic viewpoint. We also prove a new topological $R$-equisingularity condition for families for which 
the reduced critical set undergoes a flat $\mu$-constant deformation of reduced curves in the sense of~\cite{BuG}.

We have shown that, if in a family of germs with $1$-dimensional critical set the L\^e numbers with respect to a generic coordinate 
system are constant, then the family is topologically $R$-equisingular. However topological equisingularity does not imply 
the constancy of the L\^e numbers, even if the dimension of the critical set is $1$ 
(in particular the L\^e numbers are not topological invariants). In~\cite{BG} the following counterexample was
constructed:

\begin{example}
\label{example1}
Define germs $f,g_t:(\CC^3,0)\to\CC$ by
\[f(x,y,z):=x^{15}+y^{10}+z^6\]
\[g_t(x,y,z):=xy+tz.\]
\[F_t:=f^2-g_t^{12}=(f+g^6_t)(f-g^6_t).\]
The family $F_t$ has critical set of dimension $1$, it is topologically $R$-equisingular, but the L\^e numbers with respect to a 
generic coordinate system are not independent on $t$.
\end{example}

Let us mention that $(f,g_t):(\CC^3,0)\to\CC$ is an example due to Henry (appearing in~\cite{BuG}) of a family of i.c.i.s. 
with constant Milnor number and non-constant multiplicity. In~\cite{BG} we present a proof of the topological $R$-equisingularity
of this example based on the results of~\cite{DG}. An alternate proof of this can be obtained proving that $F_t$ is equisingular 
at the critical set. To do this we first show that, after fixing a certain radius $\epsilon>0$, the critical set of $F_t$ with 
reduced structure is given by $V(f,g_t)$ (this can be done in the same way that we treat the next 
Example~\ref{example2}). After we observe
that, as $V(f,g_t)$ defines a $\mu$-constant deformation of the {\em reduced} curve singularity $V(f,g_0)$, we have that
$V(f,g_t)$ is topologically equisingular (see~\cite{BuG}). 
This implies the first two conditions for equisingularity at the critical set.
The third condition, which deals with transversal Milnor numbers can be established by an easy inspection of the equation of $F_t$.

The following theorem shows that the previous example illustrates a general phenomenon.

\begin{theo}
\label{gbb}
Let $f:\CC\times(\CC^n,O)\to (\CC,0)$ be a family of holomorphic germs at the origin, holomorphically depending on a parameter $t$, 
and having $1$-dimensional critical set at the origin. Let $\Sigma:=V(\partial f/\partial z_1,...,\partial f/\partial z_n)$.
Suppose that all the irreducible components of $\Sigma$ at $(0,O)$ are of dimension $2$ and that the restriction 
\[\pi:\Sigma_{red}\to \CC\]
is a flat morphism with {\em reduced} fibres (where $\Sigma_{red}$ is $\Sigma$ with reduced structure and $\pi$ is the restriction 
of the projection of $\CC\times\CC^n$ to the first factor) such that the Milnor number at the origin $\mu((\Sigma_{red})_t,O)$
(in the sense of~\cite{BuG}) is independent of $t$. If, in a neighborhood of $(0,O)$ in $\CC\times\CC^n$,
the generic transversal Milnor number of $f_t$ at any point $(t,x)$ of
$\Sigma_{red}\setminus\CC\times\{O\}$ only depends on the irreducible component of $\Sigma_{red}$ to which $(t,x)$ belongs 
then the family $f$ is topologically $R$-equisingular.
\end{theo}
\begin{proof}
Using~\cite{BuG}, Section~5 one proves easily the first two conditions of equisingularity at the critical set. After this the third
condition follows now easily from the hypothesis. The result follows applying Theorem~\ref{porfin}.
\end{proof}

However, the condition given in the previous theorem is again not a characterisation of topological $R$-equisingularity. The point
is that the family of {\em reduced} critical sets does not undergo, in general, a flat deformation.
If a holomorphic family $f:\CC\times (\CC^n,O)\to (\CC,0)$ of function germs $f_t$ with $1$-dimensional critical set at the origin,  is equisingular with respect to the critical set then it is easy to see that the restriction $\pi:\Sigma_{red}\to\CC$ is flat at
$(0,O)$. However, what is not true in general is that the fibre $(\Sigma_{red})_0$ is reduced. In fact it may have embedded
components at the origin.

Consider the following deformation of a parametrisation in $\CC^3$ (with deformation paremeter $t$):
\begin{equation}
\label{parametric}
 x=s^3\quad\quad y=s^4\quad\quad z=ts.
\end{equation}
The following equations in $\CC\{t,x,y,z\}$ define the image $Z\subset\CC\times\CC^3$ of the family as a set:
\begin{equation}
\label{implicit}
ty-xz=0\quad\quad tx^3-y^2z=0\quad\quad y^3-x^4=0\quad\quad t^3x-z^3=0.
\end{equation}

\begin{example}
\label{example2}
The family 
\begin{equation}
\label{ejemplo}
f_t:=(ty-xz)^9+(tx^3-y^2z)^4+(y^3-x^4)^3+(t^3x-z^3)^{12}
\end{equation}
is equisingular at the critical set. For any $t$ its critical set is the image of the parametrisation~(\ref{parametric}). 
Hence, the family of {\em reduced} critical sets does not undergo, in general, a flat deformation (the fibre of $\pi:\Sigma_{res}\to\CC$ at $0$ has an embedded component at the origin). Observe that the critical set $(\Sigma_t)_{red}$ is 
smooth for $t\neq 0$ and singular (of multiplicity $3$) for $t=0$. The transversal Milnor number is $6$ and, 
hence the first L\^e number with respect to a generic coordinate system is $6$ for $t\neq 0$ and $18$ for $t=0$.
Clearly any stabilisation of the form 
\begin{equation}
\label{ejemplo2}
g_t:=f_t+u^a+v^b
\end{equation}
has all the above properties and, by Theorem~\ref{porfin} is topologically $R$-equisingular.
\end{example}

 Before checking the assertions of the example, we observe the following:

 \begin{remark}
 \label{singsmooth}
 We have obtained two topologically $R$-equivalent 
 functions with critical locus of dimension $1$, such that the critical set is smooth for one of them and singular at the origin 
 for the other. To the author's knowledge it is the first time that this behavior is observed. 
 \end{remark}
  
 The equation $f_t$ is quasihomogeneous of degree $30$ if we give weights $(3,4,1)$ to the variables $(x,y,z)$. This shows that 
 for any $t$ the critical set of $f_t$ is contained in the central fibre. In the weighted homogeneous plane $P(3,4,1)$ it is
 easy to see that the curve defined by $f_t=0$ has a unique singular point corresponding to the curve $C_t$ parametrised by 
 equations~(\ref{parametric}). Hence the critical set of $f_t$ is precisely the curve $C_t$.
 To show the equisingularity at the critical set of the family $f_t$ we only have to check the condition on the transversal Milnor
 number, but this follows from an easy inspection of the equations (the reader may check that the generic transversal Milnor number at
 any point of $C_t\setminus\{O\}$ is 
 controlled by the following terms of the equation: $(tx^3-y^2z)^4+(y^3-x^4)^3$). 
 After this the equisingularity at the critical set of $g_t$ follows easily. As $g_t:\CC^5\to\CC$ is defined in the good range of
 dimensions ($n\geq 5$) we can apply Theorem~\ref{porfin} and conclude that the family $g_t$ is topologically $R$-equisingular.
 
 The fact that the transversal Milnor number in $C_t\setminus\{O\}$ is controlled by the terms $(tx^3-y^2z)^4+(y^3-x^4)^3$ suggests
 to play the following game. Observe that the multiplicity of $f_t$ is $9$ and that the monomials of order $9$ appear in the terms 
 $(y^3-x^4)^3$ (for any $t$) and  $(ty-xz)^9$ (only for $t\neq 0$). Replacing the term $(ty-xz)^9$ by $(ty-xz)^8$ we expect not to
 change the generic transversal Milnor number at any point of $C_t\setminus\{O\}$, since it is controlled by two different terms.
 In fact one may check that 
 \begin{equation}
\label{ejemplo3}
h_t:=(ty-xz)^8+(tx^3-y^2z)^4+(y^3-x^4)^3+(t^3x-z^3)^{12}
\end{equation}
 is a family such that the only component of the critical set of $f_t$ meeting the origin is precisely $C_t$, 
with generic transversal Milnor
 number equal to $6$ at any point of $C_t\setminus\{O\}$. Unfortunately the family is not equisingular at the origin since for 
 $t\neq 0$ the critical locus contains precisely $1$ Morse point outside $C_t$, and this Morse point converges to the origin as 
 $t$ approaches $0$. As the family has not constant multiplicity, if this Morse point would have not appeared, then 
 we would have obtained a counterexample of Zariski's multiplicity conjecture (after adding high powers of
 new variables in order to meet the dimensional restrictions of our equisingularity results).
 
 The last examples shows that we have a great deal of flexibility in deforming the critical set in a family which is equisingular
 at the critical set (and, in the correct range of dimensions, topologically $R$-equisingular). We propose the following problem:
 
 \begin{problem}
 \label{problem}
 Construct examples of non-equimultiple, but topologically equisingular,  
 families of parametrised curves (like~\ref{parametric}), with one or several components, and families of
 functions whose critical set consists exactly with the parametrised curve and has a prescribed transversal Milnor number outside 
 the origin (like $f_t$). Observe that in these conditions the constructed family is equisingular at the critical set. 
 \end{problem}
 
 Specially in the case in which $C_t$ is not an i.c.i.s., and thus we need more equations than its codimension to define it, 
 we could have enough space to construct examples of equisingular at the critical set families which are not equimultiple.

\end{document}